-------------------------

%

%
\input amstex
\documentstyle{amsppt}
\loadbold
\def\cstar{$C^*$-algebra}

\def\<{\left<}										
\def\>{\right>}

\def\tr{\text{trace}\,}

\def\rank{\text{rank}}
\def\alg{\text{alg}}
\def\sp{\text{span}}

\def\qedd{
\hfill
\vrule height4pt width3pt depth2pt
\vskip .5cm
}

\magnification=\magstephalf

\topmatter
\title
The Curvature Invariant of a Hilbert 
Module over $\Bbb C[z_1,\dots,z_d]$
\endtitle

\author William Arveson
\endauthor

\affil Department of Mathematics\\
University of California\\Berkeley CA 94720, USA
\endaffil

\date 4 August 1998
\enddate
\thanks This research was supported by
NSF grants DMS-9500291 and DMS-9802474
\endthanks
%
%
\abstract 
A notion of curvature is introduced in 
multivariable operator theory and an analogue
of the Gauss-Bonnet-Chern theorem is established.
Applications are given to the metric structure
of graded ideals in $\Bbb C[z_1,\dots,z_d]$, 
and the existence of ``inner" sequences for 
closed submodules 
of the free Hilbert module $H^2(\Bbb C^d)$.  
\endabstract


\toc
\specialhead{Introduction} 
\endspecialhead

\subhead 1.  Free Hilbert modules and dilation theory
\endsubhead
\subhead 2.  Multipliers of free Hilbert modules
\endsubhead
\subhead 3.  Euler characteristic
\endsubhead
\subhead 4.  Curvature invariant
\endsubhead
\subhead 5.  Curvature operator: quantizing the Gauss map
\endsubhead
\subhead 6.  Graded Hilbert modules
\endsubhead
\subhead 7.  Degree
\endsubhead
\subhead 8.  Applications: inner sequences
and graded ideals
\endsubhead
\subhead 9.  Examples
\endsubhead

\specialhead{References} 
\endspecialhead
\endtoc

\endtopmatter
%

\document

\subheading{Introduction}

Let $\bar T=(T_1,\dots,T_d)$ be a $d$-tuple of mutually 
commuting operators acting on a common Hilbert space 
$H$.  $\bar T$ is called a $d$-{\it contraction} if 
$$
\|T_1\xi_1+\dots+T_d\xi_d\|^2\leq \|\xi_1\|^2+\dots+\|\xi_d\|^2
$$
for all $\xi_1,\dots,\xi_d\in H$.  The number $d$ will be fixed 
throughout this paper, and of course we are primarily 
interested in the cases $d\geq 2$.  Let $A=\Bbb C[z_1,\dots,z_d]$
be the complex unital algebra of all polynomials in $d$ 
commuting variables $z_1,\dots,z_d$.  A commuting 
$d$-tuple $T_1,\dots,T_d$ of operators in the 
algebra  $\Cal B(H)$ of all bounded operators on $H$ gives rise 
to an $A$-module structure on $H$ in the natural way,
$$
f\cdot \xi = f(T_1,\dots,T_d)\xi,\qquad f\in A,\quad \xi\in H; 
$$
and $(T_1,\dots,T_d)$ is a $d$-contraction iff $H$ is a 
{\it contractive} $A$-module in the following sense,
$$
\|z_1\xi_1+\dots+z_d\xi_d\|^2\leq \|\xi_1\|^2+\dots+\|\xi_d\|^2
$$
for all $\xi_1,\dots,\xi_d\in H$.  Thus it is equivalent to speak
of $d$-contractions or of contractive Hilbert $A$-modules, and we
will shift from one point of view to the other when it is convenient
to do so.  

For every $d$-contraction $\bar T=(T_1,\dots,T_d)$ we have 
$0\leq T_1T_1^*+\dots+T_dT_d^*\leq \bold 1$, and hence 
$$
\Delta = (\bold 1 - T_1T_1^*-\dots-T_dT_d^*)^{1/2}\tag{0.1}
$$
is a positive operator on $H$ of norm at most one.  
The {\it rank} of $\bar T$ is defined as the dimension of 
the range of $\Delta$.  Throughout this paper we will be 
primarily concerned with finite rank $d$-contractions (resp. 
finite rank contractive Hilbert $A$-modules).   

We introduce several numerical invariants for finite rank 
contractive $A$-modules $H$, the principal ones being
the curvature invariant $K(H)$, 
the Euler characteristic $\chi(H)$, and the 
degree $\deg(H)$.  All of these quantities
are real numbers (indeed, most are integers), and we develop their 
basic properties.  

We now describe the main results of this paper, starting with 
a sketch of the 
definition of the curvature invariant $K(H)$.  Let $H$ be a finite
rank contractive Hilbert 
$A$-module with associated $d$-contraction  
$(T_1,\dots, T_d)$.  For every point $z=(z_1,\dots,z_d)$ in 
complex $d$-space $\Bbb C^d$ we form the operator 
$$
T(z) = \bar z_1T_d+\dots+\bar z_d T_d\in\Cal B(H),  
\tag{0.2}
$$
$\bar z_k$ denoting the complex conjugate of the complex
number $z_k$.  Notice that the operator function 
$z\mapsto T(z)$ defines
an antilinear mapping of $\Bbb C^d$ into $\Cal B(H)$, 
and since $(T_1,\dots,T_d)$ is a $d$-contraction we have 
$$
\|T(z)\|\leq |z|=(|z_1|^2+\dots+|z_d|^2)^{1/2} 
$$
for all $z\in \Bbb C^d$.  In particular, if $z$ belongs to 
the open unit ball 
$$
B_d=\{z\in\Bbb C^d: |z|<1\}
$$
then 
$\|T(z)\|<1$ and hence $\bold 1-T(z)$ is invertible.  Thus
for every $z\in B_d$ we can define a positive operator 
$F(z)$ acting on the finite dimensional Hilbert space 
$\Delta H$ as follows,
$$
F(z)\xi = \Delta (\bold 1-T(z)^*)^{-1}(\bold 1-T(z))^{-1}\Delta\xi,
\qquad \xi\in\Delta H.  
$$

We require the boundary values of the real-valued
function $z\in B_d\mapsto\tr F(z)$, which 
exist in the following sense.  Let 
$\partial B_d = \{z\in \Bbb C^d: |z|=1\}$ be the unit sphere 
in $\Bbb C^d$ and let $\sigma$ be normalized surface measure
on $\partial B_d$.  

\proclaim{Theorem A}
For $\sigma$-almost every $\zeta\in \partial B_d$, the
limit 
$$
K_0(\zeta) = \lim_{r\uparrow 1}(1-r^2)\tr F(r\zeta)
$$
exists and satisfies $0\leq K_0(\zeta)\leq \rank(H)$.  
\endproclaim

Theorem A is proved in section 4.  We define the 
curvature invariant by integrating $K_0$ over the 
sphere
$$
K(H) = \int_{\partial B_d} K_0(\zeta)\,d\sigma(\zeta).  
\tag{0.3}
$$
Obviously, $K(H)$ is a real number satisfying 
$0\leq K(H)\leq {\text{rank}}(H)$.  

We now define the Euler characteristic $\chi(H)$ of a finite rank 
contractive $A$-module $H$.  $\chi(H)$ depends only on the 
{\it algebraic} structure of the following $A$-submodule of $H$:
$$
M_H = {\text{span}}\{f\cdot \xi: f\in A, \xi\in \Delta H\}.  
$$
Notice that we have not taken the closure in forming $M_H$.  
Note too that if $r={\text{rank}}(H)$ and $\zeta_1,\dots,\zeta_r$ 
is a linear basis for $\Delta H$, then $M_H$ is the set of 
``linear combinations"
$$
M_H=\{f_1\cdot\zeta_1+\dots+f_r\cdot \zeta_r: f_k\in A\}.  
$$
In particular, $M_H$ is a finitely generated $A$-module.  

It is a consequence of Hilbert's syzygy theorem 
for ungraded modules (cf. Theorem 182 of
\cite{18} or Corollary 19.8 of 
\cite{14}) that $M_H$ has a finite free resolution; that 
is, there is an exact sequence of $A$-modules 
$$
0\rightarrow F_n\rightarrow \dots \rightarrow F_2
\rightarrow F_1\rightarrow M_H\rightarrow 0
\tag{0.4}
$$
where $F_k$ is a free module of finite rank $\beta_k$,
$$
F_k = \underbrace{A\oplus \dots\oplus A}_{\beta_k {\text{times}}}.  
$$
The alternating sum of the ``Betti numbers" of this 
free resolution 
$$
\sum_{k=1}^n (-1)^{k+1}\beta_k 
$$
does not depend on the particular finite free 
resolution of $M_H$, and hence we may define 
the {\it Euler characteristic} of $H$ by 
$$
\chi(H) = \sum_{k=1}^n(-1)^{k+1}\beta_k,\tag{0.5}
$$
where $\beta_k$ is the rank of $F_k$ in any 
finite free resolution of $M_H$ of 
the form (0.4).  

One of the more notable results in the 
Riemannian geometry of surfaces is the Gauss-Bonnet 
theorem, which asserts that if $M$ is a compact 
oriented Riemannian $2$-manifold and 
$$
K:M\to \Bbb R
$$
is its Gaussian curvature function, then 
$$
\frac{1}{2\pi}\int_MK\, dA = \beta_0-\beta_1+\beta_2
\tag{0.6}
$$
where $\beta_k$ is the $k$th Betti number of $M$.  
In particular, the integral of $K$ depends only on 
the {\it topological} type of $M$.  This remarkable 
theorem was generalized by Shiing-Shen Chern to compact
oriented even-dimensional Riemannian manifolds in 
1944 \cite{6}.  

We will establish the following 
result in section 6, which we view as an analogue of 
the Gauss-Bonnet-Chern theorem for graded 
Hilbert $A$-modules.  By a {\it graded} Hilbert $A$-module
we mean a pair $(H,\Gamma)$ where $H$ is a 
(finite rank, contractive) Hilbert $A$-module and 
$\Gamma: \Bbb T\to \Cal B(H)$ is a strongly continuous 
unitary representation of the circle group such that 
$$
\Gamma(\lambda)T_k\Gamma(\lambda)^{-1} = \lambda T_k,
\qquad k=1,2,\dots,d, \lambda\in \Bbb T,
$$
$T_1,\dots,T_d$ being the $d$-contraction associated with 
the module structure of $H$.  Thus, graded Hilbert 
$A$-modules are precisely those whose underlying 
operator $d$-tuple 
$(T_1,\dots,T_d)$ possesses circular symmetry. 
$\Gamma$ is called the gauge group of $H$.  

\proclaim{Theorem B}
Let $H$ be a graded (contractive, finite rank) Hilbert 
$A$-module for which the spectrum of the gauge 
group is bounded below.  Then $K(H) = \chi(H)$, and
in particular $K(H)$ is an integer.  
\endproclaim 

We remark that the hypothesis on 
the spectrum of the gauge group
is equivalent to several other natural ones, see
Proposition 6.4.  Theorem B depends on the following 
asymptotic formulas for $K(H)$ and $\chi(H)$, which 
are valid for finite rank contractive Hilbert $A$-modules, 
graded or not.  For such an $H$, let $(T_1,\dots, T_d)$
be its associated $d$-contraction and define a completely 
positive normal map $\phi: \Cal B(H)\to\Cal B(H)$ by 
$$
\phi(A) = T_1AT_1^*+\dots +T_dAT_d^*.  
$$
Since $H$ is contractive and finite rank,
$\bold 1-\phi(\bold 1)$  is a positive finite
rank operator, and a simple argument shows that
$\bold 1-\phi^n(\bold 1)$ is a positive finite
rank operator for every $n=1,2,\dots$.

\proclaim{Theorem C}
For every contractive finite rank Hilbert $A$-module 
$H$,
$$
\chi(H) = d! \lim_{n\to \infty}
\frac{{\text{rank}}\,\,(\bold 1-\phi^{n+1}(\bold 1))}{n^d}. 
$$
\endproclaim

\proclaim{Theorem D}
For every contractive finite rank Hilbert $A$-module 
$H$,
$$
K(H) = d! \lim_{n\to \infty}
\frac{{\text{trace}}\,\,(\bold 1-\phi^{n+1}(\bold 1))}{n^d}. 
$$
\endproclaim
Theorems C and D are proved in sections 3 and 5.  The number
$K(H)$ is actually the trace of a certain 
self-adjoint trace-class operator 
$d\Gamma$, which exists for any finite 
rank contractive Hilbert module.  
While the trace of this operator 
is therefore always nonnegative, it is noteworthy that 
$d\Gamma$ itself is never a positive operator.  Indeed, 
we have found it useful to think of $d\Gamma$ 
as a higher dimensional operator-theoretic counterpart of 
the differential of the Gauss 
map  $\gamma: M\to S^2$ of an oriented $2$-manifold 
$M\subseteq\Bbb R^3$ (cf. \cite{9}, pp 136--146).  
We have glossed
over some details in order to make the essential point; see
section 5 for a more comprehensive discussion.
In any case, the formula 
$$
K(H) = {\text{trace}}\,\,d\Gamma
$$
is an essential component in the proofs of Theorems 
B and D (see Theorem 5.13 et seq).  

These results have concrete implications 
about the invariant subspaces of $H^2$ and 
the algebraic structure of graded ideals 
in the polynomial algebra $\Bbb C[z_1,\dots,z_d]$. 
The applications are discussed in section 8, 
and are briefly summarized as follows.  

The dilation theory described in section 1 implies
that for every closed submodule $M$ of the free 
Hilbert module $H^2$, there is a (finite or infinite)
sequence $\Phi=\{\phi_1,\phi_2,\dots\}$ of multipliers 
of $H^2$ whose associated multiplication operators 
$M_{\phi_n}$ satisfy
$$
P_M = M_{\phi_1}M_{\phi_1}^*+ M_{\phi_2}M_{\phi_2}^*+\dots,
\tag{0.9}
$$
$P_M$ denoting the orthogonal projection of $H^2$ 
on $M$.  Every multiplier $\phi$ can be regarded 
as a bounded holomorphic function defined on the open
unit ball $B_d=\{z\in\Bbb C^d: |z|<1\}$ and has a 
radial limit function $\tilde\phi: \partial B_d\to\Bbb C$
defined almost everywhere $(d\sigma)$ by
$$
\tilde\phi(z)=\lim_{r\to 1}\phi(rz), \qquad z\in\partial B_d.
$$
Formula (0.9) implies that the boundary values satisfy 
$\sum_n|\tilde\phi_n(z)|^2\leq 1$ almost everywhere 
on $\partial B_d$, and $\Phi$ is called an 
{\it inner sequence} if we have equality 
$$
\sum_n|\tilde\phi_n(z)|^2=1
$$
almost everywhere on $\partial B_d$.  

In dimension $d=1$, a familiar theorem of Beurling 
implies that there is a single multiplier $\phi$ which 
satisfies $P_M=M_\phi M_\phi ^*$, and such a multiplier 
must be an inner function.  By contrast, in dimension 
$d\geq 2$, there may be no single multiplier $\phi$ 
satisfying $P_M=M_\phi M_\phi^*$; indeed, in most cases 
the sequences $\{\phi_1,\phi_2,\dots\}$ of 
(0.9) are necessarily infinite (see the corollary of 
Theorem F below).  Moreover, we do not know if these 
infinite sequences associated with $M$ are inner 
in general.  This problem is associated 
with the fact that in dimension $d\geq 2$ the
canonical operators $S_1,\dots,S_d$ 
associated with $H^2$ do not form a subnormal
$d$-tuple, and are not usefully 
related to the $L^2$ space of any measure. 
Making use of the curvature invariant and 
a theorem of Auslander and Buchsbaum on the 
vanishing of the Euler characteristic of finitely 
generated modules over polynomial rings, we 
establish the following result which appears to 
cover many cases of interest.   

\proclaim{Theorem E}
Let $M$ be a closed submodule of $H^2$ which contains 
a nonzero polynomial.  Then every sequence of multipliers
$\Phi$ satisfying (0.9) is an inner sequence.  
\endproclaim

Let $I$ be an ideal in $\Bbb C[z_1,\dots,z_d]$. 
Hilbert's basis theorem implies that there  
is a finite set of 
elements $\phi_1,\dots, \phi_r\in I$ which generates
$I$ in the sense that 
$$
I = \{f_1\cdot\phi_1+\dots+f_r\cdot \phi_r: f_k\in 
\Bbb C[z_1,\dots,z_d]\}.  
$$
If $I$ is {\it graded} in the sense that it is spanned 
by its homogeneous polynomials, then one can find a 
set $\phi_1,\dots,\phi_r$ of generators
such that 

\roster
\item"{(0.7.1)}"
each $\phi_k$ is a homogeneous polynomial of 
some degree $n_k$, and  
\item"{(0.7.2)}"
$\{\phi_1,\dots,\phi_r\}$ is linearly independent.  
\endroster
Of course, systems of generators satisfying the conditions
(0.7) are by no means unique.  

We want to relate sets of generators of 
graded ideals to the 
natural norm on $\Bbb C[z_1,\dots,z_d]$, 
obtained by restricting the Hilbert 
space norm on $H^2$ to polynomials.  
To do that effectively we must 
consider infinite generators.  Let $\Phi$
be a (perhaps infinite) linearly independent set of 
homogeneous polynomials in a graded ideal $I$ which 
is {\it contractive} in the sense that whenever 
$\phi_1,\dots,\phi_r$ are distinct elements of 
$\Phi$ and $f_1,\dots,f_r\in\Bbb C[z_1,\dots,z_d]$ 
we have
$$
\|f_1\cdot\phi_1+\dots+f_r\cdot\phi_r\|^2\leq 
\|f_1\|^2+\dots+\|f_r\|^2.  
$$
One can scale down any set $\Phi$ of polynomials
so as to achieve this condition.  A set $\Phi$ 
satisfying the above conditions is said 
to be a {\it metric basis} for $I$ if 
every polynomial $g$ of degree 
$n$ in $I$ can be represented as a sum
$g=f_1\cdot\phi_1+\dots+f_r\cdot\phi_r$ where
$f_1,\dots,f_r$ and $\phi_1,\dots,\phi_r$ are as 
above and, in addition, satisfy
\roster
\item"{(0.8.1)}"
$\deg f_k+\deg\phi_k\leq n$, \qquad $k=1,\dots,r$, and 
\item"{(0.8.2)}"
$\|g\|^2=\|f_1\|^2+\dots+\|f_r\|^2$.  
\endroster
Condition (0.8.1) controls the degrees
of $\phi_1,\dots,\phi_r$ and (0.8.2) asserts that
the norms of $f_1,\dots,f_r$ are as small as 
the contractive hypothesis allows.  

In section 8 we show that every graded ideal 
in $\Bbb C[z_1,\dots,z_d]$ has a metric basis,
that the elements of a metric basis are mutually
orthogonal, 
and that any two metric bases are equivalent in a 
natural sense.  Thus, by giving up the 
requirement of finite generation of ideals, one obtains 
a uniqueness result for infinite generators satisfying 
(0.8.1) and (0.8.2).

Let $\Phi=\{\phi_1, \phi_2,\dots\}$ be a metric basis
for a graded ideal $I$ in the polynomial 
algebra $\Bbb C[z_1,\dots,z_d]$, and 
let $\sigma$ be the natural measure on 
the unit sphere $\partial B_d\subseteq \Bbb C^d$.  
Theorem E implies that 
for $\sigma$-almost every point $\zeta\in\partial B_d$
we have 
$$
\sum_n |\phi_n(\zeta)|^2 = 1.  
$$
One cannot expect the preceding ``almost 
everywhere" equation to hold everywhere on the unit 
sphere.  Indeed, if 
$$
V = \{z\in\Bbb C^d: \phi_1(z)=\phi_2(z)=\dots=0\}
$$
is the variety of common zeros of the polynomials 
$\phi_k$ (i.e., the zero set of the 
ideal $I$) then, since each $\phi_k$ is a homogeneous 
polynomial, $V$ is invariant under multiplication 
by positive scalars.  So whenever $V$ contains something 
other than $(0,0,\dots,0)$ it must intersect the unit 
sphere in $\Bbb C^d$, and in that case 
$V\cap\partial B_d$ is a nonvoid compact set of 
measure zero on which the $\phi_k$ all vanish.  

The following result implies that 
a {\it finite} generating set for
a graded ideal cannot be a metric basis except in 
a few insignificant cases.  

\proclaim{Theorem F}
Let $I$ be a graded ideal in $\Bbb C[z_1,\dots,z_d]$ 
whose metric basis is a finite set 
$\{\phi_1,\dots, \phi_n\}$.  Then $I$ is of finite 
codimension in $\Bbb C[z_1,\dots,z_d]$ and each of 
the canonical coordinates $z_1,\dots,z_d$ is 
nilpotent modulo $I$.  
\endproclaim

Given any finite rank (contractive) Hilbert module 
H over $\Bbb C[z_1,\dots,z_d]$ and a closed submodule
$K\subseteq H$, then both $K$ and its quotient 
$H/K$ are contractive Hilbert modules.  It is 
quite easy to see that 
${\text{rank }}H/K\leq {\text{rank }}H$, and 
hence $H/K$ is also of finite rank.  
However one does not have control
over the rank of the submodule
$K$.  Indeed, Theorem F has the following 
consequence, which stands in rather
stark contrast with the
assertion of Hilbert's basis theorem.  

\proclaim{Corollary}
Let $K$ be a (nonzero, closed) graded submodule 
of the free Hilbert module $H^2$ which is of 
infinite codimension in $H^2$.  Then
${\text{rank }}K = \infty$.  
\endproclaim

We end the paper with a 
discussion of some examples 
that serve to illustrate the properties 
of the invariants described above.  In 
particular, we show that any variety in 
complex projective space $\Bbb P^{d-1}$ gives
rise to a pure graded rank-one Hilbert 
$\Bbb C[z_1,\dots,z_d]$-module, and 
for some of these examples we compute 
all invariants in explicit 
terms (see section 9).

This work 
was initiated in order to obtain numerical invariants
for normal completely positive maps of $\Cal B(H)$.  
Notice that all of the invariants
introduced in this paper depend only on the properties
of the map $\phi: \Cal B(H)\to\Cal B(H)$
$$
\phi(A) = T_1AT_1^*+\dots+T_dAT_d^* \tag{0.10}
$$
associated with the canonical operators $T_1,\dots,T_d$ 
of the $\Bbb C[z_1,\dots,z_d]$-module structure.  
Indeed, Theorems C and D express $\chi(H)$ and 
$K(H)$ explicitly in terms of $\phi$, and 
Proposition 7.5 does the same for the secondary 
invariants $\deg(H)$, $\mu(H)$.  Thus these 
numbers are actually invariants of certain completely positive 
maps $\phi: \Cal B(H)\to \Cal B(H)$.  To be sure, not every 
weakly continuous completely positive map of $\Cal B(H)$ has the 
form (0.10) (with commuting operators $T_k$).   
Thus there remains an important question 
concerning the extent to which 
these results can be generalized to the 
case of noncommuting $d$-tuples of operators and their 
completely positive maps.  Ultimately, there is 
an associated problem of finding 
new numerical invariants for 
noncommutative dynamics, that is, for semigroups
of completely positive maps of $\Cal B(H)$ and
their relatives, $E_0$-semigroups.  Until now, we have 
had only the index \cite{2} and the geometric 
structures constructed in \cite{3}.

\subheading{1. Free Hilbert modules and dilation theory}

The algebra of polynomials $\Bbb C[z_1,\dots,z_d]$ in 
$d$ commuting variables (which we will abbreviate by
$A$ whenever it does not lead to 
confusion) has a natural inner product which can be 
defined in several ways.  Here, we define this inner 
product in terms of the relation that exists between
$\Bbb C[z_1,\dots,z_d]$ and the symmetric Fock space.  
Let $T_+(E)$ be the symmetric tensor algebra
of a $d$-dimensional complex vector space
$E$.  Writing $E^n$ for the symmetric tensor product
of $n$ copies of $E$ for $n\geq 1$ (with $E^0$ defined 
as $\Bbb C$) one finds that the algebraic direct sum
of vector spaces 
$$
T_+(E) = E^0\dotplus E^1\dotplus E^2\dotplus\dots
$$
is a commutative algebra with unit with respect to the 
multiplication defined by symmetric tensoring.  It
has the following universal property:
any linear mapping $L:E\to B$ of $E$ 
into a unital commutative
complex algebra $B$ extends uniquely to a unital 
homomorphism of complex algebras $\tilde L:T_+(E)\to B$.  
In particular, if we choose a linear basis $e_1,\dots,e_d$
for $E$ then there is a unique homomorphism of unital 
algebras
$\alpha: T_+(E)\to\Bbb C[z_1,\dots,z_d]$ defined by 
$\alpha(e_k)=z_k$ for $k=1,2,\dots,d$, and of course
in this case $\alpha$ is an isomorphism which we can 
use to identify $T_+(E)$ with $\Bbb C[z_1,\dots,z_d]$
if we wish.   

If we now fix an inner product on the 
one-particle space $E$ then $E$ becomes
a finite dimensional Hilbert space, and so does the tensor
product $E^{\otimes n}$ of $n$ copies of $E$ for every 
$n=2,3,\dots$.  Since the symmetric space $E^n$ is a 
subspace of $E^{\otimes n}$ for every $n\geq 1$ it follows
that $E^n$ is naturally a Hilbert space; and for $n=0$ 
we take the usual inner product on $E^0=\Bbb C$, 
$\<z,w\>=z\bar w$.  Thus, the algebraic direct sum 
$$
T_+(E) = \Bbb C\oplus E^1\oplus E^2\oplus \dots
$$
becomes an inner product space, which is 
dense in the symmetric Fock space over 
the Hilbert space $E$.  

Now we transport the inner product on $T_+(E)$ to an 
inner product on polynomials by picking an {\it orthonormal}
basis $e_1,\dots,e_d$ for $E$, and identifying $T_+(E)$ 
with $\Bbb C[z_1,\dots,z_d]$ by identifying $e_k$ with 
$z_k$ as above.  The completion of 
the polynomials in this inner product is a 
Hilbert space we denote by $H^2(\Bbb C^d)$ or simply
$H^2$ when, as will normally be the case, the 
dimension $d$ is fixed.  

The elements of $H^2$ can be realized as 
{\it certain} holomorphic 
functions in the open unit ball 
$$
B_d=\{z=(z_1,\dots,z_d)\in\Bbb C^d:
|z|=(|z_1|^2+\dots+|z_d|^2)^{1/2}<1\}
$$
which satisfy the following growth condition near the boundary
$$
|f(z)|\leq \frac{\|f\|}{\sqrt{1-|z|^2}}, \qquad z\in B_d.  
$$
We refer the reader to part I of 
\cite{1} for other characterizations
of this Hilbert norm on polynomials and 
for a development of the 
function-theoretic properties of the space $H^2$.  
Here, we summarize a few of its basic features.  For 
every $\alpha\in B_d$ the function $u_\alpha:B_d\to\Bbb C$
defined by
$$
u_\alpha(z) =\frac{1}{1-\<z,\alpha\>},\qquad |z|<1
$$
belongs to 
$H^2$; $H^2$ is spanned by 
$\{u_\alpha: \alpha\in B_d\}$ 
and for every $\alpha,\beta\in B_d$
we have 
$$
\align
\<u_\alpha,u_\beta\> &=\frac{1}{1-\<\beta,\alpha\>}, \\
\<f,u_\alpha\>&=f(\alpha),\qquad f\in H^2.  
\tag{1.1.2}
\endalign
$$

We also have 
$$
\|z_1f_1+\dots+z_df_d\|^2\leq \|f_1\|^2+\dots+\|f_d\|^2
$$
for every $f_1,\dots,f_d\in H^2$, so that in fact $H^2$ is 
a contractive Hilbert $A$-module.  The $d$-tuple of operators
$S_1,\dots,S_d$ obtained by multiplying by the $d$ coordinate 
functions define a $d$-contraction on $H^2$.  This 
$d$-contraction is called the $d$-shift in \cite{1}, and it 
has the property
$$
S_1S_1^*+\dots+S_dS_d^*=\bold 1 - E_0
\tag{1.2}
$$
where $E_0$ denotes the projection of $H^2$ onto the 
one-dimensional space of constant functions.  Using the 
terminology introduced in the introduction, $H^2$ is 
a contractive Hilbert $A$-module of rank one 
(the rank of Hilbert modules is defined later
in this section).  The Hilbert $A$-module $H^2$ will occupy 
a central position throughout this paper.  

We work with the category 
whose objects are {\it contractive} Hilbert modules over the 
algebra $A=\Bbb C[z_d,\dots,z_d]$ of polynomials, and 
we will refer to such modules simply as {\it Hilbert 
$A$-modules}.  Given two such modules $H_1$, $H_2$, 
$\hom(H_1,H_2)$ will denote the convex set of all 
operators $A\in\Cal B(H_1,H_2)$ which are contractions
($\|A\|\leq 1$) and which intertwine the 
respective module actions, 
$A(f\cdot\xi)=f\cdot A\xi$, $\xi\in H_1$.  Notice that
an isomorphism in $\hom(H_1,H_2)$ is necessarily a unitary 
operator that intertwines the module actions, and 
when such an operator exists we say 
that $H_1$ and $H_2$ are {\it isomorphic} and write
$H_1\cong H_2$.  

There are natural notions of (closed) submodule 
and quotient module in this category.  A closed 
submodule $K$ of a 
contractive Hilbert $A$-module $H$ is 
a contractive Hilbert module.  The quotient $H/K$ is 
of course a Banach space whose norm arises from an 
inner product on $H/K$; thus $H/K$ is also a contractive
Hilbert $A$-module.  

In more explicit operator-theoretic terms, let $T_1,\dots,T_d$
be the canonical operators associated with the Hilbert 
$A$-module $H$ defined by $T_k\xi=z_k\xi$, $\xi\in H$.  
Given a closed submodule $K\subseteq H$, the operators 
$T_1,\dots, T_d$ can be compressed to the coinvariant subspace
$K^\perp\subset H$ to obtain a $d$-contraction 
$\dot T_1,\dots \dot T_d$ which acts on  $K^\perp$ 
as follows
$$
\dot T_k\eta = P_{K^\perp}T_d\eta,\qquad \eta\in K^\perp,  
$$
$P_{K^\perp}$ denoting the orthogonal projection of $H$ 
on $K^\perp$.  
One finds that the Hilbert $A$-module structure on 
$K^\perp$ defined by $\dot T_1,\dots,\dot T_d$ is isomorphic
to the Hilbert $A$ module structure of the quotient $H/K$.  

The Hilbert module point of view has been emphasized 
by Douglas and Paulsen in their work on representations
of function algebras \cite{12}.  Significantly, 
the most important 
$d$-contractions cannot be dilated to {\it normal} 
$d$-contractions because of the failure of the 
von Neumann inequality for the unit ball in $\Bbb C^d$
(see \cite{1, theorem 3.3}).  It follows that
contractive Hilbert $A$-modules
cannot be profitably related 
to function algebras, and one must give up the 
idea of working with normal dilations.  
Instead, one seeks to relate Hilbert $A$ modules
to $H^2$ and multiples of $H^2$.  This 
dilation theory was worked out in \cite{1}, 
and is very effective for the category of 
(contractive) Hilbert $A$-modules.  
The purpose of this section 
is to reformulate those operator-theoretic 
results so that they are closer to the 
homological spirit of the central issues of this paper.  

Suppose we are given a submodule $K\subseteq M$
of a Hilbert $A$-module $M$.  The Hilbert modules $H$
which are isomorphic to the quotient $M/K$ are 
precisely those for 
which there is an exact sequence of Hilbert $A$-modules
$$
0\longrightarrow K\longrightarrow M
\underset U\to\longrightarrow 
H\longrightarrow 0
$$
where the connecting map $U$ is a coisometry, 
i.e., $UU^*=\bold 1_H$.  This leads us to an
important notion.  

\proclaim{Definition 1.4}
Let $H$ be a Hilbert $A$-module.  A dilation of 
$H$ is an exact sequence of Hilbert $A$-modules 
$$
M\underset U\to\longrightarrow H\longrightarrow 0
$$
where $U$ is a coisometry.  
\endproclaim

Notice that the kernel of $U$ is left unspecified
in this definition.  Two dilations 
$M_k\underset U_k\to\longrightarrow H\longrightarrow 0$,
$k=1,2$ are said to be {\it equivalent} if there 
is an isomorphism $W: M_1\to M_2$ of Hilbert 
$A$-modules such that $U_2W=U_1$.  There is no uniqueness
of dilations in this generality.  Indeed, if 
$M\underset U\to\longrightarrow H\longrightarrow 0$ is 
any dilation of $H$ and $N$ is an arbitrary Hilbert 
$A$-module, then we can construct an 
essentially different dilation 
$$
M\oplus N\underset V\to\longrightarrow H\longrightarrow 0
\tag{1.5}
$$
by taking $V$ to be the unique operator from 
$M\oplus N$ to $H$ which 
restricts to $U$ on $M$ and to $0$
on $N$.  

Before discussing this phenomenon, we collect 
some terminology that will be used throughout
the sequel.  Let $H$ be  Hilbert $A$-module 
and let $T_1,\dots,T_d$ be the canonical operators
associated with the module structure of $H$.  
Since 
$$
\|T_1\xi_1+\dots+T_d\xi_d\|^2\leq 
\|\xi_1\|^2+\dots+\|\xi_d\|^2
$$ 
for all $\xi\in H$, it follows that 
$$
0\leq T_1T_1^*+\dots+T_dT_d^*\leq \bold 1_H,
\tag{1.6}
$$
hence we can define a positive operator $\Delta$ 
on $H$ by 
$\Delta = (\bold 1 - T_1T_1^*-\dots-T_dT_d^*)^{1/2}$.  
The {\it rank} of $H$ is defined as the rank of 
the operator $\Delta$, 
$$
\rank(H) = \dim \overline{\Delta H}.  
$$
$\rank(H)$ can take on any of the values 
$0,1,2,\dots,\infty$, and we have $\rank(H)=0$
iff $T_1T_1^*+\dots+T_dT_d^*=\bold 1$.  

The operators $T_1,\dots,T_d$ determine a completely 
positive map $\phi: \Cal B(H)\to \Cal B(H)$ by way of 
$$
\phi(A) = T_1AT_1^*+\dots+T_dAT_d^*, \qquad A\in\Cal B(H).
$$
$\phi$ is continuous relative to the weak operator 
topology of $\Cal B(H)$, and by virtue of (1.6) 
we have $\|\phi\|=\|\phi(\bold 1)\| \leq 1$.   
It follows that $\bold 1\geq \phi(\bold 1)\geq 
\phi^2(\bold 1)\geq \dots$ is a decreasing 
sequence of positive operators and we  
write $\phi^\infty(\bold 1)$ for the limit
$$
\phi^\infty(\bold 1) = \lim_{n\to\infty}\phi^n(\bold 1).  
\tag{1.7}
$$
Of course, $\phi^\infty(A)$ is undefined for any
operator $A$ other than the identity.  If
$\phi^\infty(\bold 1)=0$ then $H$ is called a 
{\it pure} Hilbert $A$-module.  The opposite 
extreme $\phi^\infty(\bold 1)=\bold 1$ occurs 
only when $\rank(H)=0$.  

Finally, the set of all polynomials 
$\{f(T_1,\dots, T_d): f\in A\}$ in the canonical 
operators $T_1,\dots,T_d$ of $H$ is a commutative 
subalgebra of $\Cal B(H)$ which contains the identity
operator,  and we write this algebra of operators 
as $\alg(H)$.  $C^*(H)$ denotes the \cstar\ generated 
by $\alg(H)$.  $C^*(H)$ is irreducible iff $H$ cannot 
be decomposed into an orthogonal direct sum 
$H=H_1\oplus H_2$
of nonzero submodules $H_1$, $H_2$.  

Now let 
$$
M\underset U\to\longrightarrow H\longrightarrow 0
\tag{1.8}
$$
be a dilation of $H$, and suppose that $M$ can be 
decomposed into a direct sum of submodules 
$M=M_1\oplus M_2$ where the restriction of $U$ to 
$M_2$ vanishes.  In this case we say that $M_2$ is a 
{\it trivial} summand of the dilation (1.8).  For 
example, in the dilation (1.5), $N$ is a trivial 
summand.  

\proclaim{Proposition 1.9}
Let $M\underset U\to\longrightarrow H\longrightarrow 0$ 
be a dilation of a Hilbert $A$-module $H$.  The following 
are euivalent.  
\roster
\item
$M\underset U\to\longrightarrow H\longrightarrow 0$ has 
no nonzero trivial summands.  
\item
The set of vectors
$$
C^*(M)U^*H = \{AU^*\xi: A\in C^*(M), \xi\in H\}
$$
spans $M$.  
\item
For every operator $A$ in the commutant $C^*(H)^\prime$
we have
$$
AU^*=0 \implies A=0.  
$$
\endroster
If $M=H^2\otimes E$ is a free Hilbert $A$-module
(see the following paragraphs)
then these conditions are equivalent to
\roster
\item"{(4)}"
The map of $E$ to $H$ defined by 
$\zeta\mapsto U(1\otimes\zeta)$ is one-to-one.  
\endroster
\endproclaim
\demo{proof}
This is fundamentally a restatment of the equivalence
of properties (8.4.1), (8.4.1) and (8.4.3) of \cite{1}.\qedd
\enddemo

\proclaim{Definition 1.10}
A dilation $M\underset U\to\longrightarrow H\longrightarrow 0$
of $H$ is called minimal if it satisfies 
the conditions of Proposition 1.9.  
\endproclaim

\remark{Remark}
It is clear from property (2) of Proposition 1.9 that any
dilation 
$$
M\underset U\to\longrightarrow H\longrightarrow 0
$$
can be reduced to a minimal one 
$$
M_0\underset U_0\to\longrightarrow H\longrightarrow 0
$$
by replacing $M$ with the submodule 
$M_0=[C^*(M)U^*H]$ and $U$ with its restriction 
to $M_0$. 
\endremark

We now summarize the main results on the existence 
and uniqueness of nonnormal dilations for Hilbert 
$A$-modules.  A {\it free} Hilbert $A$-module is a 
finite or countably infinite direct sum of copies 
of the rank-one module $H^2$.  We write $n\cdot H^2$ 
for the direct sum of $n$ copies of $H^2$, 
$n=1,2,\dots,\infty$.  $n$ is uniquely determined by 
the module $n\cdot H^2$; indeed, a simple computation 
(which we omit) 
shows that $\rank(n\cdot H^2) = n$ for every $n$.  
Thus we will refer to $n\cdot H^2$ as the {\it 
free Hilbert $A$-module of rank} $n=1,2,\dots,\infty$.  If 
$E$ is a Hilbert space of dimension $n$ and 
we make the Hilbert space $H^2\otimes E$ into a 
Hilbert $A$-module by setting 
$$
f(g\otimes\zeta)=fg\otimes\zeta,
\qquad f\in A,\quad g\in H^2,\quad \zeta\in E
$$
then $H^2\otimes E$ is isomorphic to $n\cdot H^2$.  

At the other extreme, a Hilbert $A$-module $H$ is 
called {\it spherical} if its canonical operators 
$\{T_1,\dots,T_d\}$ are jointly normal in the sense 
that $\{T_1,\dots,T_d, T_1^*,\dots,T_d^*\}$ is a commuting
set of operators, and in addition satisfy 
$$
T_1T_1^*+\dots+T_dT_d^* = \bold 1.  
$$ 
Spherical $d$-tuples $(T_1,\dots,T_d)$ are 
the higher dimensional counterparts of unitary 
operators.  
One cannot avoid spherical modules in 
nonnormal dilation theory, 
and they represent a kind of degeneracy.   
By a {\it standard} Hilbert $A$-module we mean 
a direct sum of Hilbert $A$-modules of the form
$H  = F\oplus S$
where $F$ is free and $S$ is spherical.  
One or the other summand may be absent.  
A {\it standard} dilation of $H$ is a dilation 
$$
M\underset U\to\longrightarrow H\longrightarrow 0
$$
in which $M$ is a standard module.  
Theorem 8.5 of 
\cite{1} can now be reformulated as follows.  

\proclaim{Theorem 1.11}
Every Hilbert $A$-module $H$ has a minimal standard 
dilation 
$$
F\oplus S \underset U\to\longrightarrow H\longrightarrow 0,
$$ 
and any two such are equivalent.  If 
$$
F^\prime\oplus S^\prime \underset U^\prime\to\longrightarrow 
H\longrightarrow 0
$$
is a second standard dilation then every isomorphism 
$W: F\oplus S\to F^\prime \oplus S^\prime $ satisfying 
$U^\prime W=U$ decomposes into 
a direct sum $W = W_1\oplus W_2$, where $W_1$ is 
an isomorphism of the free summands and $W_2$ is an 
isomorphism of the spherical summands.  
\endproclaim

There is somewhat more information available 
from \cite{1, Theorem 8.5} concerning criteria for 
when one or the other of the summands $F$ or $S$ 
is missing.  

\proclaim{Theorem 1.12}
Let $H$ be a Hilbert $A$-module.  
\roster
\item
The minimal standard dilation of $H$ is free
of rank $n=1,2,\dots,\infty$
$$
n\cdot H^2 \underset U\to\longrightarrow H\longrightarrow 0
$$
iff $H$ is pure of rank $n$.  
\item
The minimal standard dilation of $H$ is spherical
$$
S \underset U\to\longrightarrow H\longrightarrow 0
$$
iff $\phi(\bold 1) = T_1T_1^*+\dots+T_dT_d^* = \bold 1$
(i.e., iff $\rank(H)=0$).  
\endroster
\endproclaim

Finally, we require more explicit information 
about the minimal 
standard dilation of a Hilbert $A$-module $H$ 
$$
F\oplus S\underset U\to\longrightarrow H\longrightarrow 0
$$
than is apparent from Theorems 1.11 and 1.12.  
Indeed, there is an explicit formula for 
$F$ and for the restriction of $U$ to $F$ 
which is described as follows.  
Consider the module $F=H^2\otimes \overline{\Delta H}$, 
where the $A$-module structure is defined by 
$$
f\cdot(g\otimes\zeta) = fg\otimes \zeta, \qquad
f,g\in A, \zeta\in\overline{\Delta H}.  
$$ 
We have already pointed out that 
$F=H^2\otimes \overline{\Delta H}$ is 
a free Hilbert $A$-module of rank $r=\rank(H)$.  
Moreover, 
Theorem 4.5 of \cite{1} implies that there is a unique 
bounded operator $U_0: H^2\otimes \overline{\Delta H}\to H$
satisfying  
$$
U_0(f\otimes\zeta) = f\cdot \Delta\zeta, 
\qquad f\in A,\zeta\in\overline{\Delta H}.
$$
$U_0$ is obviously a homomorphism of Hilbert 
$A$-modules, and the proof of Theorem 4.5 of \cite{1} 
shows that $U_0U_0^*=\bold 1_H-\phi^\infty(\bold 1_H)$.
In particular, if $H$ is pure then $U_0$ is 
a coisometry and 
$$
H^2\otimes\overline{\Delta H} 
\underset U_0\to\longrightarrow H\longrightarrow 0
\tag{1.13}
$$
provides a standard dilation of $H$ which has no spherical 
summand.  Indeed, in this case (1.13) is actually
the {\it minimal} standard dilation (condition (4) 
of Proposition 1.9 is obviously satisfied).  

If $H$ is not 
pure then (1.13) is not a standard dilation 
(nor is it even a dilation); but 
one can make it so by adding an 
appropriate spherical summand.  The details are 
as follows.  

\proclaim{Theorem 1.14}
Let $H$ be a Hilbert $A$-module, let 
$F=H^2\otimes \overline{\Delta H}$, and let 
$U_0: F\to H$ be the morphism defined in (1.13).  
Then there is a spherical module $S$ and a
morphism $U_1: S\to H$ such that 
$U_1U_1^* = \phi^\infty(\bold 1_H)$,  and such that 
if $U: F\oplus S\to H$ is defined by 
$$
U(\xi,\eta) = U_0\xi + U_1\eta
$$
then 
$$
F\oplus S \underset U\to\longrightarrow H\longrightarrow 0
$$
is the minimal standard dilation of $H$.  
\endproclaim
\demo{proof}
Let 
$$
F\oplus S\underset V\to\longrightarrow H\longrightarrow 0
$$
be a minimal standard dilation of $H$.  
Given the formulas of the preceding paragraph, it suffices 
to show that the restriction of $V$ to $F$ can be 
identified with $U_0$ in the sense that there is an 
isomorphism of Hilbert $A$-modules
$W:F\to H^2\otimes \overline{\Delta H}$ such 
that $V\restriction_F=U_0W$.  Now both $F$ and 
$H^2\otimes \overline{\Delta H}$ are free modules 
of the same rank $r=\rank(H)$, and thus we may identify
$F$ with $H^2\otimes \overline{\Delta H}$.  

Having made that identification, let $V_0$ be the 
restriction of $V$ to the free summand
$F=H^2\otimes \overline{\Delta H}$
and consider the linear
operator $A:\overline{\Delta H}\to H$ defined by
$$
A\zeta = V_0(1\otimes \zeta), 
$$
$1$ denoting the 
constant function in $H^2$.  We claim that 
$A$ has trivial kernel and $AA^*=\Delta^2$.  Granting
that for a moment, the polar decomposition provides a
unitary operator $W_0\in\Cal B(\overline{\Delta H})$
such that $A=\Delta W_0$, hence
$$
\align
V_0(f\otimes\zeta)&=f\cdot V_0(1\otimes\zeta)=f\cdot A\zeta
=f\cdot\Delta W_0\zeta \\
&= U_0(f\otimes W_0\zeta) =
U_0((\bold 1_{H^2}\otimes W_0)(f\otimes\zeta)), 
\endalign
$$
thus $W=\bold 1_{H^2}\otimes W_0$ provides the 
required automorphism of $H^2\otimes\overline{\Delta H}$ 
satisfying $V_0=U_0 W$.  

To see that $A$ is injective, pick 
$\zeta\in\overline{\Delta H}$ such that $A\zeta=0$, 
and consider the operator $Z$ defined on the 
dilation module 
$M=H^2\otimes\overline{\Delta H}\oplus S$ by 
$$
Z=\bold 1_{H^2}\otimes (\zeta\otimes \bar\zeta)\oplus 0,   
$$
$\zeta\otimes\bar\zeta$ denoting the rank-one operator 
on $\overline{\Delta H}$ associated with the vector $\zeta$. 
$Z$ is a self-adjoint operator in the commutant of 
$C^*(M)$, and since $V$ maps the vector 
$(1\otimes\zeta,0)\in M$ to $A\zeta=0$ we have $VZ=0$ and
hence $ZV^*=0$.  By the third condition of Proposition 1.9, 
$Z$ must be the zero operator and hence $\zeta = 0$.  

Finally, we show that $AA^*=\Delta^2$.  For that, let 
$T_1,\dots,T_d$ (resp. $\tilde T_1,\dots,\tilde T_d$) be 
the canonical operators associated with the Hilbert 
module $H$ (resp. $M$).  Since 
$M=H^2\otimes \overline{\Delta H}\oplus S$ and $S$ is 
spherical, we see from (1.2) that the projection of 
$M$ onto the subspace 
$1\otimes\overline{\Delta H}\oplus 0$ is given by
$\bold 1_M-\sum_{k=1}^d \tilde T_k\tilde T_k^*$.  
Since $VV^*=\bold 1_H$ and $V\tilde T_k=T_kV$ 
for every $k=1,\dots,d$ we have 
$$
V(\sum_{k=1}^d \tilde T_k\tilde T_k^* )V^* =
\sum_{k=1}^d T_kVV^*T_k^* = \sum_{k=1}^dT_kT_k^*,
$$
hence 
$$
AA^* = V(\bold 1_M-\sum_{k=1}^d \tilde T_k\tilde T_k^*)V^* 
=VV^*-\sum_{k=1}^d T_kT_k^* = \Delta^2
$$
as required.  \qedd
\enddemo

\example{1.15.  Dilations and free resolutions}
We conclude this section with some remarks concerning 
free resolutions of {\it pure} Hilbert modules.  
Theorem 1.12 shows that
$H^2$ and its multiples (free modules) 
occupy a key position in nonnormal dilation
theory, in that a Hilbert module $H$ is 
isomorphic to some quotient $F/K$ of a free module 
$F$ iff $H$ is pure.  This leads to the existence 
(and uniqueness) of minimal free resolutions in the 
category of pure Hilbert $A$-modules.  More 
precisely, suppose  
we start with a pure Hilbert module $H$.  Applying Theorem 
1.12 we find a minimal dilation of $H$ of the form 
$$
F\longrightarrow H\longrightarrow 0
$$
where $F$ is free.  
Let $K\subseteq F$ be the kernel of this dilation map.  
Free modules are pure, and submodules of pure 
modules are pure.  Hence, we can reapply 
Theorem 1.12 to $K$ itself 
and continue the sequence one step to the left.  
We then repeat
the process on the kernel of the resulting map, 
continuing indefinitely or until some kernel $K$ is zero.  

Thus, {\it every pure Hilbert $A$-module 
has a minimal free resolution}, that is, there 
is an exact sequence of Hilbert $A$-modules 
$$
\dots\longrightarrow F_3
\longrightarrow F_2\longrightarrow F_1
\longrightarrow H\longrightarrow 0
\tag{1.16}
$$
where each $F_k$ is a free module 
and where each connecting map is a partial isometry
which satisfies the conditions of Proposition 1.9.  
This minimal free resolution is unique in the sense that 
if 
$$
\dots\longrightarrow F_3^\prime
\longrightarrow F_2^\prime\longrightarrow F_1^\prime
\longrightarrow H\longrightarrow 0
$$
is another one then after working from right to 
left with Theorem 1.11 we find 
that there is a sequence of 
isomorphisms $W_k:F_k\to F_k^\prime$ 
such that each subdiagram
$$
\CD
F_{k+1}		@>>>			 F_k\\
@V W_{k+1} VV						@VV W_kV \\
F_{k+1}^\prime 	@>>>			F_k^\prime
\endCD
$$
commutes.  

Naturally, one would like to know if the resolution 
(1.16) is of finite length in the sense that 
$F_k=0$ for sufficiently large $k$.  
Moreover, if we start with a pure 
module $H$ of {\it finite} rank, then by 
analogy with Hilbert's syzygy theorem one might hope 
that the free resolution of $H$ is a) of finite length
and b) each of the free modules $F_k$ is of finite rank.  
Unfortunately, nothing like that is true in this 
category.  We will show in section 9 that 
every pure {\it graded}
Hilbert $A$-module $H$ of rank one which 
is not already isomorphic to $H^2$ has a 
minimal free 
resolution of the form 
$$
\dots
\longrightarrow F_2\longrightarrow F_1
\longrightarrow H\longrightarrow 0
$$
where $F_1$ is isomorphic to $H^2$ but where 
$F_2$ is free of {\it infinite} rank.  Thus, 
Hilbert's syzygy theorem fails
rather spectacularly for Hilbert $A$-modules.  
\endexample

Nevertheless, we will find in section 3 that
{\it algebraic} free resolutions do play a central 
role in the theory of Hilbert $A$-modules.

\subheading{2.  Multipliers of free Hilbert modules}

Elements of free modules, and homomorphisms from one free 
module to another, can be 
``evaluated" at points in the open unit ball 
$B_d$ in $\Bbb C^d$.  The purpose of this section is to 
discuss these evaluation maps and the relation between 
morphisms and multipliers.  Throughout the section, 
$\{u_z: z\in B_d\}$ will denote the 
set of functions in $H^2$ 
defined in section 1 (see (1.1.1) and (1.1.2)).  

Let $E$ be a separable Hilbert space and
consider the free Hilbert $A$-module $F=H^2\otimes E$
of rank $r=\dim E$, where the module 
structure is defined by 
$$
f\cdot(g\otimes \zeta) = fg\otimes \zeta, 
\qquad f,g\in A, \zeta\in E.  
$$
One can think of elements of $H^2\otimes E$ 
as $E$-valued holomorphic functions defined on $B_d$ 
in the following way.  

\proclaim{Proposition 2.1}
For every element $\xi\in H^2\otimes E$ and every 
$z\in B_d$ there is a unique vector $\hat \xi(z)\in E$ 
satisfying 
$$
\<\hat \xi(z),\zeta\> = \<\xi,u_z\otimes \zeta\>, \qquad \zeta\in E.  
$$
The function $\hat \xi:B_d\to E$ is (weakly) holomorphic and 
satisfies 
$$
\|\hat \xi(z)\| \leq \frac{\|\xi\|}{\sqrt{1-|z|^2}}, \qquad z\in B_d.
\tag{2.2}
$$
\endproclaim

\demo{proof}
The argument is straightforward and we merely 
sketch the details.  Since
$$
\|u_z\otimes \zeta\|  = \|u_z\|\cdot\|\zeta\| = 
\frac{\|\zeta\|}{\sqrt{1-|z|^2}}
$$
it follows that for fixed $z\in B_d$, 
$\zeta\mapsto \<\xi,u_z\otimes \zeta\>$ is a bounded 
antilinear functional on $E$.  By the Riesz lemma there
is a unique vector $\hat \xi(z)\in E$ such that 
$$
\<\hat \xi(z),\zeta\>=\<\xi,u_z\otimes \zeta\>,\qquad \zeta\in E, 
$$
and one has $\|\hat \xi(z)\|\leq \|\xi\|(1-|z|^2)^{-1/2}$.  
By the definition of $u_z$, $z\mapsto \<\xi,u_z\otimes \zeta\>$ 
is holomorphic in $B_d$, hence 
$\hat \xi: B_d\to E$ is weakly holomorphic.  The estimate 
(2.2) is immediate from the definition of $\hat \xi$.\qedd
\enddemo

\remark{Remarks}
We lighten notation by writing $\xi(z)$ for the 
value of $\xi$ at $z$, rather than the more pedantic $\hat \xi(z)$.  
Notice that the $A$-module structure of $H^2\otimes E$ is 
expressed conveniently in terms of the values of $\xi$ as 
follows
$$
(f\cdot \xi)(z) = f(z)\xi(z), \qquad f\in A,
\quad \xi\in H^2\otimes E,\quad z\in B_d.  
$$

Notice too that an element $\xi\in H^2\otimes E$ is uniquely
determined by its functional representative $z\in B_d\mapsto \xi(z)$
because $\{u_z: z\in B_d\}$ spans $H^2$ \cite{1}, hence 
$\{u_z\otimes\zeta: z\in B_d,\zeta\in E\}$ spans $H^2\otimes E$, 
and hence $\xi(z)=0$ for all $z\in B_d$ only when $\xi$ is 
the zero element of $H^2\otimes E$.  
\endremark

Similarly, any bounded homomorphism of free modules can be 
evaluated at points of the unit ball so as to obtain a multiplier.  
In more detail, let $E_1$ and $E_2$ be two separable Hilbert spaces
and let $\Phi: B_d\to \Cal B(E_1,E_2)$ be an operator-valued 
function defined on the open unit ball.  We will say that $\Phi$ 
is a {\it multiplier} if there is a bounded linear operator 
$\hat \Phi: H^2\otimes E_1\to H^2\otimes E_2$ such that 
$$
\Phi(z)\xi(z) = (\hat\Phi \xi)(z), 
\qquad \xi\in H^2\otimes E_1,\quad z\in B_d.  
$$
The multiplier norm of $\Phi$ is defined as the operator 
norm 
$$
\|\Phi\|_{\Cal M} = \|\hat\Phi\| = 
\sup_{\|\xi\|_{H^2\otimes E_1}\leq 1}
\|\hat\Phi \xi\|_{H^2\otimes E_2}.  
$$
Notice that the operator $\hat\Phi$ is a homomorphism
of the $A$-module structure of $H^2\otimes E_1$ to 
that of $H^2\otimes E_2$, since by the preceding remarks
for $f\in A$ and $\xi\in H^2\otimes E_1$ we have 
$$
\align
\hat\Phi(f\cdot \xi)(z)&= \Phi(z)((f\cdot\xi)(z))=
\Phi(z)(f(z)\xi(z))= f(z)\Phi(z)\xi(z) \\&=
f(z)\hat\Phi(\xi)(z)=(f\cdot\hat\Phi(\xi))(z)
\endalign
$$
for all $z\in B_d$, hence 
$\hat\Phi(f\cdot \xi)=f\cdot\hat\Phi(\xi)$.  
The space of all such multipliers is denoted 
$\Cal M(E_1,E_2)$, and we again simplify notation by dropping 
the circumflex over the operator $\hat \Phi$.  

Now let $E_1, E_2$ be two separable Hilbert spaces.  
A bounded linear operator 
$\Phi: H^2\otimes E_1\to H^2\otimes E_2$ satisfying
$\Phi(f\cdot\xi)=f\cdot\Phi(\xi)$ will be called 
a {\it homomorphism}, and in this 
section only the Banach space of 
all such will be written 
$\Cal Hom (H^2\otimes E_1,H^2\otimes E_2)$
(recall that we have reserved the notation $\hom(H,K)$ for 
spaces of homomoprhisms having norm at most $1$).  

Fix $\Phi\in \Cal Hom(H^2\otimes E_1,H^2\otimes E_2)$.  Then for 
every $z$ in the open unit ball $B_d$ we can define an operator
$\Phi(z)\in \Cal B(E_1,E_2)$ by 
$$
\<\Phi(z)\zeta_1,\zeta_2\> = 
\<\Phi(1\otimes\zeta_1),u_z\otimes \zeta_2\>, \qquad \zeta_k\in E_k;  
\tag{2.3}
$$
indeed, this follows from an application of the Riesz
lemma after taking note of the obvious estimate of the term on 
the right
$$
|\<\Phi(1\otimes\zeta_1),u_z\otimes \zeta_2\>|\leq
\frac{\|\Phi\|\cdot\|\zeta_1\|\cdot\|\zeta_2\|}{\sqrt{1-|z|^2}}.
$$

\proclaim{Proposition 2.4}
For every homomorphism $\Phi: H^2\otimes E_1\to H^2\otimes E_2$, 
the operator function $z\in B_d\mapsto \Phi(z)\in\Cal B(E_1,E_2)$
defined by (2.3) 
belongs to the multiplier space $\Cal M(E_1,E_2)$ and its 
associated operator is $\Phi$.  Moreover, 
\roster
\item
$\sup_{|z|<1}\|\Phi(z)\|\leq \|\Phi\|=\|\Phi\|_\Cal M$, and 
\item
the adjoint $\Phi^*\in \Cal B(H^2\otimes E_2,H^2\otimes E_1)$
of the operator $\Phi$ is related to the operator function 
$z\in B_d\mapsto \Phi(z)^*\in\Cal B(E_2,E_1)$ as follows,
$$
\Phi^*(u_z\otimes \zeta)=u_z\otimes \Phi(z)^*\zeta, 
\qquad z\in B_d, \quad \zeta\in E_2.  
$$
\endroster
\endproclaim

\demo{proof}
We claim first that for every $\zeta_2\in E_2$ and 
every $\xi\in H^2\otimes E_1$ of the form $\xi=f\otimes \zeta$
we have
$$
\<\Phi(z)\xi(z),\zeta_2\>=\<\Phi(\xi),u_z\otimes\zeta_2\>.  
\tag{2.5}
$$
Indeed, since $\xi(z)=f(z)\zeta_2$ the left side of (2.5) is 
$f(z)\<\Phi(z)\zeta_1,\zeta_2\>$ while since 
$$
\Phi(\xi) = \Phi(f\cdot(1\otimes\zeta_1))=
f\cdot \Phi(1\otimes\zeta_1),
$$
the right side of (2.5) is 
$$
\<f\cdot\Phi(1\otimes\zeta_1),u_z\otimes \zeta_2\> =
\<f(z)\Phi(1\otimes\zeta_1)(z), \zeta_2\> =
f(z)\<\Phi(z)\zeta_1,\zeta_2\>.  
$$
Hence (2.5) holds for elementary tensors $\xi$.  
Now for an arbitrary $\xi\in H^2\otimes E_1$ the left 
side of (2.5) is bounded by 
$$
|\<\Phi(z)\xi(z),\zeta_2\>|\leq 
\|\Phi(z)\|\cdot\|\xi(z)\|\cdot \|\zeta_2\| \leq
\frac{\|\Phi\|\cdot\|\xi\|\cdot\|\zeta_2\|}{\sqrt{1-|z|^2}}
$$
and the right side is bounded by 
$$
|\<\Phi(\xi),u_z\otimes \zeta_2\>|\leq 
\|\Phi\|\cdot\|\xi\|\cdot\|u_z\otimes \zeta_2\|\leq
\frac{\|\Phi\|\cdot\|\xi\|\cdot\|\zeta_2\|}{\sqrt{1-|z|^2}}.  
$$
Since $H^2\otimes E_1$ is spanned by elements of the form
$f\otimes \zeta_1$, (2.5) follows in general.  

Now by definition of $\Phi(\xi)(z)$ the right side of 
(2.5) is $\<\Phi(\xi)(z),\zeta_2\>$, and since 
$\zeta_2$ is arbitrary in (2.5) we conclude 
that $\Phi(\xi)(z)=\Phi(z)\xi(z)$.  Hence the function
$\Phi(\cdot)$ is a multiplier with associated operator
$\Phi\in\Cal B(h^2\otimes E_1,H^2\otimes E_2)$.  

We now verify formula (2) of Proposition 2.4.  Fix 
$f\in H^2$, $\zeta_k\in E_k$, $k=1,2$, and $z\in B_d$.  
Then we have
$$
\align
\<f\otimes\zeta_1,\Phi^*(u_z\otimes\zeta_2)\> &=
\<\Phi(f\otimes\zeta_1),u_z\otimes\zeta_2\>=
\<f\cdot\Phi(1\otimes\zeta_1),u_z\otimes\zeta_2\>\\
&=
\<f(z)\Phi(z)\zeta_1,\zeta_2\> =
f(z)\<\Phi(z)\zeta_1,\zeta_2\> =
\<f,u_z\>\<\zeta_1,\Phi(z)^*\zeta_2\> \\
&=
\<f\otimes\zeta_1,u_z\otimes\Phi(z)^*\zeta_2\>.  
\endalign
$$
Since $H^2\otimes E_1$ is spanned by vectors of the form
$f\otimes\zeta_1$, the required formula (2) follows.  

To prove (1) of Proposition 2.4 it suffices to show that
for every $\zeta_k\in E_k$, $k=1,2$ with $\|\zeta_k\|\leq 1$
we have $|\<\Phi(z)\zeta_1,\zeta_2\>|\leq \|\Phi\|$.  
For that, write
$$
(1-|z|^2)^{-1}|\<\Phi(z)\zeta_1,\zeta_2\>| =
\|u_z\|^2 |\<\zeta_1,\Phi(z)^*\zeta_2\>| =
|\<u_z\otimes\zeta_1,u_z\otimes\Phi(z)^*\zeta_2\>|.  
$$
By the formula (2) just established, the right side is 
$$
|\<u_z\otimes\zeta_1,\Phi^*(u_z\otimes\zeta_2)\>| \leq 
\|u_z\|^2 \|\Phi^*\| = (1-|z|^2)^{-1}\|\Phi\|,
$$
from which the assertion of (1) follows.  \qedd
\enddemo

\remark{Remarks}
Experience with one-dimensional operator theory 
might lead one to expect that the space of all 
multipliers 
$\Cal M(E_1,E_2)$ should coincide with the 
space $H^\infty(E_1,E_2)$ of all bounded 
holomorphic operator valued functions 
$$
F: B_d\to \Cal B(E_1,E_2).  
$$
However, the 
failure of von Neumann's inequality for the 
ball $B_d$ in dimension $d\geq 2$ (cf. \cite{1}, 
Theorem 3.3)
implies that this is far from true.  Indeed, if we 
consider the simplest case in which both spaces 
$E_1=E_2=\Bbb C$ consist of scalars, then 
$\Cal M(\Bbb C,\Bbb C) = \Cal M$ is the multiplier algebra
introduced in \cite{1}, and it was 
shown there that $\Cal M$ is a proper subalgebra of 
the algebra $H^\infty$ of all bounded holomorphic functions
defined on the open unit ball in $\Bbb C^d$
when $d\geq 2$.  Indeed, examples
are given in (\cite{1}, Theorem 3.3)
of continuous functions defined 
on the closed unit ball $f: \overline{B_d}\to \Bbb C$ 
which are holomorphic in the interior $B_d$,  but which 
are not multipliers.  
\endremark

We conclude this section with a few remarks about 
boundary values.  Let $\sigma$ denote the natural normalized
measure on the unit sphere $\partial B_d$ in complex $d$-space, 
and let $H^2(\partial B_d)$ denote the multivariate ``Hardy 
space" defined as the closure in $L^2(\partial B_d)$ of 
the holomorphic polynomials.  Every element $\tilde f$ 
of $H^2(\partial B_d)$ has a natural holomorphic extension 
$f$
to the interior the ball $B_d$, and for $\sigma$-almost every 
$z\in \partial B_d$ we have 
$$
\lim_{r\uparrow 1}f(r z) = \tilde f(z).  
$$
Moreover, 
$$
\lim_{r\uparrow 1}
\int_{\partial B_d}|f(r z)-\tilde f(z)|^2\,d\sigma(z) = 0
$$
(for example, see \cite{26}).  

These properties generalize to 
vector-valued functions as follows.  Let $E$ be 
a separable Hilbert space and let $H^2(\partial B_d;E)$ denote 
the closure in 
$L^2(\partial B_d;E)=L^2(\partial B_d)\otimes E$ of the linear 
span of all vector polynomials of the form $f\otimes \zeta$, 
with $f\in A$ and $\zeta\in E$.  Elements $\tilde \xi$ of 
$H^2(\partial B_d;E)$ extend in a similar way to holomorphic 
functions $f: B_d\to E$ and there is a Borel set 
$N\subseteq\partial B_d$ of measure zero such that for all 
$z\in \partial B_d\setminus N$ we have
$$
\lim_{r\uparrow 1}\|\xi(r z) -\tilde\xi(z)\| = 0, 
\tag{2.6}
$$
and moreover
$$
\lim_{r\uparrow 1}\int_{\partial B_d}
\|\xi(r z)-\tilde \xi(z)\|^2\, d\sigma(z) = 0; 
\tag{2.7}
$$
for example, one can establish this by making use 
of the radial maximal function, see 5.4.11 of \cite{26}.  

The preceding remarks lead immediately to the following 
conclusion about the boundary values of multipliers.  

\proclaim{Proposition 2.8}
Let $\Phi\in \Cal M(E_1,E_2)$ be the multiplier of a 
homomorphism $\Phi$ in $\Cal Hom(H^2\otimes E_1,H^2\otimes E_2)$.  
Then for $\sigma$-almost every $z\in \partial B_d$, the 
net of operators $r\in (0,1)\mapsto \Phi(r z)\in\Cal B(E_1,E_2)$
is uniformly bounded and 
converges in the strong operator topology of $\Cal B(E_1,E_2)$ 
to an operator $\tilde\Phi(z)$.  The 
operator function 
$\tilde \Phi:\partial B_d\to \Cal B(E_1,E_2)$ belongs 
to $L^\infty(\partial B_d;\Cal B(E_1,E_2))$ and satisfies
$$
\Sb{\text{ess sup}}\\{|z|=1}\endSb\|\tilde\Phi(z)\| = 
\sup_{|z|<1}\|\Phi(z)\| \leq \|\Phi\|_{\Cal M}.  
\tag{2.8}
$$
\endproclaim

\demo{proof}
The argument is straightforward, and we merely sketch the
details.  For fixed $\zeta\in E_1$, consider the holomorphic
vector-valued function $\xi: B_d\to E_2$ defined by 
$\xi(z) = \Phi(z)\zeta$.  $\xi$ satisfies

$$
\sup_{|z|<1}\|\xi(z)\| \leq 
\|\Phi(z)\|_\Cal M \cdot\|\zeta\|<\infty
$$
and therefore it is the restriction to $B_d$ of a unique
element $\tilde\xi\in H^2(\partial B_d; E_2)$.  Moreover, 
we have 
$$
\lim_{r\uparrow 1}\|\xi(r z)-\tilde \xi(z)\| = 0
$$
for $\sigma$-almost every $z\in \partial B_d$.  Since 
$E_1$ is separable, a standard argument shows that the exceptional 
set can be made independent of $\zeta\in E_1$, and for all 
such nonexceptional points 
$z\in\partial B_d$ the net of operators 
$r\mapsto \Phi(r z)$ is strongly convergent to a limit 
operator $\tilde \Phi(z)$ satisfying 2.9.  \qedd
\enddemo

\subheading{3.  Euler characteristic}

Throughout this section, $H$ will denote a {\it finite rank}
Hilbert $A$-module.  We will work not with $H$ itself but 
with the following linear submanifold of $H$
$$
M_H=\sp\{f\cdot\Delta\xi: f\in A, \xi\in H\}.  
$$
The definition and basic properties of the Euler 
characteristic are independent of any topology 
associated with the Hilbert space $H$, and depend
solely on the linear algebra of $M_H$.  As we have 
pointed out in the introduction, $M_H$ is a finitely 
generated $A$-module, and has finite free resolutions 
in the category of finitely generated $A$-modules
$$
0\longrightarrow F_n\longrightarrow\dots
\longrightarrow F_2\longrightarrow
F_1\longrightarrow M_H\longrightarrow 0,
$$
each $F_k$ being a sum of $\beta_k$ copies of the 
rank-one module $A$.  The alternating sum of the 
ranks $\beta_1-\beta_2+\beta_3-+\dots$ does not 
depend on the particular free resolution of $M_H$,
and we define the {\it Euler characteristic}
of $H$ by 
$$
\chi(H) = \sum_{k=1}^n(-1)^{k+1}\beta_k.  
\tag{3.1}
$$
The main result of this section is an asymptotic formula
(Theorem C) which expresses $\chi(H)$ in terms of the
sequence of defect operators $\bold 1-\phi^{n+1}(\bold 1)$,
$n=1,2,\dots$, where $\phi$ is the completely positive 
map on $\Cal B(H)$ associated with the canonical 
operators $T_1,\dots,T_d$ of $H$,
$$
\phi(A) = T_1AT_1^*+\dots+T_dAT_d^*.  
$$

The Hilbert polynomial is an invariant associated with 
finitely generated graded modules over polynomial rings
$k[x_1,\dots,x_d]$, $k$ being an arbitrary field.  We 
require something related to the Hilbert polynomial, 
which exists in greater generality than the former, but 
whose existence can be deduced rather easily from 
Hilbert's original work \cite{16}, \cite{17}.  
While this polynomial is very fundamental (indeed, its 
existence might be described as {\it the} fundamental 
result of multivariable linear algebra), it is 
less familiar to analysts than it is 
to algebraists.  

We define this polynomial in a way suited to 
our needs below.  It is convenient to 
work in terms of the following  
sequence of polynomials $q_0, q_1,\dots\in \Bbb Q[x]$, 
which are normalized so  that 
$q_k(0)=1$, and which are defined recursively by 
$$
\align
q_0(x) &=1,\tag{3.2.1}\\
q_k(x)-q_k(x-1) &= q_{k-1}(x), \qquad k\geq 1.  
\tag{3.2.2}
\endalign
$$
One finds that for $k\geq 1$, 
$$
q_k(x) = \frac{(x+1)(x+2)\dots(x+k)}{k!}.  
\tag{3.3}
$$
When $x=n$ is a positive integer, $q_k(n)$ is 
the binomial coefficient  $\binom{n+k}n$, and 
more generally $q_k(\Bbb Z)\subseteq \Bbb Z$, 
$k=0,1,2,\dots$.  

\proclaim{Theorem 3.4}
Let $V$ be a vector space over a field $k$, let 
$T_1,\dots,T_d$ be a commuting set of linear 
operators on $V$, and make $V$ into a 
$k[x_1,\dots,x_d]$-module 
by setting $f\cdot\xi = f(T_1,\dots,T_d)\xi$, 
$f\in k[x_1,\dots,x_d]$, $\xi\in V$.   

Let $G$ be a finite dimensional subspace of $V$ and,
for every $n=0,1,2,\dots$ define a finite dimensional 
subspace $M_n$ by
$$
M_n = \sp\{f\cdot \xi: f\in k[x_1,\dots,x_d], 
\quad\deg f\leq n, \quad\xi\in G\}.  
$$
Then there are integers $c_0, c_1,\dots,c_d\in \Bbb Z$
and $N\geq 1$ such that for all $n\geq N$ we have 
$$
\dim M_n = c_0 q_0(n)+c_1q_1(n)+\dots+c_dq_d(n).  
$$
In particular, the dimension function $n\mapsto \dim M_n$
is a polynomial for sufficiently large $n$.  
\endproclaim

\demo{proof}
We may obviously assume that $V=\cup_n M_n$, and hence
$V$ is a finitely generated $k[x_1,\dots,x_d]$-module.  
The fact that the function $n\mapsto \dim M_n$ is a 
polynomial of degree at most $d$ for sufficiently 
large $n$ follows from the result in section 8.4.5 of
\cite{19}; 
and the specific form of this polynomial follows from 
the discussion in \cite{19}, section 8.4.4. 
\qedd
\enddemo

\remark{Remark 3.4}
We emphasize that the dimension function $n\mapsto \dim M_n$
is generally {\it not} a polynomial for all $n\in \Bbb N$, 
but only for sufficiently large $n\in \Bbb N$.  

We also point out for the interested reader that one can give 
a relatively simple direct proof of Theorem 3.4 by an inductive 
argument on the number $d$ of operators, along lines 
similar to the proof of Theorem 1.11 of \cite{14}.  
\endremark

Suppose now that $G$ is a finite dimensional subspace 
of $V$ which {\it generates} 
$V$ as a $k[x_1,\dots,x_d]$-module
$$
V=\sp\{f\cdot \xi: f\in k[x_1,\dots,x_d],\quad \xi\in G\}.  
$$
The polynomial 
$$
p(x) = c_0q_0(x)+c_1q_1(x)+\dots+c_dq_d(x)
$$
defined by theorem 3.4 obviously depends on the generator $G$; 
however, its top coefficient $c_d$ does not.  In order 
to discuss that, it is convenient to broaden the context 
somewhat.  Let $M$ be a module over the polynomial 
ring $k[x_1,\dots,x_d]$.  A {\it filtration} of 
$M$ is an increasing sequence 
$M_1\subseteq M_2\subseteq \dots$ 
of finite dimensional linear subspaces of $M$ such that 
$$
\align
M&=\cup_nM_n\qquad\text{and}\\
x_kM_n&\subseteq M_{n+1}, \qquad k=1,2,\dots,d,\quad n\geq 1.  
\endalign
$$
The filtration $\{M_n\}$ is called {\it proper} if 
there is an $n_0$ such that 
$$
M_{n+1}=M_n+x_1M_n+\dots+ x_dM_n, \qquad n\geq n_0.  
\tag{3.5}
$$

\proclaim{Proposition 3.7}
Let $\{M_n\}$ be a proper filtration of $M$.  Then the 
limit 
$$
c=d!\lim_{n\to\infty}\frac{\dim M_n}{n^d}
$$
exists and defines a nonnegative integer $c=c(M)$ 
which is the same for all proper filtrations.  
\endproclaim

\demo{proof}
Let $\{M_n\}$ be a proper filtration, choose $n_0$ satisfying 
(3.5), and let $G$ be the generating subspace 
$G=M_{n_0}$.  One finds that for $n=0,1,2,\dots$
$$
M_{n_0+n}=\sp \{f\cdot \xi: \deg f\leq n,\quad \xi\in M_{n_0}\}
$$
and hence there is a polynomial $p(x)\in \Bbb Q[x]$ of the 
form stipulated in Theorem 3.4 such that 
$\dim M_{n_0+n} = p(n)$
for sufficiently large $n$.  Writing 
$$
p(x) = c_0q_0(x)+c_1q_1(x)+\dots+c_dq_d(x)
$$
and noting that $q_k$ is a polynomial of degree $k$ with 
leading coefficient $1/k!$, we find that 
$$
c_d=d!\lim_{n\to\infty}\frac{p(n)}{n^d} =
d!\lim_{n\to \infty}\frac{\dim M_{n_0+n}}{n^d} =
d!\lim_{n\to \infty}\frac{\dim M_n}{n^d}, 
$$
as asserted.  

Now let $\{M_n^\prime\}$ be another proper filtration.  
Since $M=\cup_nM_n^\prime$ and $M_{n_0}$ is finite 
dimensional, there is an $n_1\in\Bbb N$ such that 
$M_{n_0}\subseteq M_{n_1}^\prime$.  Since $\{M_n^\prime\}$
is also proper we can increase 
$n_1$ if necessary to arrange the condition of (3.5) on
$M_n^\prime$ for all $n\geq n_1$, and hence
$$
M_{n_1+n}^\prime = \sp \{f\cdot\xi: \deg f\leq n, \quad 
\xi\in M_{n_1}^\prime\}.  
$$
Letting $c_d^\prime$ be the top coefficient of the 
polynomial $p^\prime(x)$ satisfying 
$$
\dim M_{n_1+n}^\prime = p^\prime(n)
$$
for sufficiently large $n$, the preceding argument shows
that 
$$
c_d^\prime = d!\lim_{n\to \infty}\frac{\dim M_n^\prime}{n^d}.  
$$

On the other hand, the inclusion 
$M_{n_0}\subseteq M_{n_1}^\prime$, 
together with the condition (3.5) on both $\{M_n\}$ 
and $\{M_n^\prime\}$,  implies 
$$
\align
M_{n_0+n}&=\sp\{f\cdot\xi: \deg f\leq n,\quad \xi\in M_{n_0}\} \\
&\subseteq \sp \{f\cdot\eta: \deg f\leq n,\quad \eta\in M_{n_1}^\prime\}
= M_{n_1+n}^\prime.  
\endalign
$$
Thus we have 
$$
\lim_{n\to\infty}\frac{\dim M_n}{n^d} = 
\lim_{n\to\infty}\frac{\dim M_{n_0+n}}{n^d} \leq 
\lim_{n\to\infty}\frac{\dim M_{n_1+n}^\prime}{n^d}=
\lim_{n\to\infty}\frac{\dim M_n^\prime}{n^d},  
$$
from which we conclude that $c_d\leq c_d^\prime$. By 
symmetry we also have $c_d^\prime\leq c_d$.  
\qedd 
\enddemo

The following two results together constitute a variant 
of the Artin-Rees lemma of commutative algebra 
(cf. \cite{29}, page II-9).  Since the result we require 
is formulated differently than the Artin-Rees lemma 
(normally a statement about the behavior of 
{\it decreasing} filtrations associated with ideals 
and their relation to submodules), 
and since we have been unable to locate an appropriate 
reference, we have included complete proofs.  

Associated with any filtration $\{M_n\}$ of a 
$k[x_1,\dots,x_d]$-module $M$ there is an associated 
$\Bbb Z$-graded module $\bar M$, which is defined as 
the (algebraic) direct sum of finite dimensional 
vector spaces
$$
\bar M = \sum_{n\in \Bbb Z} \bar M_n, 
$$
where $\bar M_n=M_n/M_{n-1}$ 
for each $n\in \Bbb Z$, and where for nonpositive 
values of $n$, $M_n$ is taken as $\{0\}$.  The 
$k[x_1,\dots,x_d]$-module structure on $\bar M$ is 
defined by the commuting $d$-tuple of ``shift" operators
$T_1,\dots,T_d$, where $T_k$ is defined on each summand
$\bar M_n$ by 
$$
T_k: \xi + M_{n-1}\in M_n/M_{n-1}\mapsto 
x_k\xi + M_n \in M_{n+1}/M_n.  
$$
\remark{Remark 3.7}
For our purposes, the essential feature of this 
construction is that for every $n\geq 1$, the following 
are equivalent
\roster
\item
$M_{n+1} = M_n + x_1M_n+\dots +x_dM_d$ 
\item
$\bar M_{n+1} = T_1\bar M_n+\dots +T_d\bar M_n$.  
\endroster
\endremark

\proclaim{Lemma 3.8}
Let $\{M_n\}$ be a filtration of a $k[x_1,\dots,x_d]$-module
$M$.  The following are equivalent: 
\roster
\item
$\{M_n\}$ is proper.  
\item 
The $k[x_1,\dots,x_d]$-module $\bar M$ is finitely generated.  
\endroster
\endproclaim

\demo{proof of (1) $\implies$ (2)}
Find an $n_0\in \Bbb N$ such that 
$$
M_{n+1} = M_n + x_1M_n+\dots+ x_dM_n
$$
for all $n\geq n_0$.  From Remark 3.7 we have
$\bar M_{n+1} = T_1\bar M_n+\dots +T_d\bar M_n$ 
for all $n\geq n_0$, hence $G = \bar M_1+\dots+\bar M_{n_0}$ 
is a finite dimensional generating space for $\bar M$.  
\enddemo
\demo{proof of (2) $\implies$ (1)}
Assuming (2), we can find a finite set of homogeneous 
elements $\xi_k\in \bar M_{n_k}$, $k=1,\dots,r$ which 
generate $\bar M$ as a $k[x_1,\dots,x_d]$-module.  It 
follows that for $n\geq \max(n_1,\dots,n_r)$ we have 
$$
\bar M_{n+1} = T_1\bar M_n+\dots+T_d\bar M_n.  
$$
For such an $n$, Remark 3.7 implies that 
$$
M_{n+1}=M_n + x_1M_n+\dots+x_d M_n,
$$
hence $\{M_n\}$ is proper.  \qedd\enddemo

\proclaim{Lemma 3.9}
Let $\{M_n\}$ be a proper filtration of a 
$k[x_1,\dots,x_d]$-module $M$, let $K\subseteq M$ be a 
submodule, and let $\{K_n\}$ be the filtration induced 
on $K$ by 
$$
K_n = K\cap M_n.  
$$
Then $\{K_n\}$ is a proper filtration of $K$.  
\endproclaim
\demo{proof}
Form the graded modules 
$$
\bar M = \sum_{n\in \Bbb Z} M_n/M_{n-1}
$$
and 
$$
\bar K = \sum_{n\in \Bbb Z} K_n/K_{n-1} .  
$$
Because of the natural isomorphism
$$
\bar K_n=
K\cap M_n/K\cap M_{n-1} \cong( K\cap M_n + M_{n-1})/M_{n-1}
\subseteq M_n/M_{n-1}=\bar M_n, 
$$
$\bar K$ is isomorphic to a submodule of $\bar M$.  
Lemma 3.8 implies that $\bar M$ is finitely generated.  Thus 
by Hilbert's basis theorem (asserting that graded submodules of 
finitely generated graded modules are finitely generated), 
it follows that $\bar K$ is finitely generated.  Now apply 
Lemma 3.8 once again to conclude that $\{K_n\}$ is a 
proper filtration of $K$.  \qedd
\enddemo
We remark that the proof of Lemma 3.9 is inspired by Cartier's 
proof of the Artin-Rees lemma \cite{29}, p II-9.  

Let $M$ be a finitely generated $A$-module, 
choose a finite dimensional subspace $G\subseteq M$ 
which generates $M$ as an $A$-module, and set
$$
M_n=\sp \{f\cdot\zeta: 
f\in A,\quad \deg f\leq n,\quad \zeta\in G\}.
$$
Since $\{M_n\}$ is a proper filtration, Proposition 
3.7 implies that the number 
$$
c(M)=d!\lim_{n\to\infty}\frac{\dim M_n}{n^d}
$$
exists as an invariant of $M$ independently of 
the choice of generator $G$.  
The following result shows that this 
invariant is additive 
on short exact sequences.  

\proclaim{Proposition 3.10}
For every exact sequence 
$$
0\longrightarrow K\longrightarrow L\longrightarrow M
\longrightarrow 0
$$
of finitely generated $k[x_1,\dots,x_d]$-modules we have 
$c(L) = c(K) + c(M)$.
\endproclaim
\demo{proof}
Since $c(M)$ depends only on the isomorphism class of 
$M$, we may assume that $K\subseteq L$ is a submodule 
of $L$ and $M=L/K$ is its quotient.  Pick any proper 
filtration $\{L_n\}$ for $L$ and let $\{\dot L_n\}$ 
and $\{K_n\}$ be the associated filtrations of 
$L/K$ and $K$ 
$$
\align
\dot L_n &= (L_n + K)/K \subseteq L/K, \\
K_n &= K\cap M_n \subseteq K.  
\endalign
$$
It is obvious that $\{\dot L_n\}$ is proper, and 
Lemma 3.9 implies that $\{K_n\}$ is proper as well.  

Now for each $n\geq 1$ we have an exact sequence of 
finite dimensional vector spaces
$$
0\longrightarrow K_n\longrightarrow L_n 
\longrightarrow \dot L_n \longrightarrow 0,
$$
and hence 
$$
\dim L_n = \dim K_n + \dim \dot L_n.  
$$
Since each of the three filtrations is proper we can 
multiply the preceding equation through by $d!/n^d$ and 
take the limit to obtain
$c(L) = c(K) + c(L/K)$.
\qedd
\enddemo

\remark{Remark 3.11}
The addition property of Proposition 3.10 generalizes 
immediately to the following assertion.  For every finite
exact sequence 
$$
0\longrightarrow M_n\longrightarrow \dots
\longrightarrow M_1\longrightarrow M_0\longrightarrow 0
$$
of finitely generated $k[x_1,\dots,x_d]$-modules, 
we have 
$$
\sum_{k=0}^n (-1)^kc(M_k) = 0.  
$$
\endremark

\proclaim{Corollary}
Let $M$ be a finitely generated $k[x_1,\dots,x_d]$-module
and let 
$$
0\longrightarrow F_n\longrightarrow \dots 
\longrightarrow F_1\longrightarrow M\longrightarrow 0
$$
be a finite free resolution of $M$, where 
$$
F_k = \beta_k\cdot k[x_1,\dots,x_d]
$$
is a direct sum of $\beta_k$ copies of the rank-one 
free module $k[x_1,\dots,x_d]$.  Then
$$
c(M) = \sum_{k=1}^n(-1)^{k+1}\beta_k.  
$$
\endproclaim
\demo{proof}
Remark (3.11) implies that 
$$
c(M) = \sum_{k=1}^n(-1)^{k+1}c(F_k), 
$$
and thus it suffices to show that if 
$F=\beta \cdot k[x_1,\dots,x_d]$ 
is a free module of rank $\beta\in \Bbb N$, then 
$c(F)=\beta$.  

By the additivity property of 3.10 we have 
$$
c(\beta\cdot k[x_1,\dots,x_d]) = 
\beta \cdot c(k[x_1,\dots,x_d])
$$
and thus we have to show that $c(k[x_1,\dots,x_d])$ is $1$.  

This follows from a computation of the dimensions of 
$$
\Cal P_n = \{f\in k[x_1,\dots,x_d]: \deg f\leq n\}
$$
and it is a classical result that 
$$
\dim \Cal P_n = q_d(n) = \frac{(n+1)\dots (n+d)}{d!}
$$
(see, for example, Appendix A of \cite{1}).  Thus
$$
c(k[x_1,\dots,x_d]) = 
d!\lim_{n\to\infty}\frac{\dim \Cal P_n}{n^d}=
\lim_{n\to\infty}\frac{(n+1)\dots(n+d)}{n^d}=1
$$
and the corollary is established.\qedd
\enddemo

We now deduce the main result of this section.  Let $H$ be 
a finite-rank Hilbert module over $A=\Bbb C[z_1,\dots,z_d]$, 
and let $\phi: \Cal B(H)\to \Cal B(H)$ be its associated 
completely positive map
$\phi(A) = T_1AT_1^*+\dots+T_dAT_d^*$.  

\proclaim{Theorem C}
$$
\chi(H) = 
d!\lim_{n\to\infty}\frac{\rank\,(\bold 1-\phi^{n+1}(\bold 1))}{n^d}.  
$$
\endproclaim
 \demo{proof}
Consider the module 
$$
M_H=\sp \{f\cdot \Delta\xi: f\in A, \quad \xi\in H\}
$$
and its natural (proper) filtration 
$$
M_n = \sp \{f\cdot \Delta\xi: \deg f\leq n, \quad \xi\in H\},
\qquad n=1,2,\dots.  
$$
In view of the definition of $\chi(H)$ in terms of free 
resolutions of $M_H$, the preceding corollary implies that
$$
\chi(H) = c(M_H) = d!\lim_{n\to\infty}\frac{\dim M_n}{n^d}.  
$$
Thus it suffices to show that 
$$
\dim M_n = \rank (\bold 1-\phi^{n+1}(\bold 1))
$$
for every $n=1,2,\dots$.  For that, we will prove
$$
M_n = (\bold 1-\phi^{n+1}(\bold 1))H.  
\tag{3.12}
$$

Indeed, writing 
$$
\bold 1 - \phi^{n+1}(\bold 1) = 
\sum_{k=0}^n\phi^k(\bold 1-\phi(\bold 1)) =
\sum_{k=0}^n\phi^k(\Delta^2),\tag{3.13}
$$
we see in particular that 
$\bold 1-\phi^{n+1}(\bold 1)$ is a positive 
finite rank operator for every $n$ and hence
$$
(\bold 1-\phi^{n+1}(\bold 1))H = 
\ker (\bold 1-\phi^{n+1}(\bold 1))^\perp.  
$$
The kernel of $\bold 1-\phi^{n+1}(\bold 1)$ 
is easily computed.  We have 
$$
\ker(\bold 1-\phi^{n+1}(\bold 1)) = 
\{\xi\in H: \<(\bold 1-\phi^{n+1}(\bold 1))\xi,\xi\> = 0\}, 
$$
and by (3.13), $\<(\bold 1-\phi^{n+1}(\bold 1))\xi,\xi\> = 0$ 
iff 
$$
\sum_{k=0}^n\<\phi^k(\Delta^2)\xi,\xi\> = 0.  
$$
Since 
$$
\phi^k(\Delta^2) = \sum_{i_1,\dots,i_k=1}^d T_{i_1}\dots T_{i_k}
\Delta^2T_{i_k}^*\dots T_{i_1}^*, 
$$
the latter is equivalent to 
$$
\sum_{k=0}^n\sum_{i_1,\dots,i_k=1}^d\|\Delta
T_{i_k}^*\dots T_{i_1}^*\xi\|^2 = 0.  
$$
Thus the kernel of $\bold 1-\phi^{n+1}(\bold 1)$
is the orthocomplement of the space spanned by 
$$
\{T_{i_1}\dots T_{i_k}\Delta\eta: \eta\in H, 
\quad 1\leq i_1,\dots,i_k\leq d, \quad k=0,1,\dots,n\},   
$$ 
namely $M_n=\sp\{f\cdot\Delta\eta: \deg f\leq n,\quad
\eta\in H\}$.  
This shows that 
$$
\ker (\bold 1-\phi^{n+1}(\bold 1)) = M_n^\perp, 
$$
from which formula (3.12) is evident.  
\qedd
\enddemo

\remark{Remark 3.14}
Closed submodules of finite-rank Hilbert $A$-modules 
neet not have finite rank (see section 8, Corollary
of Theorem F).  However, 
if $H_0$ is a submodule of a finite rank Hilbert
module $H$ which is of finite codimension in $H$, then
$\rank(H_0)<\infty$.  Indeed, if $P_0$ is the projection
of $H$ onto $H_0$, then 
$$
\rank(H_0)=\rank(\bold 1_{H_0}-\phi_{H_0}(\bold 1_{H_0}))=
\rank(P_0-\phi_H(P_0)).  
$$
Since $P_0-\phi_H(P_0) = (\bold 1_H-\phi_H(\bold 1_H))
-P_0^\perp + \phi_H(P_0^\perp)$, we have 
$$
\rank(H_0)\leq \rank(H) + \rank(P_0^\perp)+
\rank(\phi_H(P_0^\perp))<\infty.  
$$

On the other hand, given a submodule $H_0\subseteq H$
with $\dim(H/H_0)<\infty$, the algebraic module 
$M_{H_0}$ is not a submodule of $M_H$, nor is it 
conveniently related to $M_H$.  Thus there is 
no direct way of relating 
$\chi(H)$ to $\chi(H_0)$ by way of their 
definitions.  Nevertheless, 
Theorem C implies the following.  
\endremark

\proclaim{Corollary 1: stability of Euler characteristic}
Let $H_0$ be a closed submodule of a finite rank Hilbert 
$A$-module $H$ such that $\dim(H/H_0)<\infty$.  Then 
$\chi(H_0)=\chi(H)$.  
\endproclaim
\demo{proof}
By estimating as in Remark 3.14 we have
$$
\align
\rank(\bold 1_{H_0}-\phi_{H_0}^{n+1}(\bold 1_{H_0}))\leq 
&\rank(\bold 1_H -\phi_H^{n+1}(\bold 1_H)) +\\
&\rank P_0^\perp + \rank(\phi_H^{n+1}(P_0^\perp)), 
\endalign
$$
$P_0$ denoting the projection of $H$ on $H_0$.  
Similarly,
$$
\align
\rank(\bold 1_H-\phi_H^{n+1}(\bold 1_H))\leq 
&\rank(\bold 1_{H_0} -\phi_{H_0}^{n+1}(\bold 1_{H_0})) +\\
&\rank P_0^\perp + \rank(\phi_H^{n+1}(P_0^\perp)), 
\endalign
$$
Thus we have the inequality
$$
|\rank(\bold 1_H-\phi_H^{n+1}(\bold 1_H))-
\rank(\bold 1_{H_0}-\phi_{H_0}^{n+1}(\bold 1_{H_0}))| 
\leq\rank P_0^\perp + \rank(\phi_H^{n+1}(P_0^\perp)). 
\tag{3.15}
$$

One estimates the right side as follows.  Note that
$$
\<\phi^{n+1}(P_0^\perp)\xi,\xi\> = 
\sum_{i_1,\dots,i_{n+1}=1}^d\<P_0^\perp 
T_{i_{n+1}}^*\dots T_1^*\xi,T_{i_{n+1}}^*\dots T_1^*\xi\>
$$
vanishes iff $\xi$ belongs to the kernel of every operator
of the form $P_0^\perp f(T_1,\dots,T_d)^*$ where 
$f\in E_{n+1}H^2$ is a homogeneous polynomial of 
degree $n+1$.  Hence the range of the positive finite
rank operator $\phi^{n+1}(P_0^\perp)$ is the orthocomplement
of all such vectors $\xi$, and is therefore spanned 
linearly by the ranges of all operators 
$f(T_1,\dots,T_d)P_0^\perp$, $f\in E_{n+1}H^2$, i.e., 
$$
\sp\{f\cdot \zeta: f\in E_{n+1}H^2,
\quad \zeta\in P_0^\perp H\}.  
$$
It follows that 
$$
\rank(\phi^{n+1}(P_0^\perp))\leq 
\dim(E_{n+1}H^2)\cdot \rank P_0^\perp =
q_{d-1}(n+1) \rank P_0^\perp.  
$$
Thus (3.15) implies that 
$$
|\frac{\rank(\bold 1_H-\phi_H^{n+1}(\bold 1_H))}{n^d} -
\frac{\rank(\bold 1_{H_0}-\phi_{H_0}^{n+1}(\bold 1_{H_0}))}{n^d}|
$$
is at most
$$
\rank(P_0^\perp)\frac{1+q_{d-1}(n+1)}{n^d}.  
$$
Since $q_{d-1}(x)$ is a polynomial of degree $d-1$, the 
latter tends to zero as $n\to\infty$, and the conclusion 
$|\chi(H)-\chi(H_0)|=0$ follows 
from Theorem C after taking the limit on $n$.  
\qedd
\enddemo

For algebraic reasons, the Euler characteristic of a 
finitely generated $A$-module must be nonnegative
(\cite{18}, Theorem 192) and hence $\chi(H)\geq 0$
for every finite rank Hilbert $A$-module $H$.  
One also has the following 
upper bound, which we collect here for later use.

\proclaim{Corollary 2}
For every finite rank Hilbert $A$-module $H$, 
$0\leq \chi(H)\leq \rank(H)$.  
\endproclaim
\demo{proof}
Let $M_H$ be the algebraic module associated with 
$H$ and let $M_1\subseteq M_2\subseteq \dots$ be 
the proper filtration of it defined by 
$$
M_n=\sp \{f\cdot\xi: 
f\in A,\quad \deg f\leq n,\quad \xi\in\Delta H\}
$$
$\Delta$ denoting the square root of 
$\bold 1_H-T_1T_1^*-\dots-T_dT_d^*$.  Clearly 
$$
\dim M_n\leq 
\dim \{f\in A: \deg f\leq n\}\cdot\dim \Delta H =
q_d(N)\cdot \rank(H).  
$$
From the corollary of Prop 3.10 which identifies $\chi(M_H)$
with 
$$
c(M_H)=d!\lim_{n\to\infty}\frac{\dim M_n}{n^d} \leq 
d!\lim_{n\to\infty}\frac{q_d(n)}{n^d}\cdot \rank(H)=
\rank(H),
$$
the inequality follows.  
\qedd
\enddemo

\subheading{4. Curvature invariant}
Let $H$ be a Hilbert $A$-module with canonical operators 
$T_1,\dots,T_d$.  For every $z\in \Bbb C^d$ we define 
the operator $T(z)\in\Cal B(H)$ as in (0.2),
$$
T(z)=\bar z_1T_1+\dots+\bar z_dT_d.  
$$
We have pointed out in the introduction that 
$\|T(z)\|\leq |z|$, and hence $\bold 1 -T(z)$ is 
invertible for all $z$ in the open unit ball 
$B_d$.  Thus we can define an operator-valued function 
$F:B_d\to \Cal B(\overline{\Delta H})$ as follows:
$$
F(z)\xi = 
\Delta (\bold 1-T(z)^*)^{-1}(\bold 1-T(z))^{-1}\Delta\xi, 
\qquad \xi\in\overline{\Delta H}.    
\tag{4.1}
$$

Assuming that $\rank(H)<\infty$, then $F(z)$ is a positive
operator acting on a finite dimensional Hilbert space
and we may take its trace to obtain 
a numerical function defined for all $z$ the open ball.  
We show in 
Theorem A below that the latter function has ``renormalized" 
boundary values 
$$
K_0(z)=\lim_{r\uparrow 1}(1-r^2){\text{trace}}\,F(rz)
$$
for almost every point $z\in\partial B_d$ relative to 
the natural measure $d\sigma$ on $\partial B_d$.  
The curvature invariant $K(H)$ is defined by integrating 
the function $K_0(\cdot)$ over $\partial B_d$, and 
in general $K(H)$ is a real number  satisfying
$0\leq K(H)\leq \rank(H)$.  Significantly, 
the curvature invariant is sufficiently 
sensitive that it detects
precisely when the closed submodule of $H$ generated by 
$\Delta H$ is a free Hilbert module;
the criterion for freeness is that the curvature
should be the {\it maximum\/} possible
value $K(H)=\rank(H)$.  

Fix a Hilbert $A$-module $H$, of arbitrary positive 
rank, and form the free Hilbert $A$-module 
$H^2\otimes \overline{\Delta H}$ of rank 
$r=\rank H$.  We have seen in section 1 that 
there is a unique operator 
$U_0\in\hom(H^2\otimes \overline{\Delta H},H)$ satisfying 
$$
U_0(f\otimes \zeta) = f\cdot\Delta\zeta,\qquad 
f\in A,\quad \zeta\in H.  
\tag{4.2}
$$
$U_0$ is a coisometry {\it only} when $H$ is pure; but 
it is a contraction in general and hence
$\bold 1-U_0^*U_0$ is a positive operator in 
$\Cal B(H^2\otimes\overline{\Delta H})$ of norm at 
most $1$.  The following result implies that 
$\bold 1-U_0^*U_0$ can be associated with a 
multiplier, and that fact 
is essential to the proof of Theorems A and D.  

\proclaim{Proposition 4.3}
There is a free Hilbert module $F=H^2\otimes E$ and a
multiplier $\Phi\in\Cal M(E,\overline{\Delta H})$ whose
associated homomorphism satisfies 
$$
U_0^*U_0 + \Phi\Phi^* = 
\bold 1_{H^2\otimes\overline{\Delta H}}.
$$
\endproclaim
\demo{proof}
Let $F_0$ be the free module $H^2\otimes \overline{\Delta H}$. 
By Theorem 1.14 there is a spherical module $S_0$ and 
a map $U_1\in\hom(S_0,H)$ such that 
$$
U: (\xi,\eta)\in F_0\oplus S_0\mapsto U_0\xi + U_1\eta\in H
$$
defines a minimal dilation 
$$
F_0\oplus S_0\underset U\to\longrightarrow H\longrightarrow 0.  
$$  

We consider the kernel $K=\ker U$ of the dilation 
map $U$.  $K$ is a Hilbert 
submodule of $F_0\oplus S_0$ and therefore it too has 
a minimal dilation 
$$
F\oplus S\underset V\to\longrightarrow K\longrightarrow 0,
$$
where $V\in\hom(F\oplus S,F_0\oplus S_0)$ satisfies 
$VV^*=P_K$, $P_K$ denoting the projection of 
$F_0\oplus S_0$ onto $K$.  Define 
$\Phi\in\hom(F,F_0)$ by 
$\Phi = P_{F_0}V\restriction_F$, 
$P_{F_0}$ denoting the projection of $F_0\oplus S_0$ 
onto the first summand.  We have to show that 
$$
U_0^*U_0 = \bold 1_{F_0} - \Phi\Phi^*,
\tag{4.4}
$$
and for that we require
\proclaim{Lemma 4.5}
Let $F_1,F_2$ be free Hilbert $A$-modules, let 
$S_1,S_2$ be spherical Hilbert $A$-modules, and let 
$V\in\hom(F_1\oplus S_1,F_2\oplus S_2)$.  Then 
$VS_1\subseteq S_2$.  
\endproclaim

Assuming for the moment that Lemma 4.5 has been proved,
we establish (4.4) as follows.  Since $U_0$ is the restriction 
of $U$ to the free summand $F_0\subseteq F_0\oplus S_0$ we 
have $U_0^*U_0 = P_{F_0}U^*U\restriction_{F_0}$.  
Since $U$ is a coisometry, $U^*U$ is the projection 
$\bold 1 - P_K$ onto the orthocomplement of 
$K=\ker U \subseteq F_0\oplus S_0$, hence
$$
U_0^*U_0 = P_{F_0}(\bold 1-P_K)\restriction_{F_0} =
\bold 1_{F_0}-P_{F_0}P_K\restriction_{F_0}.  
\tag{4.6}
$$
Applying Lemma 4.5 to the map 
$V\in\hom(F\oplus S,F_0\oplus S_0)$ we find that 
$VS\subseteq S_0$ or, equivalently, that 
$P_{F_0}V = P_{F_0}VP_F$.  Since $VV^*=P_K$, 
the right side of (4.6) becomes 
$$
\bold 1_{F_0}-P_{F_0}VV^*\restriction_{F_0} =
\bold 1_{F_0}-P_{F_0}VP_FV^*\restriction_{F_0} =
\bold 1_{F_0}-\Phi\Phi^*,
$$
as required.  

\demo{proof of Lemma 4.5}
Let $T_1,\dots,T_d$ be the canonical operators associated
with the Hilbert module $F_1\oplus S_1$, and 
let $\phi_1: \Cal B(F_1\oplus S_1)\to\Cal B(F_1\oplus S_1)$
be their associated completely positive map.  Then 
$\phi_1(\bold 1)\geq \phi_1^2(\bold 1)
\geq\phi_1^3(\bold 1)\geq\dots$ is 
a decreasing sequence of projections with limit projection
$$
\lim_{n\to\infty}\phi_1^n(\bold 1)=0\oplus \bold 1_{S_1}.  
$$
Similarly, if $\phi_2$ is the corresponding completely 
positive map on $\Cal B(F_2\oplus S_2)$ then 
$\phi_2^n(\bold 1) \downarrow 0\oplus \bold 1_{S_2}$.  
Since $V\in\hom(F_1\oplus S_1,F_2\oplus S_2)$ we have 
$$
V\phi_1^n(\bold 1)V^* = \phi_2^n(VV^*)\leq 
\phi_2^n(\bold 1) 
$$
for every $n$, hence 
$$
V(0\oplus \bold 1_{S_1})V^* = 
\lim_{n\to\infty}V\phi_1^n(\bold 1)V^* \leq 
\lim_{n\to\infty}\phi_2^n(\bold 1) = 0\oplus \bold 1_{S_2},   
$$
from which  $VS_1\subseteq S_2$ follows.  
\qedd\enddemo
\enddemo

From Proposition 4.3 we obtain
\proclaim{Corollary}
Let $F:B_d\to \Cal B(\overline{\Delta H})$ be the function 
(4.1) and let $\Phi\in\Cal M(E,\overline{\Delta H})$ be the 
multiplier of Proposition 4.3.  Then for all $z\in B_d$ 
we have
$$
(1-|z|^2)F(z) = \bold 1 - \Phi(z)\Phi(z)^*.  
$$
\endproclaim
\demo{proof}
Fix $\alpha\in B_d$, $\zeta_1,\zeta_2\in \overline{\Delta H}$.  
From (4.1) we can write
$$
\<F(\alpha)\zeta_1,\zeta_2\> = 
\<(\bold 1-T(\alpha))^{-1}\Delta\zeta_1,
(\bold 1-T(\alpha))^{-1}\Delta\zeta_2\>.
\tag{4.7}
$$
Consider the operator $U_0: H^2\otimes\overline{\Delta H}\to H$
given by $U_0(f\otimes\zeta)=f\cdot \Delta\zeta$.  Notice that 
for the element $u_\alpha\in H^2$ defined by 
$$
u_\alpha(z) = (1-\<z,\alpha\>)^{-1}, \qquad z\in B_d
$$
we have 
$$
U_0(u_\alpha\otimes \zeta)=(\bold 1-T(\alpha))^{-1}\Delta \zeta.  
\tag{4.8}
$$
Indeed, the sequence of polynomials $f_n\in H^2$ defined by
$$
f_n(z) = \sum_{k=0}^n\<z,\alpha\>^k
$$
converges in the $H^2$-norm to $u_\alpha$ since
$$
\|u_\alpha-f_n\|^2=\sum_{k=n+1}^\infty |\alpha|^{2k}\to 0
$$
as $n\to\infty$.  Since
$$
U_0(f_n\otimes\zeta) = f_n\cdot\Delta\zeta = 
\sum_{k=0}^nT(\alpha)^k\Delta\zeta, 
$$
formula (4.8) follows by taking the limit as $n\to\infty$.  

From (4.8) we find that 
$$
\<F(\alpha)\zeta_1,\zeta_2\> = 
\<U_0(u_\alpha\otimes\zeta_1),U_0(u_\alpha\otimes\zeta_2)\>=
\<U_0^*U_0 u_\alpha\otimes \zeta_1, u_\alpha\otimes\zeta_2\>.  
$$
By Proposition 4.3 we have $U_0^*U_0 = \bold 1-\Phi\Phi^*$, 
and using the formula 
$\Phi^*(u_\alpha\otimes\zeta)=
u_\alpha\otimes\Phi(\alpha)^*\zeta$
of Proposition 2.4 we can write
$$
\align
\<F(\alpha)\zeta_1,\zeta_2\> &= 
\<(\bold 1-\Phi\Phi^*)u_\alpha\otimes\zeta_1,
u_\alpha\otimes\zeta_2\> \\&=
\<u_\alpha\otimes\zeta_1,u_\alpha\otimes\zeta_2\>-
\<u_\alpha\otimes\Phi(\alpha)^*\zeta_1,
u_\alpha\otimes\Phi(\alpha)^*\zeta_2\>\\ &=
\|u_\alpha\|^2(\<\zeta_1,\zeta_2\>-\<\Phi(\alpha)^*\zeta_1,
\Phi(\alpha)^*\zeta_2\>) \\&=
(1-|\alpha|^2)^{-1}
\<(\bold 1-\Phi(\alpha)\Phi(\alpha)^*\zeta_1,\zeta_2\>, 
\endalign
$$
and the corollary follows after multiplying through 
by $1-|\alpha|^2$.
\qedd
\enddemo

\proclaim{Lemma 4.9}
Let $E_1,E_2$ be Hilbert spaces with $E_2$ finite dimensional, 
let $\Phi\in\Cal M(E_1,E_2)$ be a multiplier, let 
$\tilde\Phi: \partial B_d\to \Cal B(E_1,E_2)$ be its 
boundary function (see Proposition 2.8), and let 
$\sigma$ be normalized measure on $\partial B_d$.  Then for 
$\sigma$-almost every point $z\in\partial B_d$ we have 
$$
\lim_{r\uparrow 1}{\text{trace}}\,(\Phi(rz)\Phi(rz)^*) =
{\text{trace}}\,(\tilde\Phi(z)\tilde\Phi(z)^*).  
$$
\endproclaim 
\demo{proof}
For any Hibert-Schmidt operators $A,B\in\Cal B(E_1,E_2)$ we 
have ${\text{trace}}\,AA^* = {\text{trace}}\,A^*A$ and 
$$
|\sqrt{{\text{trace}}(AA^*)}-\sqrt{{\text{trace}}\,(BB^*)}|\leq 
({\text{trace}}\,(|A-B|^2))^{1/2}
$$
(where we have written the usual $|X|^2$ for $X^*X$).  Thus it
suffices to show that 
$$
\lim_{r\uparrow 1}{\text{trace}}\,|\Phi(rz)-\tilde\Phi(z)|^2 = 0
$$
almost everywhere ($d\sigma$).  

Now since $E_2$ is finite-dimensional, the operator norm on 
$\Cal B(E_1,E_2)$ is equivalent to the Hilbert-Schmidt norm 
and since the space $\Cal L^2(E_1,E_2)$ of Hilbert-Schmidt 
operators is a Hilbert space with its natural norm 
$$
\|A\|_2 = ({\text{trace}}\,A^*A)^{1/2}, 
$$
we may consider that the multiplier $\Phi$ is a bounded 
holomorphic function from $B_d$ to the Hilbert space 
$\Cal L^2(E_1,E_2)$.  The required assertion now follows 
from (2.6).  
\qedd
\enddemo

\proclaim{Theorem A}
Let $H$ be a Hilbert $A$-module of finite positive rank, 
let $F: B_d\to \Cal B(\Delta H)$ be the operator 
function defined by (4.1), and let $\sigma$ denote normalized 
measure on the sphere $\partial B_d$.  Then for $\sigma$-almost
every $z\in B_d$, the limit 
$$
K_0(z) = \lim_{r\uparrow 1}(1-r^2){\text{trace}}\,F(r\cdot z)
$$
exists and satisfies
$$
0\leq K_0(z)\leq \bold \rank \,H.  
$$

Moreover, the extreme case 
$K_0(z)=\rank\,H$ (a.e.) occurs if, and 
only if, the closed submodule of $H$ generated by 
the range of $\Delta$ is free.  
\endproclaim
\demo{proof}
The corollary of Proposition 4.3 implies that for all 
$z\in\partial B_d$ and every $r\in(0,1)$ we have 
$$
(1-r^2)F(rz) = \bold 1-\Phi(rz)\Phi(rz)^*.  
$$
From Lemma 4.9 we may conclude that for almost every 
$z\in\partial B_d$, 
$$
\lim_{r\uparrow 1}(1-r^2){\text{trace}}\,F(rz) =
\rank H -{\text{trace}}\,\tilde\Phi(z)\tilde\Phi(z)^*,
$$
and hence the limit function $K_0(\cdot)$ is expressed 
in terms of $\tilde\Phi$ by 
$$
K_0(z) = 
{\text{trace}}\,(\bold 1_{\Delta H} - 
\tilde\Phi(z)\tilde\Phi(z)^*).  
\tag{4.10}
$$
Since $\|\tilde\Phi(z)\|\leq 1$ we have 
$0\leq \bold 1_{\Delta H}-\tilde\Phi(z)\tilde\Phi(z)^*
\leq \bold 1_{\Delta H}$ and hence 
$$
0\leq K_0(z)\leq 
\tr (\bold 1_{\Delta H}-\tilde\Phi(z)\tilde\Phi(z)^*)\leq
\tr \bold 1_{\Delta H} = \rank (H)
$$
for almost every $z\in\partial B_d$.  

Now consider the extreme case in which 
$$
K_0(z) = \rank\,H, \qquad {\text{a.e. }}d\sigma(z).  
\tag{4.11}
$$
From formula (4.10) it 
follows that ${\text{trace}}\,(\tilde\Phi(z)\tilde\Phi(z)^*) = 0$ 
for almost every $z\in \partial B_d$.  
Since the trace is faithful and since
$\Phi$ is uniquely determined by its boundary values, 
we conclude that $\Phi=0$. 
Proposition 4.3 implies that $U_0$ must be an 
isometry.  Thus $U_0$ is an isomorphism of 
$H^2\otimes\overline{\Delta H}$ onto the closed 
submodule of $H$ 
$$
\overline{M_H}=\overline{\text{span}}\{f\cdot\Delta\xi: 
f\in A,\quad \xi\in H\}
$$
generated by $\Delta H$.  

Conversely, if $\overline{M_H}$ is isomorphic to 
a free Hilbert module
$H^2\otimes E$ then by a direct computation one finds
that $\Delta$ is identified with the projection onto 
$1\otimes E\subseteq H^2\otimes E$,
hence $\Delta H$ is identified with $E$ and 
$U_0$ is clearly an isometry. By Proposition 4.3, 
the multiplier $\Phi$ must be $0$, and formula (4.10)
shows that $K_0(z)$ is the constant function with 
value $\tr \bold 1_{\Delta H}=\rank\,H$.  \qedd 
\enddemo

We now define the curvature invariant of a Hilbert 
$A$-module $H$ of finite rank,
$$
K(H) = \int_{\partial B_d} K_0(z)\, d\sigma(z).  
\tag{4.12}
$$
The basic property of $K(H)$ is that it is sensitive 
enough to detect exactly when a finite rank pure 
Hilbert module is free.  

\proclaim{Theorem 4.13}
For every finite rank Hilbert $A$-module $H$, we have 
$$
0\leq K(H)\leq \rank(H).  
$$
If $H$ is pure, then $K(H)=\rank(H)$ if, and only if, 
$H\cong H^2\otimes E$ is a free Hilbert $A$-module.  
\endproclaim
\demo{proof}
This is immediate from the preceding discussion after noting
that for a pure Hilbert $A$-module $H$, the map
$U_0: H^2\otimes \overline{\Delta H}\to H$ of (4.2) 
is a coisometry, and in particular 
$$
H=U_0(H^2\otimes\overline{\Delta H})=\overline{M_H} =
\overline{\text{span}}\{f\cdot\Delta \xi: f\in A,\quad\xi\in H\}
$$
is generated as a Hilbert $A$-module by the range of $\Delta$.  
\qedd
\enddemo

The curvature invariant also detects ``inner sequences".  
More precisely, let $M$ be a closed submodule of the rank-one
free Hilbert $A$-module $H^2$.  From the corollary of 
Theorem 5.9 below (Theorem 5.9 is proved
independently of the discussion
to follow) there is a Hilbert space $E$ and a multiplier 
$\Phi\in\Cal M(E,\Bbb C)$ whose morphism satisfies 
$\Phi\Phi^*=P_M$.  Choosing an orthonormal basis 
$e_1,e_2,\dots$ for $E$ we obtain 
a sequence $\{\phi_n\}$ of 
multipliers of $H^2$ as follows
$$
\phi_n(z)=(\Phi e_n)(z)=\Phi(z)e_n, 
\qquad n=1,2,\dots.  
$$
The associated multiplication operators 
$V_nf=\phi_n\cdot f$ satisfy
$$
V_1V_1^*+V_2V_2^*+\dots = \Phi\Phi^*=P_M. 
\tag{4.14}
$$
For definiteness of notation, we can assume that the 
sequence $\phi_1,\phi_2,\dots$ is infinite by adding 
harmless zero functions if it is not. 

Note that 
$$
\sup_{|z|<1} \sum_{n=1}^\infty |\phi_n(z)|^2\leq 1.
\tag{4.15}
$$ 
Indeed, if $\{u_\alpha: \alpha\in B_d\}$ denotes the 
family of functions in $H^2$ defined in (1.1), then 
$v_\alpha = (1-|\alpha|^2)^{1/2}u_\alpha$ is a unit 
vector in $H^2$ which is an eigenvector for the adjoints
of multiplication operators associated with 
any multiplier; thus for the operators $V_n$ we have
$V_n^*v_\alpha = \overline{\phi_\alpha(z)}v_\alpha$.  
Using (4.14) we find that 
$$
\sum_{n=1}^\infty |\phi_n(\alpha)|^2 = 
\sum_{n=1}^\infty \|V_n^* v_\alpha\|^2=
\sum_{n=1}^\infty \<V_nV_n^*v_\alpha,v_\alpha\>=
\<P_M v_\alpha,v_\alpha\>\leq 1,
$$
and (4.15) follows.  

Therefore, the boundary functions 
$\tilde\phi_n:\partial B_d\to\Bbb C$ must 
satisfy
$$
\sum_{n=1}^\infty |\tilde \phi_n(z)|^2\leq 1
$$
almost everywhere with respect to the natural normalized
measure $\sigma$ on $\partial B_d$ and we will say that 
$(\phi_n)$ is an {\it inner sequence} 
if equality holds almost everywhere 
$$
\sum_{n=1}^\infty |\tilde \phi_n(z)|^2=1,
\qquad \text{a. e. }(d\sigma) \text{ on }\partial B_d.  
$$
Significantly, we do not know if every sequence 
$(\phi_n)$ associated with an invariant subspace of 
$H^2$ as in (4.14) must be an inner sequence (see
section 8).  The following result shows the relevance
of the curvature invariant for this problem.  

\proclaim{Theorem 4.16}
Let $M\subseteq H^2$ be a closed submodule and let 
$\phi_1,\phi_2,\dots$ be a sequence in the multiplier 
algebra of $H^2$ whose multiplication operators
satisfy
$$
P_M=M_{\phi_1}M_{\phi_1}^*+M_{\phi_2}M_{\phi_2}^*+\dots.
$$
Then $(\phi_n)$ is an inner sequence 
iff $K(H^2/M)=0$.  
\endproclaim
\demo{proof}
Let $F$ be the direct sum of an infinite number of 
copies of $H^2$ and define $\Phi\in\hom(F,H^2)$ 
by 
$$
\Phi(f_1,f_2,\dots)=\sum_{n=1}^\infty \phi_n\cdot f_n.
$$
Then $\Phi\Phi^*=P_M$.  Letting 
$U: H^2\to H^2/M$ be the natural projection, 
then $U$ defines a minimal dilation of $H^2/M$
and we have 
$$
U^*U + \Phi\Phi^*=\bold 1_{H^2}.  
$$
Writing 
$$
K(H^2/M) = \int_{\partial B_d} K_0(z)\,d\sigma(z),
$$
we see from formula (4.10) that in this case 
$$
K_0(z) = 1-\tilde\Phi(z)\tilde\Phi(z)^* =
1-\sum_{n=1}^\infty |\tilde\phi_n(z)|^2, 
$$
and therefore $K(H^2/M)=0$ iff 
$\sum_n|\tilde\phi_n(z)|^2=1$ almost everywhere 
on $\partial B_d$.
\qedd
\enddemo

We will discuss applications of Theorem 4.16 in 
section 8.

\subheading{5. Curvature operator: quantizing the Gauss map}

Let us recall a convenient description of the Gaussian curvature
of a compact oriented Riemannian $2$-manifold $M$.  It is not 
necessary to do so, but for simplicity we will assume that 
$M\subseteq \Bbb R^3$ can be embedded in $\Bbb R^3$ 
in such a way that it 
inherits the usual metric structure of $\Bbb R^3$.  After
choosing one of the two orientations of $M$ (as a nondegenerate
$2$-form) we normalize it in the obvious way to obtain a 
continuous field of {\it unit} normal vectors at every 
point of $M$.  

For every point $p$ of $M$ one can translate the normal vector 
at $p$ to the origin of $\Bbb R^3$ (without changing its
direction), and the endpoint of that translated vector is a 
point $\gamma(p)$ on the unit sphere $S^2$.  This defines
the Gauss map 
$$
\gamma: M\to S^2
\tag{5.1}
$$
of $M$ to the sphere.  Now fix $p\in M$.  The tangent 
plane $T_pM$ is obviously parallel to the corresponding 
tangent plane $T_{\gamma(p)}S^2$ of the sphere (they 
have the same normal vector)
and hence both are cosets
of the same $2$-dimensional subspace $V\subseteq \Bbb R^3$:
$$
T_pM = p + V, \qquad T_{\gamma(p)}S^2 = \gamma(p) + V.  
$$
Thus the differential $d\gamma(p)$ defines a linear 
operator on the two-dimensional vector space $V$, and 
the Gaussian curvature $K(p)$ of $M$ at $p$ is defined
as the determinant of this operator
$K(p) = \det d\gamma(p)$.  
$K(p)$ does not depend on the choice of orientation.  
The Gauss-Bonnet theorem asserts that 
the average value of $K(\cdot)$
is the alternating sum of the Betti numbers of $M$
$$
\frac{1}{2\pi}\int_M K(p) = \beta_0-\beta_1+\beta_2.  
$$ 

In this section we define a curvature {\it operator} 
associated with any finite rank Hilbert $A$-module 
$H$.  We discuss how it can be viewed as a quantized 
(higher-dimensional) analogue of the differential of the 
Gauss map $\gamma: M\to S^2$, we show that the curvature
operator is of trace class, that its trace agrees 
with the curvature invariant $K(H)$ of section 4, 
and we establish a key asymptotic formula for $K(H)$.  

These considerations are best 
formulated in the following general 
setting.  Let $C$ be a finite dimensional Hilbert space
and let $F=H^2\otimes C$ be the free Hilbert 
$A$-module of rank $r=\dim C$.  

We may consider both $H^2\otimes C$ and the ``Hardy" space 
$H^2(\partial B_d; C)$ of $C$-valued functions as 
spaces of vector-valued holomorphic functions defined 
in the open unit ball $B_d$.  For $H^2\otimes C$ this 
is described in Proposition 2.1.  $H^2(\partial B_d; C)$
is defined as the subspace of 
$L^2(\partial B_d,d\sigma; C)$ obtained by closing 
the space of $C$-valued holomorphic polynomials 
$f: \Bbb C^d\to C$ in the 
$L^2(\partial B_d,d\sigma; C)$-norm; and there is
a natural way of extending functions in 
$H^2(\partial B_d;C)$ holomorphically to the open ball
$B_d$ \cite{26}.  Theorem 4.3 of \cite{1}
implies that these two spaces of holomorphic functions
are related as follows
$$
H^2\otimes C \subseteq H^2(\partial B_d; C).  
\tag{5.2}
$$
Moreover, 
the inclusion map of (5.2) is a compact 
operator when $d\geq 2$.  We will write $F$ instead of 
$H^2\otimes C$, $\partial F$ instead of 
$H^2(\partial B_d; C)$,  
and $b: F\to \partial F$ 
for the inclusion map of (5.2). The 
nature of $b$ and $b^*b$ will be described
more precisely in Proposition 5.7 below.  

We define a linear map 
$\Gamma: \Cal B(F)\to\Cal B(\partial F)$ as follows
$$
\Gamma(X) = bXb^*.  \tag{5.3}
$$

\remark{Remark 5.4}
We first record some simple observations about the 
operator mapping $\Gamma$.  It is obvious that $\Gamma$
is a normal completely positive linear map.  $\Gamma$ is 
also an order isomorphism because $b$ is injective.  Indeed,
if $\Gamma(X)\geq 0$ then $\<X\xi,\xi\>\geq 0$ for every 
$\xi$ in the range $b^*(\partial F)$, and $b^*(\partial F)$
is dense in $F$ because $b$ has trivial kernel.  
A similar argument shows that $\Gamma$ is in fact a 
complete order isomorphism.  However, 
in dimension $d\geq 2$ the range of 
$\Gamma$ is a linear space of compact operators 
which is norm-dense in $\Cal K(\partial F)$ but 
proper: $\Gamma(\Cal B(F))\neq \Cal K(\partial F)$.  
Indeed, if the range of $\Gamma$ were norm-closed
then the closed graph theorem would imply that 
$\Gamma$ is a linear isomorphism of the Banach 
space $\Cal B(F)$ onto $\Cal K(\partial F)$, 
which is obviously absurd since $\Cal B(F)$ is 
inseparable.  
\endremark

$\partial F = H^2(\partial B_d; C)$ is a Hilbert $A$-module 
whose canonical operators $Z_1,\dots,Z_d$ are defined 
by $Z_kf (z) = \<z,e_k\>f(z)$ for $z\in\partial B_d$,  
$e_1,\dots,e_d$ being an orthonormal basis for $\Bbb C^d$.  
$(Z_1,\dots,Z_d)$ is a pure subnormal $d$-contraction
for which 
$$
Z_1^*Z_1+\dots+Z_d^*Z_d = \bold 1,
$$
while
$$
Z_1Z_1^*+\dots+Z_dZ_d^* = \bold 1 - \tilde E_0
$$
$\tilde E_0$ denoting the projection of 
$H^2(\partial B_d; C)$ onto the finite dimensional 
space of constant $C$-valued functions, and 
$\phi_{\partial F}(A) = Z_1AZ_1^*+\dots+Z_dAZ_d^*$ 
is a normal completely positive map on 
$\Cal B(H^2(\partial B_d; C)$.  

\proclaim{Definition 5.4}
$d\Gamma: \Cal B(F)\to \Cal B(\partial F)$ is 
defined as the following linear map
$$
d\Gamma(X) = \Gamma(X) - \phi_{\partial F}(\Gamma(X)) =
\Gamma(X)-\sum_{k=1}^d Z_k\Gamma(X)Z_k^*.  
$$
\endproclaim

We now define the curvature operator of a 
finite rank Hilbert $A$-module $H$.  Let 
$U_0: H^2\otimes \Delta H\to H$ be the homomorphism defined
in section 1,
$$
U_0(f\otimes \zeta) = f\cdot \Delta\zeta, 
\qquad f\in A,\quad\zeta\in \Delta H.
$$

\proclaim{Definition 5.5}
Let $H$ be a finite rank Hilbert $A$-module and
take $F=H^2\otimes \Delta H$, 
$\partial F = H^2(\partial B_d;\Delta H)$ above.  
The curvature 
operator of $H$ is defined as the self-adjoint
operator
$$
d\Gamma(U_0^*U_0)\in \Cal B(H^2(\partial B_d; \Delta H)).  
$$
\endproclaim

\remark{Remarks}
We have found it useful to think of the operator 
$\Gamma(U_0^*U_0)$ as a higher-dimensional 
``quantized" analogue 
of the Gauss map $\gamma: M\to S^2$ of (5.1), and 
of the curvature operator $d\Gamma(U_0^*U_0)$ as 
its ``differential".  Of course, this is only an 
analogy.  But we will also find that 
$d\Gamma(U_0^*U_0)$ belongs to the trace class, and 
$$
\tr d\Gamma(U_0^*U_0) = K(H),
$$
the term $K(H)$ on the right being analogous to 
the average Gaussian curvature 
$$
\frac{1}{2\pi}\int_M K = \frac{1}{2\pi}\int_M \det \gamma(p).  
$$
On the other hand, $K(H)$ is defined in section 4 as the 
integral of the trace (not the determinant) of an 
operator-valued function, and thus this analogy must 
not be carried to extremes.  

We also remark that the curvature operator can be defined
in somewhat more concrete terms as follows.  Let $T(z)$
denote the operator function of $z\in\Bbb C^d$ defined
in (0.2).  $T(z)$ is invertible for $|z|<1$, and hence 
every vector $\xi\in H$ gives rise to a function 
$\hat\xi: B_d \to \Delta H$ by way of 
$$
\hat\xi(z) = \Delta (\bold 1-T(z)^*)^{-1}\xi,
\qquad z\in B_d.  
$$
It is a fact that $\hat\xi$ belongs to 
$\partial F = H^2(\partial F; \Delta H)$, and thus 
we have defined a linear mapping 
$B: \xi\in H\mapsto \hat\xi\in\partial F$.  Indeed, 
the reader can verify that $B$ is related to 
$b$ and $U_0$ by $B=bU_0^*$, and hence the curvature
operator of Definition 5.5 is identical with
$$
BB^*-\phi_{\partial F}(BB^*) = 
BB^* - \sum_{k=1}^d Z_kBB^*Z_k^*.  
$$
We will not have to make use of the operator $B$
in the sequel.  
\endremark

Returning now to the general setting $F=H^2\otimes C$, 
$\partial F=H^2(\partial B_d; C)$, where $C$ is a 
finite dimensional Hilbert space, we work out 
the basic properties of the operator mapping 
$$
d\Gamma: \Cal B(F)\to\Cal B(\partial F).
$$
The essential properties of the inclusion map 
$b: F\to \partial F$ are summarized as follows.  We will 
write $E_n$, $n=0,1,2,\dots$ for the projection of 
$F=H^2\otimes C$ onto its subspace of homogeneous 
(vector-valued) polynomials of degree $n$, and we 
have 
$$
\align
\tr E_n &= \dim\{f\in H^2: f(\lambda z) = \lambda^nf(z), 
 \lambda\in \Bbb C, z\in\Bbb C^d\}\cdot\dim C \\
&= q_{d-1}(n)\cdot\dim C,
\endalign
$$ 
where $q_{d-1}(x)$ is the $d-1$st term in the sequence 
of polynomials $q_0(x),q_1(x),\dots$ of (3.3).  Since this 
polynomial will play an important part in the remainder 
of this section, we reiterate its definition here: 
for $d=1$, $q_{d-1}(x)=1$ and otherwise
$$
q_{d-1}(x) = \frac{(x+1)\dots(x+d-1)}{(d-1)!}, \qquad d\geq 2.  
\tag{5.6}
$$
Let $\tilde E_d$ be the corresponding sequence of 
projections acting on the Hilbert $A$-module $\partial F$.  
We will also write $N$ and $\tilde N$ for the respective 
number operators on $F$ and $\partial F$, 
$$
N=\sum_{n=0}^\infty n E_n,\quad \tilde N = 
\sum_{n=0}^\infty n\tilde E_n.  
$$

\proclaim{Proposition 5.7}
Let $b: F\to \partial F$ be the natural inclusion.  Then
\roster
\item
$bE_n = \tilde E_n b, \qquad n=0,1,2,\dots.  $
\item
$b\in hom(F,\partial F).  $
\item
$b^*b = q_{d-1}(N)^{-1} = 
\sum_{n=0}^\infty \frac{1}{q_{d-1}(n)}E_n.$  
\endroster
\endproclaim
\demo{proof}
Properties (1) and (2) are immediate from the definition 
of $b$.  Property (3) follows from a direct comparison 
of the norms in $H^2$ and $H^2(\partial B_d)$.  Indeed, if 
$f,g\in H^2$ are both homogeneous polynomials of degree
$n$ of the specific form
$$
f(z) = \<z,\alpha\>^n, \quad g(z) = \<z,\beta\>^n, \qquad 
\alpha,\beta\in \Bbb C^d, 
$$
then $\<f,g\>_{H^2} = \<\beta,\alpha\>^n$, whereas 
if we consider $f,g$ as elements of $H^2(\partial B_d)$ then 
we have 
$$
\<bf,bg\> = \<f,g\>_{H^2(B_d)} = 
q_{d-1}(n)^{-1}\<\beta,\alpha\>^n, 
$$
see Proposition 1.4.9 of \cite{26}.  Since 
$E_nH^2$ is spanned by such $f,g$ we find that for 
all $f,g\in E_nH^2$, 
$$
\<bf,bg\> = q_{d-1}(n)^{-1}\<f,g\>_{H^2}.   
$$
Thus 
$$
E_nb^*bE_n = q_{d-1}(n)^{-1}E_n = q_{d-1}(N)^{-1}E_n,
$$
and (3) follows for the case $C=\Bbb C$ from the fact 
that $b^*b$ commutes with $E_n$ and $\sum_nE_n = \bold 1$.  

If we now tensor both $H^2$ and $H^2(\partial B_d)$ with 
the finite dimensional space $C$ then we obtain 
(3) in general after noting that 
$\dim(K_1\otimes K_2) = \dim K_1\cdot\dim K_2$ for finite
dimensional vector spaces $K_1, K_2$.  
\qedd
\enddemo

\remark{Remark 5.8}
In the one-variable case $d=1$, $q_{d-1}(x)$ is the constant
polynomial $1$, and hence 5.7 (3) asserts the obvious fact 
that $b$ is a unitary operator; i.e., there is no difference 
between $H^2$ and $H^2(S^1)$ in dimension one.  

In dimension $d\geq 2$ however, $q_{d-1}(x)$ is a polynomial 
of degree $d-1\geq 1$ and hence 
$$
b^*b = q_{d-1}(N)^{-1}
$$
is a positive compact operator.  Significantly, the 
operator $b^*b$ is never trace class.  
Indeed, the computations of \cite{1, Appendix A} imply 
that $b^*b \in \Cal L^p$ 
iff $p > \frac{d}{d-1} > 1$.  
\endremark

Returning now to the discussion of $d\Gamma$, we begin 
by giving a description of the cone of all operators 
$X$ whose ``differential" $d\Gamma(X)$ is positive.  

\proclaim{Theorem 5.9}
Let $F=H^2\otimes C$ be a free Hilbert module.  Then 
for any $X\in \Cal B(F)$ the following are equivalent.  
\roster
\item
$d\Gamma(X)\geq 0$.  
\item
There is a sequence $\Phi_1,\Phi_2,\dots$ of 
bounded endomorphisms 
of $F$ such that 
$$
X = \sum_{n=1}^\infty \Phi_n\Phi_n^*.  
$$
\item
There is a free Hilbert module $\tilde F$ and a 
bounded homomorphism of Hilbert $A$-modules 
$\Phi:\tilde F\to F$ such that 
$X=\Phi\Phi^*$.  
\endroster
\endproclaim

\demo{proof}
We prove the sequence of implications 
(3)$\implies$ (2)$\implies$ (1)$\implies$ (3).  

\demo{proof of (3)$\implies$(2)}
Let $\tilde F$ be a free Hilbert module and let 
$\Phi: \hom(\tilde F,F)$ satisfy $\Phi\Phi^*=X$. 
We may assume that $\tilde F$ has infinite rank 
by adding a direct summand $H^2\otimes \ell^2$ 
to it, and by extending $\Phi$ to the larger free 
Hilbert module by making it zero on the summand 
$H^2\otimes \ell^2$.   

Assuming that this has been arranged, and taking note
of the isomorphism $\tilde F \cong F\oplus F\oplus \dots$,
we may assume that $\tilde F = F\oplus F\oplus\dots$ is 
a direct sum of copies of $F$.  Defining endomorphisms 
$\Phi_n\in\Cal B(F)$ by simply restricting $\Phi$ to the 
$n$th summand, we deduce the required representation
$$
X = \Phi\Phi^* = \sum_{n=1}^\infty \Phi_n\Phi_n^*.  
$$
\enddemo

\demo{proof of (2)$\implies$(1)}
Since $d\Gamma$ is a normal operator mapping, it suffices
to show that $d\Gamma(\Phi\Phi^*)\geq 0$ for every 
endomorphism $\Phi: F\to F$.  For that, write 
$$
\Gamma(\Phi\Phi^*) = b\Phi\Phi^*b^* = (b\Phi)(b\Phi)^*,
$$
where $b\Phi\in\hom(F,\partial F)$ is the composite 
homomorphism of Hilbert $A$-modules.  Thus if
$T_1,\dots,T_d$ are the canonical operators of $F$ 
and $Z_1,\dots,Z_d$ are those of $\partial F$, then 
we have $Z_kb\Phi = b\Phi T_k$ for every $k$, hence
$$
\sum_{k-1}^dZ_k\Gamma(\Phi\Phi^*)Z_k^* =
b\Phi(\sum_{k=1}^d T_kT_k^*)(b\Phi)^*,
$$
and therefore 
$$
d\Gamma(\Phi\Phi^*) = 
b\Phi(\bold 1-\sum_{k=1}^d T_kT_k^*)(b\Phi)^*
$$
is a positive operator because $(T_1,\dots,T_d)$ 
is a $d$-contraction.  
\enddemo

\demo{proof of (1)$\implies$(3)}
Let $T_1,\dots,T_d$ be the canonical operators of $F$, 
and suppose $d\Gamma(X)\geq 0$.  We claim first that 
$$
0\leq \sum_{k=1}^d T_kXT_k^*\leq X.  
\tag{5.10}
$$
To see that, note that for every $k=1,\dots,d$ we have 
$bT_k=Z_kb$, $Z_1,\dots,Z_d$ being the canonical operators
of $\partial F$.  Hence
$$
0\leq d\Gamma(X) = bXb^*-\sum_{k=1}^dZ_kbXb^*Z_k^* =
bXb^*-b(\sum_{k=1}^dT_kXT_k^*)b^* = 
\Gamma(X-\sum_{k=1}^dT_kXT_k^*).  
$$
Since $\Gamma$ is an order isomorphism the latter implies 
$X-\sum_kT_kXT_k^*\geq 0$, or 
$$
X-\phi(X)\geq 0, 
\tag{5.11}
$$
$\phi$ being the completely positive map of $\Cal B(F)$ 
defined by 
$$
\phi(A) = \sum_{k=1}^dT_kAT_k^*.
$$
Free Hilbert $A$-modules are pure, hence 
$\phi^n(\bold 1)\downarrow 0$ as $n\to\infty$.  It follows 
that for every positive operator $A\in \Cal B(F)$ 
we have $0\leq \phi^n(A)\leq \|A\|\phi^n(\bold 1)$, 
and hence $\phi^n(A)\to 0$ in the strong operator 
topology of $\Cal B(F)$, as $n\to\infty$.  By taking linear 
combinations we find that $\lim_{n\to\infty}\phi^n(A) = 0$
in the strong operator topology 
for every $A\in\Cal B(F)$.  

Returning now to equation (5.11), we find that 
$$
X-\phi^{n+1}(X) = \sum_{k=1}^n\phi^k(X-\phi(X)) \geq 0
$$
for every $n=1,2,\dots$ and since $\phi^{n+1}(X)$ must
tend strongly to $0$ by the preceding paragraph, 
we conclude that  $X\geq 0$ 
by taking the limit on $n$ in the preceding inequality.  
Thus we may add the positive operator $\phi(X)$ to 
(5.11) to obtain the desired inequality (5.10).  

Now consider the closed subspace $K\subseteq F$ obtained
by closing the range of the positive operator $X^{1/2}$.  
We will make $K$ into a {\it pure} Hilbert $A$-module 
as follows.  

We claim first that there is a unique 
$d$-contraction $\tilde T_1,\dots,\tilde T_d$ acting 
on $K$ such that 
$$
T_k X^{1/2} = X^{1/2} \tilde T_k, \qquad k=1,2,\dots,d.  
$$
Indeed, the uniqueness of $\tilde T_1,\dots,\tilde T_d$ is 
clear from the fact that if $\xi\in K$ and $X^{1/2}\xi=0$, 
then $\xi=0$; that is simply because 
$K$ is the closure of the range of $X^{1/2}$, 
hence the kernel of the restriction of $X^{1/2}$ to $K$ 
is trivial.  

In order to construct the operators $\tilde T_k$ it is easier 
to work with adjoints, and we will 
define operators $A_k = \tilde T_k^*$ as 
follows.  Fix $k=1,\dots,d$ and $\xi\in F$.  Then 
by (5.10) we have 
$$
\|X^{1/2}T_k^*\xi\|^2 \leq 
\sum_{k=1}^d\|X^{1/2}T_k^*\xi\|^2 = 
\sum_{k=1}^d\<T_kXT_k^*\xi,\xi\> \leq 
\<X\xi,\xi\> = \|X^{1/2}\xi\|^2,  
$$
and hence there is a unique contraction $A_k\in \Cal B(K)$ 
such that 
$$
A_kX^{1/2} = X^{1/2}T_k^*, \qquad k=1,\dots,d. 
\tag{5.12} 
$$
As in the previous estimate (5.12) implies 
$$
\sum_{k=1}^d\|A_kX^{1/2}\xi\|^2 \leq \|X^{1/2}\xi\|^2,
\qquad \xi\in F,
$$
and hence $A_1^*A_1+\dots+A_d^*A_d\leq \bold 1_K$.  
Since the $T_k^*$ mutually commute, (5.12) implies 
that the $A_k$ must mutually commute, and we 
set $\tilde T_k=A_k^*$, $k=1,\dots,d$.

Next, we claim that $(\tilde T_1,\dots,\tilde T_d)$ is 
a pure $d$-contraction in the sense that if 
$\tilde \phi: \Cal B(K)\to \Cal B(K)$ is the map 
defined by 
$$
\tilde \phi(A) = \sum_{k=1}^d\tilde T_kA\tilde T_k^*, 
$$
then $\tilde \phi^n(\bold 1_K)\downarrow 0$ as 
$n\to \infty$.  Since $\{\tilde\phi^n(\bold 1_K):n\geq 0\}$ 
is a uniformly bounded sequence of positive operators, 
the claim will follow if we show that 
$$
\lim_{n\to\infty}\<\tilde\phi^n(\bold 1_K)\eta,\eta\>=0
$$
for all $\eta$ in the dense linear manifold $X^{1/2}F$ 
of $K$.  But for $\eta = X^{1/2}\xi$, $\xi \in F$, we 
have 
$$
\<\tilde\phi^n(\bold 1_K)X^{1/2}\xi,X^{1/2}\xi\> 
= \<X^{1/2}\tilde\phi^n(\bold 1_K)X^{1/2}\xi,\xi\>.  
$$
Since $X^{1/2}\tilde T_k = T_kX^{1/2}$ for all $k$ 
it follows that $X^{1/2}\tilde\phi^n(\bold 1_K)X^{1/2} =
\phi^n(X)$ for every $n=0,1,2,\dots$, hence
$$
\<\tilde\phi^n(\bold 1_K)X^{1/2}\xi,X^{1/2}\xi\> = 
\<\phi^n(X)\xi,\xi\> \leq \|X\|\<\phi^n(\bold 1)\xi,\xi\>
$$
and the right side decreases to zero as $n\to\infty$ 
because the free module $F$ is pure.  

Using the operators $\tilde T_1,\dots,\tilde T_d\in \Cal B(K)$
we now consider $K$ to be a pure Hilbert $A$-module.  By the 
basic dilation theory (see Theorem 1.14) there is a free Hilbert 
$A$-module $\tilde F$ and a {\it coisometry} 
$U\in \hom(\tilde F,K)$.  Thus if we denote the canonical 
operators of $\tilde F$ by $S_1,\dots,S_d$, then we have 
$US_k=\tilde T_kU$ for $k=1,\dots,d$.  Now define a linear 
operator $\Phi\in\Cal B(\tilde F,F)$ by $\Phi = X^{1/2}U$.  
We have 
$$
\Phi S_k=X^{1/2}US_k=X^{1/2}\tilde T_k U=
T_kX^{1/2}U=T_k \Phi,
$$
so that $\Phi\in\hom(\tilde F,F)$.  Finally, since 
$UU^*=\bold 1_K$, it follows that 
$$
\Phi\Phi^* = X^{1/2}UU^*X^{1/2} = X,
$$
and the proof of Theorem 5.9 is complete.  
\qedd
\enddemo
\enddemo

\remark{Remark}
Notice that the proof of Theorem 5.9 made no 
use of the finite dimensionality of $C$, and 
in fact this result is valid {\it verbatim} 
for free Hilbert $A$-modules of arbitrary 
rank $1,2,\dots,\infty$.  
\endremark

We digress momentarily to record the following 
(a fact associated with the dilation theory of 
\cite{1, Theorem 8.5}) for later use.  

\proclaim{Corollary}
Let $M$ be a closed submodule of $H^2\otimes C$ 
and let $P_M$ be the projection onto $M$.  Then 
there is a Hilbert space $E$ and a partial 
isometry $\Phi\in\hom(H^2\otimes E,H^2\otimes C)$
such that $P_M=\Phi\Phi^*$.  
\endproclaim
\demo{proof}
Let $S_1,\dots,S_d$ be the canonical operators 
of $H^2\otimes C$.  Since
$$
S_1P_MS_1^*+\dots+S_dP_MS_d^*\leq P_M
$$
and since $bS_k=Z_kb$, $k=1,\dots,d$ we have 
$$
d\,\Gamma(P_M)=b(P_M-\sum_{k=1}^dS_kP_MS_k^*)b^*\geq 0
$$
and the conclusion follows from condition 
(3) of Theorem 5.9.
\qedd
\enddemo

The importance of the cone 
$$
\Cal P = \{X\in \Cal B(F): d\Gamma(X)\geq 0\}
$$
for this work is that for every operator $X$ in 
the complex vector space spanned by $\Cal P$, 
$d\Gamma(X)$ is trace class, and $\tr\,d\Gamma(X)$
can be expressed in terms of $X$ by way of an 
asymptotic formula.  

\proclaim{Theorem 5.13}
Let $F=H^2\otimes C$, where $C$ is a finite dimensional Hilbert
space.  For every operator $X$ in the complex linear span 
of $\{X\in\Cal B(F): d\Gamma(X)\geq 0\}$, 
$d\Gamma(X)$ belongs to $\Cal L^1(\partial F)$ and 
$$
\tr (d\Gamma(X)) = \dim C\cdot
\lim_{n\to\infty}\frac{\tr (XE_n)}{\tr E_n},
$$
$E_0,E_1,\dots$ being the sequence of spectral projections 
of the number operator of $F$.  
\endproclaim

Theorem 5.13 depends on a general identity, 
which we establish first.  

\proclaim{Lemma 5.14}
Let $F=H^2\otimes C$ be as in Theorem 5.13.  Then for 
every $X\in\Cal B(F)$ and $n=0,1,2,\dots$ we have 
$$
\tr (d\Gamma(X)\tilde P_n)=\dim C\cdot
\frac{\tr (XE_n)}{\tr E_n},
$$
where $\tilde P_n = \tilde E_0+\tilde E_1+\dots+\tilde E_n$,
$\{\tilde E_n\}$ being the spectral projections of the number 
operator of $\partial F$.  
\endproclaim

\remark{Remark}
Notice that all of the operators 
$E_n, XE_n, d\Gamma(X)\tilde P_n$ appearing in Lemma 5.14 are
of finite rank.  Note too that traces on the right refer to 
the Hilbert space $F$, while the trace on the left refers 
to the Hilbert space $\partial F$.  
\endremark

\demo{proof of Lemma 5.14}
Let $\dot E_n$ be the projection of $H^2$ onto 
its space of homogeneous polynomials of degree $n$.  Then
$E_n = \dot E_n\otimes\bold 1_C$, and hence 
$$
\tr E_n = \tr \dot E_n\cdot \dim C = 
q_{d-1}(n) \cdot \dim C,
$$
for all $n=0,1,\dots$ where $q_{d-1}(x)$ is the polynomial 
of (5.6).  Thus we have to show that 
$$
\tr (d\Gamma(X)\tilde P_n) = \frac{\tr(XE_n)}{q_{d-1}(n)},
\qquad n=0,1,\dots.  \tag{5.15}
$$ 

Using Proposition  5.7 we have
$q_{d-1}(n)^{-1}E_n = b^*bE_n = b^*\tilde E_nb$
and the right side of (5.15) can be rewritten
$$
\frac{\tr(XE_n)}{q_{d-1}(n)} =
\tr(Xb^*\tilde E_nb) = \tr(\Gamma(X)\tilde E_n).
$$
Setting $Y=\Gamma(X)\in\Cal B(\partial F)$, equation
(5.15) becomes
$$
\tr((Y-\phi_{\partial F}(Y))\tilde P_n)=
\tr(Y\tilde E_n),
\tag{5.16}
$$
where $\phi_{\partial F}$ is the completely positive
map on $\Cal B(\partial F)$ defined by the canonical 
operators of $\partial F$, 
$\phi_{\partial F}(A) = Z_1AZ_1^*+\dots+Z_dAZ_d^*$.

Notice that since 
$Z_1^*Z_1+\dots+Z_d^*Z_d=\bold 1_{\partial F}$, 
$\phi_{\partial F}$ leaves the trace invariant in the 
sense that for $A\in\Cal L^1(\partial F)$ we have
$$
\tr(\phi_{\partial F}(A))=
\tr(\sum_{k=1}^dZ_k^*Z_kA)=\tr(A).  
$$
Moreover, the relations
$Z_k^*\tilde P_n=\tilde P_{n-1}Z_k^*$ imply that 
$\phi_{\partial F}(A)\tilde P_n=
\phi_{\partial F}(A\tilde P_{n-1})$ for $n=0,1,\dots$,
where of course $\tilde P_{-1}$ is taken as $0$.  Thus
$$
\tr(\phi_{\partial F}(Y)\tilde P_n) =
\tr(\phi_{\partial F}(Y\tilde P_{n-1})) =
\tr(Y\tilde P_{n-1})
$$
and the left side of (5.16) can be written 
$$
\tr(Y\tilde P_n)-\tr(\phi_{\partial F}(Y)\tilde P_n) =
\tr(Y\tilde P_n)-\tr(Y\tilde P_{n-1})
$$
which agrees with the right side of (5.16) because
$\tilde P_n-\tilde P_{n-1}=\tilde E_n$.
\qedd
\enddemo

\demo{proof of Theorem 5.13}
It suffices to show that for any operator $X$ 
in $\Cal B(F)$ for which $d\Gamma(X)\geq 0$, we 
must have $\tr(d\Gamma(X))<\infty$ as well as 
the limit formula of 5.13.  From Lemma 5.14
we have 
$$
\tr(d\Gamma(X)\tilde P_n) = \dim C\cdot
\frac{\tr(XE_n)}{\tr E_n}, \qquad n=0,1,2,\dots.  
\tag{5.17}
$$
Theorem 5.9 implies that $X$ must be a positive 
operator and hence 
$$
0\leq \frac{\tr (XE_n)}{\tr E_n} \leq \|X\|
$$
for every $n$ because $\rho_n(A)=\tr (AE_n)/\tr E_n$ is
a state of $\Cal B(F)$.  Since the projections 
$\tilde P_n$ increase to $\bold 1_{\partial F}$ with 
increasing $n$ we conclude from (5.17) that 
$$
\tr (d\Gamma(X)) = 
\sup_{n\geq 0}\tr (d\Gamma(X)\tilde P_n) \leq
\dim C\cdot \|X\| <\infty.  
$$
Moreover, since in this case
$$
\tr (d\Gamma(X)) = 
\lim_{n\to\infty}\tr (d\Gamma(X)\tilde P_n), 
$$
we may infer the required 
limit formula directly 
from (5.17) as well.
\qedd
\enddemo

In view of Theorem 5.9, 
the following lemma shows how to compute the trace
of $d\Gamma(X)$ in the most important cases.  

\proclaim{Lemma 5.18}
Let $F_k = H^2\otimes C_k$, $k=1,2$, where $C_2$ is 
finite dimensional, and let $\Phi\in\hom(F_1,F_2)$.  
Considering $\Phi$ as a multiplier in 
$\Cal M(C_1,C_2)$ with boundary value function 
$\tilde\Phi: \partial B_d\to \Cal B(C_1,C_2)$ we 
have 
$$
\tr d\Gamma(\Phi\Phi^*) = 
\int_{\partial B_d}
\tr(\tilde\Phi(z)\tilde\Phi(z)^*)\,d\sigma(z).  
$$
\endproclaim 

\remark{Remark}
Note that for $\sigma$-almost every $z\in\partial B_d$,
$\tilde\Phi(z)\tilde \Phi(z)^*$ is a positive operator 
in $\Cal B(C_2)$, and since $C_2$ is finite dimensional
the right side is well defined and dominated by 
$\|\Phi\|^2\cdot\dim C_2$.  
\endremark

\demo{proof of Lemma 5.18}
Consider the linear operator 
$A: C_1\to H^2(\partial B_d;C_2)$ defined by
$A\zeta = b(\Phi(1\otimes \zeta))$, $\zeta\in C_1$.  
We claim first that 
$$
d\Gamma(\Phi\Phi^*) = AA^*.  
\tag{5.19}
$$
Indeed, since $b\Phi\in\hom(F_1,H^2(\partial B_d,C_2))$
we have 
$$
\sum_{k=1}^d Z_kb\Phi\Phi^*b^*Z_k^* = 
\sum_{k=1}^d Z_k(b\Phi)(b\Phi)^*Z_k^* = 
b\Phi(\sum_{k=1}^dT_kT_k^*)(b\Phi)^*,
$$
where $T_1,\dots,T_d$ are the canonical operators of 
$F_1=H^2\otimes C_1$, and hence 
$$
d\Gamma(\Phi\Phi^*) = 
b\Phi(b\Phi)^*-\sum_{k=1}^dZ_kb\Phi\Phi^*b^*Z_k^*=
b\Phi(\bold 1_{F_1}-\sum_{k=1}^dT_kT_k^*)(b\Phi)^*.  
$$
The operator $\bold 1_{F_1}-\sum_k T_kT_k^*$ is the 
projection of $F_1=H^2\otimes C_1$ onto its space of 
$C_1$-valued constant functions and, denoting by 
$[1]$ the projection of $H^2$ onto the one dimensional
space of constants $\Bbb C\cdot 1$, the preceding 
formula becomes 
$$
d\Gamma(\Phi\Phi^*) = b\Phi([1]\otimes\bold 1_{C_1})(b\Phi)^*
= AA^*, 
$$
as asserted in (5.19).  

Now fix an orthonormal basis $e_1, e_2,\dots$ for $C_1$.  
By formula (5.19) we can evaluate the trace of 
$d\Gamma(\Phi\Phi^*)$ in terms of the vector functions 
$Ae_n\in H^2(\partial B_d,C_2)$ as follows,
$$
\align
\tr\, d\Gamma(\Phi\Phi^*) &= 
\tr_{H^2(\partial B_d;C_2)}(AA^*) =
\tr_{C_1}(A^*A) \tag{5.20}\\
&= \sum_n\<A^*Ae_n,e_n\>=\sum_n\|Ae_n\|_{H^2(\partial B_d;C_2)}^2. 
\endalign 
$$

Turning now to the term on the right in Lemma 5.18, 
we first consider $Ae_n=b(\Phi(1\otimes e_n))$ as 
a function from the open ball $B_d$ to $C_2$.  In 
terms of the multiplier $\Phi(\cdot)$ of $\Phi$ we 
have 
$$
Ae_n(z) = b\Phi(1\otimes e_n)(z) = \Phi(z)e_n
$$
and hence the boundary values $\tilde{Ae_n}$
of $Ae_n$ are given by
$\tilde{Ae_n}(z) = \tilde\Phi(z)e_n$ 
for $\sigma$-almost every $z\in\partial B_d$.  Thus 
for such $z\in\partial B_d$ we have 
$$
\tr_{C_2}(\tilde\Phi(z)\tilde\Phi(z)^*) =
\tr_{C_1}(\tilde\Phi(z)^*\tilde\Phi(z)) =
\sum_n\|\tilde\Phi(z)e_n\|^2 = \sum_n\|\tilde{Ae_n}(z)\|^2.  
$$
Integrating the latter over the sphere we obtain
$$
\int_{\partial B_d}\tr_{C_2} (\tilde\Phi(z)\tilde\Phi(z)^*)\,d\sigma
= \sum_n\|\tilde{Ae_n}\|^2_{H^2(\partial B_d;C_2)}
$$
and from (5.20) we see that this coincides with 
$\tr\,d\Gamma(\Phi\Phi^*)$.
\qedd
\enddemo

We now establish the main asymptotic formula for $K(H)$.  

\proclaim{Theorem D}
For every finite rank Hilbert $A$-module $H$, the 
curvature operator $d\Gamma(U_0^*U_0)$ belongs to the 
trace class $\Cal L^1(H^2(\partial B_d;\Delta H))$, 
and we have 
$$
K(H) = \tr\,d\Gamma(U_0^*U_0) = 
d!\lim_{n\to\infty}\frac{\tr (\bold 1-\phi^{n+1}(\bold 1))}
{n^d}
$$
where $\phi:\Cal B(H)\to\Cal B(H)$ is the canonical 
completely positive map associated with the 
$A$-module structure of $H$.  
\endproclaim

Let $\Delta=(\bold 1-\phi(\bold 1))^{1/2}$.  We will actually
prove a slightly stronger assertion, namely
$$
K(H) = \tr\,d\Gamma(U_0^*U_0) = 
(d-1)!\lim_{n\to\infty}\frac{\tr (\phi^n(\Delta^2))}
{n^{d-1}}.
\tag{5.21}
$$
We first point out that it suffices to prove (5.21).  
For that, let $a_k=\tr \phi^k(\Delta^2)$, $k=0,1,2,\dots$. 
Since 
$$
\bold 1-\phi^{n+1}(\bold 1) = 
\sum_{k=0}^n \phi^k(\bold 1-\phi(\bold 1))=
\sum_{k=0}^n \phi^k(\Delta^2)
$$
and since for every $r=1,2,\dots$
$$
q_r(n)=\frac{(n+1)\dots(n+r)}{r!}\sim \frac{n^r}{r!}, 
$$
we have
$$
d!\frac{\tr (\bold 1-\phi^{n+1}(\bold 1))}{n^d} \sim
\frac{\tr (\bold 1-\phi^{n+1}(\bold 1))}{q_d(n)} =
\frac{a_0+a_1+\dots+a_n}{q_d(n)}
$$
while 
$$
(d-1)!\frac{\tr \phi^n(\Delta^2)}{n^{d-1}} \sim
\frac{\tr \phi^n(\Delta^2)}{q_{d-1}(n)} =
\frac{a_n}{q_{d-1}(n)}
$$
Thus the following elementary lemma allows one to deduce
Theorem D from (5.21).  

\proclaim{Lemma 5.22}
Let $d=1,2,\dots$ and let $a_0,a_1,\dots$ be a sequence of 
real numbers such that 
$$
\lim_{n\to\infty}\frac{a_n}{q_{d-1}(n)}=L\in\Bbb R.  
$$
Then
$$
\lim_{n\to\infty}\frac{a_0+a_1+\dots+a_n}{q_d(n)}=L.  
$$
\endproclaim
\demo{proof of Lemma 5.22}
Choose $\epsilon>0$.  By hypothesis, there is an 
$n_0\in\Bbb N$ such that 
$$
(L-\epsilon)q_{d-1}(k)\leq a_k\leq (L+\epsilon)q_{d-1}(k),
\qquad k\geq n_0.  
\tag{5.23}
$$
By the recursion formula (3.2.2) we have 
$$
\sum_{k=n_0}^nq_{d-1}(k)=
\sum_{k=n_0}^n(q_d(k)-q_d(k-1)) = q_d(n)-q_d(n_0-1).  
$$
Thus if we sum (5.26)  from $n_0$ to $n$ and 
divide through by 
$q_d(n)$ we obtain
$$
(L-\epsilon)(1-\frac{q_d(n_0-1)}{q_d(n)}) \leq 
\frac{a_{n_0}+\dots+a_n}{q_d(n)} \leq
(L+\epsilon)(1-\frac{q_d(n_0-1)}{q_d(n)}).  
$$
Since $q_d(n)\to\infty$ as $n\to\infty$, the latter
inequality implies 
$$
L-\epsilon\leq 
\liminf_{n\to\infty}\frac{a_0+\dots+a_n}{q_d(n)} \leq 
\limsup_{n\to\infty}\frac{a_0+\dots+a_n}{q_d(n)} \leq
L+\epsilon,
$$
and since $\epsilon$ is arbitrary, 5.24 follows.
\qedd
\enddemo
\demo{proof of Theorem D}
Let $U_0:H^2\otimes \Delta H\to H$ be the homomorphism 
$U_0(f\otimes \zeta)=f\cdot\Delta\zeta$ defined in 
section 1 and discussed above.  We claim that
for every $n=0,1,\dots$
$$
\tr \phi^n(\Delta^2) = \tr (U_0^*U_0E_n),
\tag{5.24}
$$
$E_n\in\Cal B(H^2\otimes\Delta H)$ 
being the projection onto the space of homogeneous
polynomials of degree $n$.  Indeed, from the discussion 
preceding (1.13) we have
$$
U_0U_0^*=\bold 1-\phi^\infty(\bold 1),
$$
and since 
$\phi(\phi^\infty(\bold 1))=\phi^\infty(\bold 1)$
we can write 
$$
\Delta^2=\bold 1-\phi(\bold 1) = 
(\bold 1-\phi^\infty(\bold 1))-
\phi(\bold 1-\phi^\infty(\bold 1)) =
U_0U_0^*-\phi(U_0U_0^*).  
$$
Thus
$$
\phi^n(\Delta^2)=\phi^n(U_0U_0^*)-\phi^{n+1}(U_0U_0^*).
\tag{5.25}
$$
Write $F=H^2\otimes\Delta H$, and let 
$\phi_F: \Cal B(F)\to\Cal B(F)$ denote its canonical 
completely positive map.  Since $U_0\in\hom(F,H)$ we 
have 
$$
\phi^k(U_0U_0^*)=U_0\phi_F^k(\bold 1_F)U_0^*
$$
for every $k=0,1,\dots$.  Moreover, 
$$
\phi_F^n(\bold 1_F)-\phi_F^{n+1}(\bold 1_F) =
\phi_F^n(\bold 1_F-\phi_F(\bold 1_F)) =
\phi_F^n(E_0)=E_n,
$$
so that (5.25) implies 
$$
\phi^n(\Delta^2) = U_0E_nU_0^*.  
$$
The formula (5.24) follows immediately since 
$$
\tr_H(U_0E_nU_0^*)=\tr_F(U_0^*U_0E_n).  
$$

By Proposition 4.3 there is a free module 
$\tilde F$ and $\Phi\in\hom(\tilde F,F)$ such that 
$$
U_0^*U_0 = \bold 1_F-\Phi\Phi^*.  
$$
Since both $d\Gamma(\bold 1_F)$ and 
$d\Gamma(\Phi\Phi^*)$ are positive operators
by Theorem 5.9 
(indeed $d\Gamma(\bold 1_F)$ is the projection 
of $H^2(\partial B_d; \Delta H)$ onto its subspace 
of constant functions), it follows from Theorem 
5.13 that the curvature operator $d\Gamma(U_0^*U_0)$ 
is trace class and, in view of (5.24), satisfies
$$
\tr\,d\Gamma(U_0^*U_0) = 
\lim_{n\to\infty}\frac{\tr(U_0^*U_0E_n)}{q_{d-1}(n)}=
\lim_{n\to\infty}\frac{\tr \phi^n(\Delta^2)}{q_{d-1}(n)}.
\tag{5.26}
$$

Finally, we use (5.18) together with 
$U_0^*U_0=\bold 1_F-\Phi\Phi^*$ to evaluate the left side
of (5.26) and we find that 
$$
\tr\,d\Gamma(U_0^*U_0) =
\int_{\partial B_d} \tr 
(\bold 1_{\Delta H}-\tilde\Phi(z)\tilde\Phi(z)^*)\,d\sigma(z).
$$
Formula (4.10) shows that the term 
on the right is $K(H)$.
\qedd
\enddemo

As in the case of the Euler characteristic, the 
asymptotic formula of Theorem D leads to the following
result on stability under finite dimensional perturbations.  

\proclaim{Corollary 1: stability of curvature}
Let $H_0$ be a closed submodule of a 
finite rank Hilbert $A$-module $H$ such 
that $\dim(H/H_0)<\infty$.  Then $K(H)=K(H_0)$.  
\endproclaim
\demo{proof}
This is proved by estimating exactly as in the proof of 
the corollary of Theorem C.  Indeed, the estimates here 
are simpler because the trace is a linear functional.  
One finds that 
$$
|\tr(\bold 1_H-\phi_H^{n+1}(\bold 1_H)) -
\tr(\bold 1_{H_0}-\phi_{H_0}^{n+1}(\bold 1_{H_0}))|
\leq\dim(H/H_0)(1+q_{d-1}(n+1)).  
$$
As in the proof of the corollary of Theorem 
C, one can multiply through by $d!/n^d$ and take the 
limit on $n$ to obtain the required relation 
$|K(H)-K(H_0)|\leq 0$.
\qedd
\enddemo

We also point out the following application to invariant
subspaces of the $d$-shift $S_1,\dots,S_d$ acting on 
$H^2$.  In dimension $d=1$ the invariant subspaces of 
the simple unilateral shift define submodules which are 
isomorphic to $H^2$ itself, and 
in particular they all have 
rank one.  In higher dimensions, on the other hand, 
we can never have that behavior for 
the ranks of submodules of finite codimension.  

\proclaim{Corollary 2}
Suppose that $d\geq 2$, and let 
$M$ be a proper closed submodule of $H^2$ of finite 
codimension.  Then $\rank\,M\geq 2$.  
\endproclaim
\demo{proof}
Theorem 4.13 implies that $K(H^2)=1$, and hence corollary 
1 above implies that $K(M)=1$ as well.  
By \cite{1,Lemma 7.14},
no proper submodule of $H^2$ can be a free Hilbert 
module in dimension $d>1$, hence 
by the extremal property of $K(M)$ (Theorem 4.13)
we must have $\rank(M)>K(M)=1$.  
\qedd
\enddemo
\remark{Remark}
Of course, the ranks of finite codimensional submodules 
of $H^2$ must be finite by Remark 3.14, 
and they can be arbitrarily large.
\endremark

Since the trace of a positive operator $A$ is dominated by 
$\|A\|\cdot\rank(A)$, Theorems C and D together imply that 
$K(H)\leq \chi(H)$, and we conclude 

\proclaim{Corollary 3}
For every finite rank Hilbert $A$-module $H$, 
$$
0\leq K(H)\leq \chi(H)\leq \rank(H).  
$$
\endproclaim

In the next section we will show that $K(H)=\chi(H)$ for 
graded Hilbert modules, but in section 9 we give examples
of ungraded Hilbert modules for which $K(H)<\chi(H)$.  
The inequality of Corollary 3 is useful; a 
significant application is given in Theorem E of section 8.

\subheading{6. Graded Hilbert modules}

In this section we prove 
an analogue of the Gauss-Bonnet-Chern theorem for
Hilbert $A$-modules.  The most general setting in 
which one might hope for such a result is the class
of finite rank {\it pure} Hilbert $A$-modules.  By 
the discussion of section 1, these are the Hilbert 
$A$-modules which are isomorphic to quotients 
$F/M$ of finite rank free modules $F$ by closed 
submodules $M$.  However, in Proposition 9.2 we  
give examples of submodules $M\subseteq H^2$ 
for which $K(H^2/M)\neq \chi(H^2/M)$.  In this 
section we establish the result (Theorem B) under 
the additional hypothesis that $H$ is {\it graded}.  
Examples are obtained by taking $H=F/M$ where 
$F$ is free of finite rank and $M$ is a closed 
submodule generated by a set of homogeneous 
polynomials (perhaps of different degrees).  
In particular, one can associate
such a module $H$ with any algebraic 
variety in complex projective space $\Bbb P^{d-1}$
(see section 9).

By a {\it graded} Hilbert space we mean a pair 
$H,\Gamma$ where $H$ is a (separable) Hilbert space 
and $\Gamma: \Bbb T\to\Cal B(H)$ 
is a strongly continuous unitary representation of 
the circle group
$\Bbb T=\{\lambda\in\Bbb C: |\lambda|=1\}$.  
$\Gamma$ is called the 
{\it gauge group} of $H$.  Alternately, one may 
think of the structure $H,\Gamma$ as a 
$\Bbb Z$-graded Hilbert space by considering 
the spectral subspaces $\{H_n: n\in\Bbb Z\}$ of 
$\Gamma$,
$$
H_n = \{\xi\in H: \Gamma(\lambda)\xi=\lambda^n\xi,
\quad \lambda\in\Bbb T\}.  
$$
The spectral subspaces give rise to an orthogonal 
decomposition 
$$
H=\dots \oplus H_{-1}\oplus H_0\oplus H_1\oplus\dots .
\tag{6.1}
$$
Conversely, given an orthogonal decomposition 
of a Hilbert space $H$ of the form (6.1), one can 
define an associated gauge group $\Gamma$ by
$$
\Gamma(\lambda) = \sum_{n=-\infty}^{\infty}
\lambda^n E_n
\qquad \lambda\in\Bbb T
$$
$E_n$ being the orthogonal projection onto $H_n$.  

A Hilbert $A$-module is said to be {\it graded} if 
there is given a distinguished gauge group $\Gamma$
on $H$ which is related to the canonical operators 
$T_1,\dots,T_d$ of $H$ by
$$
\Gamma(\lambda)T_k\Gamma(\lambda)^{-1} = \lambda T_k,
\qquad k=1,\dots,d,\quad \lambda\in\Bbb T.
\tag{6.2}
$$
Thus, graded Hilbert $A$-modules are those whose 
operators admit minimal (i.e., circular) symmetry.  
Letting $H_n$ be the $n$th spectral subspace of 
$\Gamma$, (6.2) implies that each operator is of 
degree one in the sense that 
$$
T_k H_n\subseteq H_{n+1}, \qquad 
k=1,\dots,d,\quad n\in\Bbb Z.  
\tag{6.3}
$$
Conversely, given a $\Bbb Z$-graded Hilbert 
space which is also an $A$-module satisfying 
(6.3), then it follows that the corresponding 
gauge group
$$
\Gamma(\lambda) = \sum_{n=-\infty}^\infty 
\lambda^n E_n
$$
satisfies (6.2), and moreover that the spectral 
projections $E_n$ of $\Gamma$ satisfy 
$T_kE_n=E_{n+1}T_k$ for $k=1,\dots,d$.  Thus it 
is equivalent to think in terms of gauge groups 
satisfying (6.2), or of $\Bbb Z$-graded Hilbert 
$A$-modules with degree-one operators satisfying 
(6.3).  Algebraists tend to prefer the latter 
description because it generalizes to fields 
other than the complex numbers.  On the other 
hand, the former description is more convenient 
for operator theory on complex Hilbert spaces, 
and in this section we work mainly with 
gauge groups and (6.2).  

Let $H$ be a graded Hilbert $A$-module.  A 
linear subspace $S\subseteq H$ is said to be 
{\it graded} if $\Gamma(\lambda)S\subseteq S$ 
for every $\lambda\in\Bbb T$.  
If $K\subseteq H$ is a graded (closed) 
submodule of $H$ then $K$ is a graded 
Hilbert $A$-module, and the gauge group of 
$K$ is of course the corresponding 
subrepresentation of $\Gamma$.  Similarly, 
the quotient $H/K$ of $H$ by a graded submodule 
$K$ is graded in an obvious way.  We require 
the following observation, asserting that several 
natural hypotheses on graded Hilbert modules
are equivalent.  

\proclaim{Proposition 6.4}
For every graded finite rank Hilbert $A$-module
$H$, the following 
are equivalent.  
\roster
\item
The spectrum of the gauge group $\Gamma$ is bounded 
below.  
\item
$H$ is pure in the sense that its associated 
completely positive map of $\Cal B(H)$ 
$\phi(A)=T_1AT_1^*+\dots+T_dAT_d^*$ satisfies
$\phi^n(\bold 1)\downarrow 0$ as $n\to\infty$.  
\item
The algebraic submodule
$$
M_H = \sp\{f\cdot\Delta \zeta: f\in A,\quad \zeta\in\Delta H\}
$$ 
is dense in $H$.  
\item
There is a finite-dimensional graded linear subspace 
$G\subset H$ which generates $H$ as a Hilbert $A$-module.  
\endroster
Moreover, if (1) through (4) are satisfied then the 
spectral subspaces of $\Gamma$, 
$$
H_n=\{\xi\in H: \Gamma(\lambda)\xi=\lambda^n\xi\}
\qquad n\in\Bbb Z,
$$ 
are all finite dimensional.  
\endproclaim
\demo{proof}
We prove that 
(1)$\implies$ (2) $\implies$(3)$\implies$(4)$\implies$(1).  
Let $E_n$ be the projection onto the $n$th spectral subspace
$H_n$ of $\Gamma$ and let 
$T_1,\dots,T_d$ be the canonical operators of $H$.  
From the commutation formula (6.2) it follows that 
$T_1T_1^*+\dots+T_dT_d^*$ commutes with 
$\Gamma(\lambda)$ and hence
$$
\Gamma(\lambda)\Delta = \Delta\Gamma(\lambda),\qquad
\lambda\in\Bbb T, 
\tag{6.5}
$$
where $\Delta = (\bold 1-T_1T_1^*-\dots-T_dT_d^*)^{1/2}$.
\demo{proof of (1)$\implies$(2)}
The hypothesis (1) implies that there is an integer
$n_0$ such that $E_n=0$ for $n<n_0$.  
By the preceding remarks we have
$T_kE_p=E_{p+1}T_k$ for every $p\in\Bbb Z$.  Thus
$\phi(E_p)=E_{p+1}\phi(\bold 1)\leq E_{p+1}$, 
and hence $\phi^n(E_p)\leq E_{p+n}$.  Writing 
$$
\phi^n(\bold 1)=\phi^n(\sum_{p=n_0}^\infty E_p) =
\sum_{p=n_0}^\infty\phi^n(E_p) \leq \sum_{p=n_0+n}^\infty E_p,
$$
the conclusion $\lim_n\phi^n(\bold 1)=0$ is apparent.  
\enddemo
\demo{proof of (2)$\implies$(3)}Assuming $H$ is pure, 
(1.13) implies that
the natural map $U_0\in\hom(H^2\otimes\Delta H,H)$
defined by $U_0(f\otimes \zeta)=f\cdot\Delta \zeta$ satisfies
$U_0U_0^*=\bold 1$, and therefore
$\overline{M_H}=U_0(H^2\otimes\Delta H)=H$.  
\enddemo
\demo{proof of (3)$\implies$(4)}
Assuming (3), notice that $G =\Delta H$ satisfies 
condition (4).  Indeed,  $G$ is finite 
dimensional because $\rank(H)<\infty$, it is graded 
because of (6.5), and it generates $H$ as a
closed $A$-module
because the $A$-module $M_H$ generated by $G$ is dense
in $H$.   
\enddemo
\demo{proof of (4)$\implies$(1)}
Let $G\subseteq H$ satisfy (4).  The restriction of 
$\Gamma$ to $G$ is a finite direct sum of irreducible 
subrepresentations, and hence there are integers 
$n_0\leq n_1$ such that 
$$
G = G_{n_0}\oplus G_{n_0+1}\oplus\dots\oplus G_{n_1}
$$
where $G_k=G\cap H_k$.  In particular, 
$G\subseteq H_{n_0}+H_{n_0+1}+\dots$.  Since the 
space $H_{n_0}+H_{n_0+1}+\dots$ is invariant under 
the operators $T_1,\dots,T_d$ by (6.3), we have
$$
H=\overline{\sp A\cdot G}
\subseteq H_{n_0}+H_{n_0+1}+\dots.  
$$ 
Thus $H=H_{n_0}+H_{n_0+1}+\dots$, hence
the spectrum of $\Gamma$ is bounded below by $n_0$.

The finite dimensionality of all of the spectral subspaces 
of $\Gamma$ follows from condition (4), together with the 
fact that for every $n=0,1,2,\dots$, 
the space $\Cal P_n$ of operators 
$\{f(T_1,\dots,T_d)\}$ where 
$f$ is a homogeneous polynomial of degree 
$n$ is finite dimensional and 
and $\Cal P_n$ maps $H_k$ into $H_{k+n}$.  
\qedd\enddemo
\enddemo

\proclaim{Theorem B}
For every finite rank graded Hilbert $A$-module
$H$ satisfying the conditions of Proposition
6.4 we have $K(H)=\chi(H)$.  
\endproclaim
\demo{proof} 
Because 
of the stability properties of 
$\chi(\cdot)$ and $K(\cdot)$ established in the corollaries
of Theorems C and D, it suffices to exhibit a closed 
submodule $H_0\subseteq H$ of finite codimension for 
which $K(H_0) = \chi(H_0)$.  $H_0$ is constructed as 
follows.  

Let $\{E_n: n\in\Bbb Z\}$ be the spectral projections 
of the gauge group
$$
\Gamma(\lambda)=\sum_{n=-\infty}^\infty \lambda^nE_n.  
$$
Since $\Delta$ is a finite rank operator in the 
commutant of $\{E_n: n\in \Bbb Z\}$, we must have 
$E_n\Delta = \Delta E_n =0$ for all but a finite 
number of $n\in\Bbb Z$, and hence there are integers 
$n_0\leq n_1$ such that 
$$
\Delta = \Delta_{n_0}+\Delta_{n_0+1}+\dots+\Delta_{n_1}
\tag{6.6}
$$
$\Delta_k$ denoting the finite rank positive operator 
$\Delta E_k$.  

We claim that for all $n\geq n_1$ we have 
$$
\phi(E_n)=E_{n+1}.  
\tag{6.7}
$$
Indeed, since $H$ is pure (Proposition 6.4 (2))
we can assert that 
$$
\bold 1_H = \sum_{p=0}^\infty \phi^p(\Delta^2)
\tag{6.8}
$$
because 
$$
\sum_{p=0}^n\phi^p(\Delta^2)=
\sum_{p=0}^n\phi^p(\bold 1_H-\phi(\bold 1_H)) =
\bold 1_H-\phi^{n+1}(\bold 1_H)
$$
converges strongly to $\bold 1_H$ as $n\to\infty$.  
Multiplying (6.8) on the left with $E_n$ we find
that
$$
E_n = \sum_{p=0}^\infty E_n\phi^p(\Delta^2), 
\qquad n\in\Bbb Z.
\tag{6.9}
$$
Using (6.6) we have 
$$
E_n\phi^p(\Delta^2) = 
\sum_{k=n_0}^{n_1}E_n\phi^p(\Delta_k^2).  
$$
Now $\Delta_k^2\leq E_k$ and hence 
$\phi^p(\Delta_k^2)\leq E_{k+p}$ for every 
$p=0,1,\dots$.  Thus for $n\geq n_1$,
$$
\sum_{p=0}^\infty\sum_{k=n_0}^{n_1}
E_n\phi^p(\Delta_k^2) =
\sum_{k=n_0}^{n_1}\phi^{n-k}(\Delta_k^2) =
\phi^{n-n_1}(\sum_{k=n_0}^{n_1}\phi^{n_1-k}(\Delta_k^2)).
$$
This shows that when $n\geq n_1$, $E_n$ has the form
$$
E_n = \phi^{n-n_1}(B), 
\tag{6.10}
$$
where $B$ is the operator 
$$
B=\sum_{k=n_0}^{n_1}\phi^{n_1-k}(\Delta_k^2), 
$$
and (6.7) follows immediately from (6.10).

Now consider the submodule $H_0\subseteq H$ defined by
$$
H_0=\sum_{n=n_1}^\infty E_nH.  
$$
Notice that $H_0^\perp$ is finite dimensional.  
Indeed, that is apparent from the fact that
$$
H_0^\perp = \sum_{n=-\infty}^{n_1-1}E_nH
$$
because by Proposition 6.4 (1) only a finite number 
of the projections $\{E_n: n<n_1\}$ can be nonzero
(indeed, here one can show that $E_n=0$ for 
$n<n_0$), and Proposition 6.4 also implies
that $E_n$ is finite dimensional for all $n$.  

Let $\phi_0: \Cal B(H_0)\to\Cal B(H_0)$
be the completely positive map of $\Cal B(H_0)$ associated
with the operators 
$T_1\restriction_{H_0},\dots,T_d\restriction_{H_0}$.  
Then for every $k=0,1,\dots$ we have 
$$
\phi_0^k(\bold 1_{H_0}) = \sum_{n=n_1}^\infty\phi^k(E_n).  
$$
From (6.7) we have $\phi^k(E_n)=E_{n+k}$ for 
$n\geq n_1$, and hence 
$$
\phi_0^k(\bold 1_{H_0}) = \sum_{p=n_1+k}^\infty E_p.  
$$
It follows that 
$$
\bold 1_{H_0}-\phi_0^{k+1}(\bold 1_{H_0}) =
E_{n_1}+E_{n_1+1}+\dots+E_{n_1+k}
$$
is a projection for every $k=0,1,\dots$.  
Thus for every $k\geq 0$,
$$
\tr (\bold 1_{H_0}-\phi_0^{k+1}(\bold 1_{H_0})) =
\rank (\bold 1_{H_0}-\phi_0^{k+1}(\bold 1_{H_0})),
$$
and the desired formula 
$K(H_0)=\chi(H_0)$ follows 
immediately from Theorems C and D 
after multiplying through by $d!/k^d$ and taking 
the limit on $k$.
\qedd
\enddemo

\subheading{7.  Degree}

Theorem C shows that the Euler characteristic 
(of a finite rank Hilbert $A$-module) vanishes 
whenever the rank function 
$\rank(\bold 1-\phi^{n+1}(\bold 1))$ grows relatively 
slowly.  In such cases there are other numerical 
invariants which must be nontrivial and 
which can be calculated explicitly in certain cases.  
In this brief section we define these secondary 
invariants and summarize their basic properties.

Let $H$ be a finite rank Hilbert $A$-module.  
Consider the algebraic submodule
$$
M_H=\sp\{f\cdot\Delta \xi: f\in A, \quad \xi\in H\}
$$
and its natural filtration $\{M_n: n=0,1,2,\dots\}$
$$
M_n=
\sp\{f\cdot\Delta \xi: \deg f\leq n, \quad \xi\in H\}.
$$
By Theorem 3.4 
there are integers 
$c_0,c_1,\dots,c_d$ such that 
$$
\dim M_n = c_0q_0(n)+c_1q_1(n)+\dots+c_dq_d(n)
\tag{7.1}
$$
for sufficiently large $n$.  Let 
$k$ be the degree of the polynomial on the 
right of (7.1).  We observe first that 
the pair $(k,c_k)$ depends only on the algebraic
structure of $M_H$.  

\proclaim{Proposition 7.2}
Let $M$ be a finitely generated $A$-module,  
let $\{M_n: n\geq 1\}$ be a proper filtration 
of $M$, and suppose $M\neq \{0\}$.  
Then there is a unique integer $k$, 
$0\leq k\leq d$, such that the limit 
$$
\mu(M)=k!\lim_{n\to\infty}\frac{\dim M_n}{n^k}
$$
exists and is nonzero.  $\mu(M)$ is a positive
integer and the pair $(k,\mu(M))$ 
does not depend on the particular filtration $\{M_n\}$.  
\endproclaim
\demo{proof}
By Theorem 3.4 
there are integers 
$c_0,c_1,\dots,c_d$ such that 
$$
\dim M_n = c_0q_0(n)+c_1q_1(n)+\dots+c_dq_d(n)
$$
for sufficiently large $n$.  Let 
$k$ be the degree of the polynomial on the 
right.  Noting that $q_r(x)$ is a polynomial 
of degree $r$ with leading coefficient 
$1/r!$, it is clear that this $k$ is the 
unique integer with the stated property
and that 
$$
\mu=k!\lim_{n\to\infty}\frac{\dim M_n}{n^k}=c_k
$$
is a (necessarily positive) integer.  

To see that $(k,\mu)$ does not depend on the 
filtration, let $\{M_n^\prime\}$ be a second 
proper filtration.  $\{M_n^\prime\}$ gives 
rise to a polynomial $p^\prime(x)$ of degree
$k^\prime$ which satisfies 
$\dim M_n^\prime = p^\prime(n)$ for sufficiently
large $n$.  As in the proof of Proposition 
3.7, there is an integer $p$ such that 
$\dim M_n\leq \dim M_{n+p}^\prime$ for 
sufficiently large $n$.  Thus 
$$
0<k!\lim_{n\to\infty}\frac{\dim M_n}{n^k} \leq
k!\limsup_{n\to\infty}\frac{\dim M_{n+p}^\prime}{n^k}.  
$$
Now if $k$ were greater than $k^\prime$ then the 
term on the right would be $0$.  Hence $k\leq k^\prime$
and, by symmetry, $k = k^\prime$.  

We may now argue exactly as in the proof of 
Proposition 3.7 to conclude that the leading 
coefficients of the two polynomials must 
be the same, hence $\mu=\mu^\prime$.  
\qedd
\enddemo

\proclaim{Definition 7.3}
Let $H$ be a Hilbert $A$-module of finite positive
rank.  The degree of the polynomial (7.1) associated 
with any proper filtration of the algebraic module 
$M_H$ is called the degree of $H$, and is written
$\deg(H)$.  
\endproclaim

We will also write $\mu(H)$ for the integer
$$
\mu(H) = 
\deg(H)! \lim_{n\to\infty}\frac{\dim M_n}{n^{\deg(H)}}
$$
associated with the degree of $H$.  
If $M_H$ is finite dimensional and not 
$\{0\}$ then the 
sequence of dimensions $\dim M_n$ associated 
with any proper filtration $\{M_n\}$ is eventually 
a nonzero constant, hence $\deg(H)=0$ and 
$\mu(H)=\dim(H)$; conversely, if 
$\deg (H)=0$ then $M_H$ is finite dimensional.  
In particular, $\deg{H}$ {\it is a positive integer 
satisfying} $\deg(H)\leq d$
{\it whenever the algebraic submodule}
$M_H$ {\it is infinite dimensional}.  

Note too that $\deg(H)=d$ iff the Euler characteristic
is positive, and in that case we have 
$\mu(H)=\chi(H)$.  In general, there is no obvious 
relation between $\deg(H)$ and $\rank(H)$, or 
between $\mu(H)$ and $\rank(H)$.  In particular, 
$\mu(H)$ can be arbitrarily large.  The operator-theoretic 
significance of the invariant $\mu(H)$ is not yet 
understood.  An example for which 
$1<\deg(H)<d$ is worked out in section 9.

Finally, let $\phi$ be the completely positive map associated
with the canonical operators $T_1,\dots,T_d$, 
$$
\phi(A) = T_1AT_1^*+\dots+T_dAT_d^*, \qquad A\in\Cal B(H), 
$$
and consider the generating function of the of the 
sequence of integers 
$\rank(\bold 1-\phi^{n+1}(\bold 1))$, $n=0,1,2,\dots$, defined
as the formal power series
$$
\hat\phi(t) = 
\sum_{n=0}^\infty \rank(\bold 1-\phi^{n+1}(\bold 1))t^n.
\tag{7.4}
$$
We require the following description of 
$\deg(H)$ and $\mu(H)$ in terms of $\hat\phi(t)$.

\proclaim{Proposition 7.5}
The series $\hat\phi(t)$ converges for every 
$t$ in the open unit disk of the complex plane.  
There is a 
polynomial $p(t) = a_0+a_1t+\dots+a_st^s$ and a 
sequence $c_0, c_1,\dots,c_d$ of real numbers, 
not all of which are $0$,  such that 
$$
\hat\phi(t) = p(t) + \frac{c_0}{1-t} + \frac{c_1}{(1-t)^2} +
\dots + \frac{c_d}{(1-t)^{d+1}}, \qquad |t|<1.  
$$
This decomposition is unique, and $c_k$ belongs to 
$\Bbb Z$ for every $k=0,1,\dots, d$.  $\deg(H)$ is 
the largest $k$ for which $c_k\neq 0$, and $\mu(H)=c_k$.  
\endproclaim
\demo{proof}
The proof of Theorem C shows that 
$\rank(\bold 1-\phi^{n+1}(\bold 1))=\dim M_n$,
where $\{M_n: n=1,2,\dots\}$ is the natural filtration 
of $M_H$,
$$
M_n=
\sp\{f\cdot \xi: \deg(f)\leq n,\quad \xi\in\Delta H\}.
$$
Since each $q_r(x)$ is a polynomial 
of degree $r$, formula (7.1) implies that there is 
a constant $K>0$ such that
$$
\dim M_n\leq Kn^d, \qquad n=1,2,\dots,
$$
and this estimate implies that the power series
$\sum_n \dim M_nt^n$ converges absolutely for every 
complex number $t$ in the open unit disk.  

Note too that for every $k=0,1,\dots,d$ the generating
function for the sequence $q_k(n)$, $n=0,1,\dots$ is 
given by 
$$
\hat{q_k}(t) = \sum_{n=0}^\infty q_k(n) t^n = 
(1-t)^{-k-1},
\qquad |t|<1.  
\tag{7.6}
$$
Indeed, the formula is obvious for $k=0$ since $q_0(n)=1$
for every $n$; and 
for positive $k$ the recurrence formula 3.2.2, together 
with $q_k(0)=1$, implies that 
$$
(1-t)\hat{q_k}(t) = \hat{q_{k-1}}(t), 
$$
from which (7.6) follows immediately.

Using (7.1) and (7.6) we find that 
there is a polynomial $f(x)$ such that 
$$
\hat\phi(t) = f(t) + \sum_{k=0}^d\frac{c_k}{(1-t)^{k+1}},
\tag{7.7}
$$
as asserted.  

(7.7) implies that $\hat\phi$ extends to a meromorphic
function in the entire complex plane, having a single 
pole at $t=1$.  The uniqueness of the representation 
of (7.7) follows from the uniqueness of 
the Laurent expansion of an analytic function around
a pole.  The remaining assertions of Propostion 7.5 are now 
obvious from the relation that exists between (7.1) 
and (7.6).  
\qedd
\enddemo

\subheading{8. Applications: inner sequences and 
graded ideals}

Let $M$ be a closed submodule of the free Hilbert 
module $H^2=H^2(\Bbb C^d)$ in dimension $d\geq 2$.  
Theorem 5.9 implies that there is a (finite or 
infinite) sequence of multipliers $\phi_n$ of 
$H^2$ whose associated multiplication operators 
$M_{\phi_n}$ satisfy 
$$
\sum_n M_{\phi_n}M_{\phi_n}^* = P_M, 
\tag{8.1}
$$
$P_M$ denoting the orthogonal projection of 
$H^2$ on $M$.  

\remark{Remarks}
If $\Phi=\{\phi_1,\phi_2,\dots\}$ is a finite set then
we require that it be linearly independent.  If $\Phi$ 
is infinite then the appropriate condition is that for 
every sequence $\lambda=(\lambda_n)\in\ell^2$ we have 
$$
\sum_{n=1}^\infty \lambda_n\phi_n = 0 \implies
\lambda_1=\lambda_2=\dots=0.
\tag{8.2}
$$
It is best to think of conditions (8.1) and (8.2) in
terms of dilation theory.  (8.1) asserts that the 
operator $A\in\hom(H^2\otimes\ell^2,H^2)$ defined 
by $A(f\otimes \lambda)=f\cdot\sum_n\lambda_n\phi_n$
gives rise to a dilation of the pure Hilbert 
$A$-module $M$
$$
H^2\otimes\ell^2\underset A\to\longrightarrow M\longrightarrow 0
$$
and (8.2) asserts that it is the minimal dilation of 
$M$ in that condition (4) of Proposition 1.9 is satisfied.  

The uniqueness of minimal dilations leads to the following 
description of all possible sequences 
$\Psi=\{\psi_1,\psi_2,\dots\}$ which also satisfy (8.1)
and (8.2). For definiteness, suppose that $\Phi$ is infinite. 
Then $\Psi$ must also be infinite and it is related to 
$\Phi$ through an infinite unitary matrix of scalars 
$(\lambda_{ij})$ as follows,
$$
\psi_k=\sum_{j=1}^\infty \lambda_{kj}\phi_j, \qquad k=1,2,\dots.  
$$
If $\Phi$ is a finite set with $N$ elements then so is 
$\Psi$, and the connecting matrix $(\lambda_{ij})$
in this 
case is an element of $\Cal U_N(\Bbb C)$.  Since 
we will obtain stronger results for graded submodules 
later in the section, we omit these details.  
\endremark

We have seen in (4.15) that any sequence of multipliers 
$\Phi=\{\phi_1,\phi_2,\dots\}$ satisfying (8.1) obeys 
$\sum_n|\phi(z)|^2\leq 1$ for every $z\in B_d$, 
and hence the associated sequence of boundary functions
$\tilde\phi_n: \partial B_d\to \Bbb C$ satisfies 
$\sum_n|\tilde\phi_n(z)|^2\leq 1$ almost everywhere $d\sigma$
on the boundary $\partial B_d$.  $\Phi$ is called an 
{\it inner sequence} if equality holds
$$
\sum_n|\tilde \phi_n(z)|^2=1
$$
almost everywhere $(d\sigma)$ on $\partial B_d$.  

We do not know if every nonzero closed submodule 
$M\subseteq H^2$ is associated with an inner 
sequence via (8.1); but the following result covers many 
cases of interest.  

\proclaim{Theorem E}
Let $M$ be a closed submodule of $H^2$ which contains a 
nonzero polynomial.  Then every sequence 
$\Phi=\{\phi_1,\phi_2,\dots\}$ satisfying (8.1) is an 
inner sequence.  
\endproclaim

\demo{proof}
Consider the rank-one Hilbert module $H=H^2/M$.  The 
natural projection $U: H^2\to H^2/M$ provides a minimal
dilation 
$$
H^2\underset U\to\longrightarrow H^2/M\longrightarrow 0
$$
hence the algebraic submodule of $H$ is given by 
$$
M_H=U(A) = (A+M)/M \cong A/A\cap M.  
$$
Thus the annihilator of $M_H$ is $A\cap M\neq \{0\}$.  
A theorem of Auslander and Buchsbaum (Corollary 20.13 of
\cite{14}, or Theorem 195 of \cite{18}) implies that $\chi(M_H)=0$.  
By Corollary 2 of Theorem C we have 
$K(H)=\chi(H)=\chi(M_H)=0$, and the assertion follows 
from Theorem 4.16.
\qedd
\enddemo

We turn now to the discussion of graded ideals in 
$A$ and graded submodules of $H^2$.  
Consider the free Hilbert module $H^2$ as a graded 
module with gauge group
$$
(\Gamma_0(\lambda)f)(z_1,\dots,z_d)=
f(\lambda z_1,\dots,\lambda z_d) \qquad \lambda\in \Bbb T,
\quad f\in H^2.  
$$
Every closed graded submodule $M\subseteq H^2$ decomposes into 
an orthogonal direct sum 
$$
M=M_0\oplus M_1\oplus M_2\oplus\dots
$$
where $M_n=\{f\in M: \Gamma_0(\lambda)f=\lambda^nf \}$ 
consists of all elements of $M$ which are homogeneous
polynomials of degree $n$.  The (nonclosed) linear
span 
$$
I=M_0 + M_1 + M_2 +\dots
$$
is a graded ideal in the polynomial algebra 
$A=\Bbb C[z_1,\dots,z_d]$ which is dense in $M$ in 
the norm of $H^2$.  This association 
$I \leftrightarrow M=\bar I$ between graded 
polynomial ideals and graded submodules of $H^2$ 
is bijective.  

These remarks show that one can apply the Hilbert 
space methods of this paper to gain information about
the structure of graded ideals in the algebra $A$.  
We base our approach to such questions on the notion 
of metric basis.  We show that metric bases are 
inner sequences, and that they are necessarily 
infinite whenever the ideal is of infinite 
codimension in $A$.  
This has direct implications for
infinite dimensional pure 
rank-one graded Hilbert modules $H$: such 
a Hilbert module $H$ is either
free (i.e., $H\cong H^2$) or its minimal resolution
into free Hilbert modules
$$
\dots\longrightarrow F_2\longrightarrow F_1
\longrightarrow H\longrightarrow 0
$$
becomes infinite at $F_2$ in that 
$F_2\cong H^2\otimes \ell^2$ is a free Hilbert 
module of infinite rank.  

Let $I$ be a graded ideal in $A=\Bbb C[z_1,\dots,z_d]$.  
By a {\it generator} for $I$ we mean a (finite or infinite)
set $\Phi$ of 
homogeneous polynomials of $I$ 
which is linearly independent (i.e., every finite subset
is linearly independent) and which generates $I$ in 
the sense that every element $g\in I$ can be written 
as a finite linear combination
$$
g=f_1\cdot\phi_1+f_1\cdot\phi_2+\dots+f_r\cdot\phi_r. 
\tag{8.3} 
$$
where $f_1,f_2,\dots\in A$ and 
$\phi_1,\phi_2,\dots,\phi_r\in\Phi$.  
Of course, Hilbert's basis theorem implies that 
every ideal of $A$ has a finite generator; 
but since we want to relate generators closely
to the natural norm on $A$, we must work
with infinite generators.  

A generator $\Phi$ is said to be {\it contractive}
if for every finite set $\phi_1,\dots,\phi_r$ of 
distinct elements of $\Phi$ and $f_1,\dots,f_r\in A$
we have
$$
\|\phi_1\cdot f_1+\dots+\phi_r\cdot f_r\|^2\leq 
\|f_1\|^2+\dots+\|f_r\|^2.  
\tag{8.4}
$$
Any generator can be made into a contractive
one by scaling down its individual members appropriately.  

\proclaim{Definition 8.5}
Let $I$ be a nonzero graded ideal in $A$.  A metric 
basis for $I$ is a contractive generator with the following
additional property; every polynomial $g\in I$ of degree
$n$ admits a representation of the the form (8.3) where
$\phi_1,\dots,\phi_r$ are distinct elements of $\Phi$ 
and where 
$f_1,\dots,f_r,\phi_1,\dots,\phi_r$ together satisfy
\roster
\item
$\deg \phi_k + \deg f_k \leq n$, $k=1,\dots,r$ and 
\item
$\|g\|^2 = \|f_1\|^2+\dots +\|f_r\|^2$.  
\endroster
\endproclaim

\remark{Remarks}
Property (1) implies that the degrees of 
$\phi_1,\dots,\phi_r$ are as small as they could 
possibly be, and (2) asserts that the norm
of the $r$-tuple $(f_1,\dots,f_r)$ is as small as the 
contractivity hypothesis allows.  
The price one has to pay for these
two favorable conditions is 
that metric bases are typically infinite
(see Theorem F below), and in such 
cases the lengths $r$ of the expressions 
appearing in (8.3)
are unbounded as $g$ varies.  
Nevertheless, if one uses a metric 
basis $\Phi$ to represent 
polynomials $g$ of degree at most 
$n$ as in (8.3),
then one has good control over the length of such 
expressions in terms of $n$.  More precisely, 
for every polynomial $g$ of degree at most 
$n$, there is a minimum length $r=r(g)$ of 
expressions of the form (8.3).  In remark 
8.10 we exhibit explicit upper bounds for 
the sequence $N_1, N_2,\dots$ defined by
$$
N_n=\max\{r(g): g\in A,\quad\deg g\leq n\}
$$
in terms of $n$ and the defect operator 
$\Delta$ of the Hilbert $A$-module 
$\bar I\subseteq H^2$.  
\endremark

Another feature of metric bases is that they 
can be written down quite explicitly, as 
we now show by exhibiting a metric basis 
for an arbitrary graded ideal $I\neq\{0\}$. 
Let $M=\bar I$ be the corresponding
graded submodule of $H^2$, and let $T_1,\dots,T_d$ be the 
canonical operators of $M$, $T_kf = z_kf$, $f\in M$.  Let 
$\Delta = (\bold 1_M-T_1T_1^*-\dots-T_dT_d^*)^{1/2}$ be 
the defect operator of the Hilbert $A$-module $M$, and 
let $G$ be the Hilbert space $G=\overline{\Delta M}$.  
The restriction of $\Delta$ to $G$ is a positive contraction
with trivial kernel.  Since
the gauge group $\Gamma_0$ of $H^2$ leaves $M$ invariant and
commutes with $\Delta$, it follows that 
$$
\Gamma_G(\lambda) = \Gamma_0(\lambda)\restriction_G,
\qquad \lambda\in\Bbb T
$$
defines a grading on $G$, and we obtain a decomposition 
$$
\Delta\restriction_G=\Delta_0+\Delta_1+\Delta_2+\dots
$$
where $\Delta_n=\Delta P_{G_n}$ is a positive finite 
rank operator whose restriction to the space 
$G_n\subseteq G$ of homogeneous polynomials of degree 
$n$ in $G$ is positive definite.  For each $n=0,1,2,\dots$
let $\Psi_n$ be an orthonormal basis for $G_n$ consisting 
of {\it eigenvectors} of $\Delta$.  
Note that every element $\psi\in\Psi$ is associated with 
a positive eigenvalue $\lambda$ of $\Delta$, $\Delta \psi=\lambda\psi$.  
We define $\Phi_n=\Delta\Psi_n$, and 
$$
\Phi=\Phi_0\cup\Phi_1\cup\Phi_2\cup\dots =\Delta\Psi.
\tag{8.6}
$$
By construction, the elements of $\Phi$ are nonzero mutually 
orthogonal homogeneous polynomials, and in particular $\Phi$ 
is a linearly independent set.  

\proclaim{Proposition 8.7}
The set $\Phi$ defined by (8.6) is a metric basis for $I$.  
\endproclaim
\demo{proof}
The proof is a direct application of the dilation theory 
summarized in section 1.  Consider the free Hilbert 
module $F=H^2\otimes G=H^2\otimes\overline{\Delta M}$, and 
define the operator $U\in\hom(F,M)$ by
$$
U(f\otimes g)=f\cdot\Delta g, \qquad f\in A\quad g\in G.  
$$
Because $M$ is a pure Hilbert $A$-module, $U$ is a coisometry.

Since $M\subseteq H^2$ we can think of $U$ as an element 
of $\hom(F,H^2)$, and then we have 
$UU^* = P_M$, $P_M$ denoting the projection of $H^2$ onto $M$.  
We can make $F=H^2\otimes G$ into a graded Hilbert module 
by introducing the gauge group
$$
\Gamma_F(\lambda)=\Gamma_0(\lambda)\otimes\Gamma_G(\lambda),
\qquad \lambda\in\Bbb T.  
$$
Notice that $U$ is a graded operator in the sense that 
$\Gamma_0(\lambda)U=U\Gamma_F(\lambda)$, $\lambda\in\Bbb T$.  
Indeed, for $f\in A$ and $g\in G$ 
$$
\align
U(\Gamma_F(\lambda)(f\otimes g))&=
U(\Gamma_0(\lambda)f\otimes\Gamma_0(\lambda)g)=
\Gamma_0(\lambda)f\cdot\Delta\Gamma_0(\lambda)g \\
&=
\Gamma_0(\lambda)f\cdot (\Gamma_0(\lambda)\Delta g)=
\Gamma_0(\lambda)(f\cdot \Delta g)=
\Gamma_0(\lambda)U(f\otimes g),   
\endalign
$$
and the assertion follows because $H^2\otimes G$ is 
spanned by such elements $f\otimes g$.  Letting 
$\{E_n: n= 0,1,\dots\}$ and $\{\tilde E_n: n=0,1,\dots\}$
be the spectral projections of $\Gamma_0$ and $\Gamma_F$,
$$
\Gamma_0(\lambda)=\sum_{n=0}^\infty \lambda^nE_n,
\qquad
\Gamma_F(\lambda)=\sum_{n=0}^\infty \lambda^n\tilde E_n
$$
then we have 
$$
U\tilde E_n = E_n U, \qquad n=0,1,2,\dots.  
\tag{8.8}
$$

We claim now that $\Phi$ satisfies the contractivity 
condition (8.4).  
For that, let $\phi_1,\dots,\phi_r$ be distinct 
elements of $\Phi$.  Then there are positive 
numbers $\lambda_1,\dots,\lambda_r$ and an 
orthonormal set $\psi_1,\dots,\psi_r\in\Psi$
such that $\phi_k=\Delta\psi_k=\lambda_k\psi_k$,
$k=1,\dots,r$.  Thus for $f_1,\dots,f_r\in A$ we 
have 
$$
\sum_{k=1}^r f_k\cdot\phi_k = 
\sum_{k=1}^r f_k\cdot\Delta\psi_k =
U(\sum_{k=1}^r f_k\otimes \psi_k).  
$$
Since $\|U\|\leq 1$ and 
$\psi_1,\dots,\psi_r$ is an orthonormal 
set in $G$ we have
$$
\|\sum_{k=1}^r f_k\cdot\phi_k\|^2 \leq 
\|\sum_{k=1}^r f_k\otimes \psi_k\|^2 =
\sum_{k=1}^r\|f_k\|^2,
$$
and the assertion follows.  

It remains to show that for every 
polynomial $p\in I$
with $\deg p=n$, there are elements $\phi_1,\dots,\phi_r\in\Phi$
and $f_1,\dots,f_r\in A$ satisfying 
$p=f_1\cdot\phi_1+\dots+f_r\cdot\phi_r$
along with conditions (1) and (2) of Definition 8.5.  
Fixing $p$, define an element $\xi\in H^2\otimes G$ by
$$
\xi = (\tilde E_0+\tilde E_1+\dots+\tilde E_n)U^*p.  
$$
Note first that $U\xi=p$. Indeed, since $UU^*=P_M$ we 
have $UU^*p=p$ and using (8.8) we have 
$$
U\xi = U((\sum_{k=0}^n\tilde E_k)U^*p)=
(\sum_{k=0}^n E_k)UU^*p = \sum_{k=0}^n E_kp=p.  
$$

It is also clear from its definition that $\xi$ belongs
to $\tilde E_0F + \tilde E_1 F+\dots \tilde E_n F$.  Since 
the gauge group of $F$ is a tensor product 
$\Gamma_F(\lambda)=\Gamma_0(\lambda)\otimes\Gamma_G(\lambda)$, 
we have 
$$
\tilde E_kF = \sum_{p=0}^k H^2_p\otimes G_{k-p}
$$
and hence 
$$
\tilde E_0F + \tilde E_1 F+\dots \tilde E_n F
\subseteq \sum\{H^2_p\otimes G_q: p,q\geq 0,\quad p+q\leq n\}.  
$$
Noting that 
$H^2_p\otimes G_q$ is spanned by elements of the 
form $f\otimes \psi$ with $f\in H^2_p$ and $\psi\in\Psi_q$,
it follows that $\xi$ has a decomposition 
of the form
$$
\xi=f_1\otimes\psi_1+f_2\otimes\psi_2+\dots+f_r\otimes \psi_r
\tag{8.9}
$$
where $\psi_1,\dots,\psi_r$ are distinct elements of 
$\Psi$ and where $\deg f_k+\deg\psi_k\leq n$.  
Setting $\phi_k=\Delta\psi_k\in\Phi$ 
we have $\deg\phi_k=\deg\psi_k$
and 
$$
p=U\xi=\sum_{k=1}^rf_k\cdot\Delta\psi_k=\sum_{k=1}^rf_k\cdot\phi_k.  
$$

Finally, we claim that $\|p\|^2=\|f_1\|^2+\dots+\|f_r\|^2$.  
Indeed, since $\psi_1,\dots,\psi_r$ are orthonormal in $G$, 
(8.9) implies that 
$$
\|f_1\|^2+\dots+\|f_r\|^2=\|\xi\|^2=
\|(\tilde E_0+\dots+\tilde E_n)U^*p\|^2\leq 
\|U^*p\|^2\leq \|p\|^2, 
$$
while since $\Phi$ satisfies (8.4) we must have the 
opposite inequality 
$$
\|f_1\|^2+\dots+\|f_r\|^2\geq \|p\|^2,
$$ 
and the proof is complete.
\qedd
\enddemo

\remark{Remark 8.10}Let $I$ be a graded ideal in $A$ and 
let $g\in I$ be a polynomial of degree at most $n$.  
We can estimate the length
$r$ of a representation $g=
f_1\cdot\phi_1+\dots+f_r\cdot\phi_r$
satisfying the conditions of Definition 8.5 as follows.  
Let $\Delta$ be the defect operator of 
$M=\bar I$ and let $\Delta_k$ denote its restriction 
to the subspace $M_k=I_k$ of homogeneous polynomials 
of degree $k$.  Since by property (1) of Definition 
8.5 the elements $\phi_1,\dots,\phi_r$ must all 
belong to $\Phi_0\cup\Phi_1\cup\dots\cup\Phi_n$ and since
the cardinality of $\Phi_k$ is the rank of $\Delta_k$
we conclude that 
$$
r\leq \rank \Delta_0+\dots+\rank \Delta_n.  
$$
\endremark

We now take up the issue of uniqueness of metric bases.  
Let $S=\{\xi_1,\dots,\xi_m\}$ and 
$T=\{\eta_1,\dots,\eta_m\}$ be two finite linearly 
independent sets of vectors in a Hilbert space $H$.  
We say that $S\sim T$ if $m=n$ and there is a unitary
matrix $(u_{ij})\in M_n(\Bbb C)$ such that 
$$
\eta_i = \sum_{j=1}^n u_{ij}\xi_j, \qquad 1=1,\dots,n. 
\tag{8.11} 
$$
More specifically, if 
$\Phi$ and $\Psi$ are two linearly independent 
sets of homogeneous polynomials in $H^2$ then we say 
that $\Phi$ and $\Psi$ are equivalent 
(written $\Phi\sim\Psi$) if $\Phi_n\sim\Psi_n$ for 
every $n=0,1,2,\dots$ where $\Phi_n$ (resp. $\Psi_n$)
denotes the set of all elements of $\Phi$ 
(resp. $\Psi$) which are homogeneous
of degree $n$.  We now show that any two metric bases 
for a graded ideal $I$ of $A$ are equivalent.  That 
result is based on the following.  We use the notation
$\xi\otimes\bar\eta$ to denote the rank-one operator 
$$
\zeta\mapsto \<\zeta,\eta\>\xi
$$
associated with a pair of vectors $\xi,\eta$ in a 
Hilbert space.  

\proclaim{Lemma 8.12}
Let $S=\{\xi_1,\dots,\xi_m\}$ and $T=\{\eta_1,\dots,\eta_n\}$
be linearly independent sets of vectors in a Hilbert space
$H$ such that
$$
\sum_{k=1}^m\xi_k\otimes \bar\xi_k =
\sum_{j=1}^n\eta_j\otimes \bar\eta_j.  
$$
Then $S\sim T$ in the sense of (8.11).  
\endproclaim
\demo{proof}
Consider the linear map $A:\Bbb C^m\to H$ defined by 
$$
Az = z_1\xi_1+\dots+z_m\xi_m.  
$$
After endowing $\Bbb C^m$ with its usual inner product, 
one finds that the adjoint $A^*: H\to\Bbb C^m$ is given 
by
$$
A^*\zeta= (\<\zeta,\xi_1\>,\dots,\<\zeta,\xi_m\>),
$$
and therefore
$$
AA^*=\sum_{k=1}^m \xi_k\otimes\bar \xi_k.  
$$
Because of linear independence we have 
$\ker A=\{0\}$, and hence $A^*H=\Bbb C^m$.  

By the polar decomposition, there is a unique isometry
$U:\Bbb C^m\to H$ having range $A\Bbb C^m$ such that 
$A=PU$, where $P=(AA^*)^{1/2}$.  Doing the same with 
the set $T$ we obtain a similar operator $B: \Bbb C^n\to H$
having trivial kernel such that $BB^*=AA^* = P^2$, and 
an isometry $V: \Bbb C^n\to H$ having 
range $B\Bbb C^n$ such that $B=PV$.  

Since $A\Bbb C^m=B\Bbb C^n=PH$, the operator 
$W=U^*V\in\Cal B(\Bbb C^n,\Bbb C^m)$ is unitary
and satisfies $AW=B$.  Thus $m=n$ and, letting 
$(w_{ij})$ be the matrix of $W$ relative to the usual
orthonormal basis $e_1,\dots,e_n$ for $\Bbb C^n$ we 
find that 
$$
\eta_i=Be_i=AWe_i=\sum_{j=1}^nw_{ij}Ae_j =
\sum_{j=1}^nw_{ij}\xi_j,
$$
as required.
\qedd
\enddemo

\remark{Remarks}
We point out that the converse of Lemma 8.12 is true
in the sense that if $\{\xi_1,\dots,\xi_n\}$ is 
any set of vectors in a Hilbert space and $\eta_1,\dots,\eta_n$
is related to $\{\xi_1,\dots,\xi_n\}$ by way of a unitary 
$n\times n$ matrix as in (8.11), then one 
verifies by direct computation that 
$$
\sum_{k=1}^n\xi_k\otimes\bar\xi_k=
\sum_{k=1}^n\eta_k\otimes\bar\eta_k.  
$$

Let $P$ be a finite rank positive operator on a Hilbert 
space $H$.  There is an orthonormal basis
$\{e_1,\dots,e_r\}$ for the 
subspace $PH$ consisting of eigenvectors for $P$, 
$Pe_j=\lambda_je_j$, $j=1,\dots,r$, and one has 
$\lambda_j>0$ for every $j$.  Setting 
$\xi_k=\sqrt{\lambda_k}e_k$ we have 
$P=\xi_1\otimes\bar \xi_1+\dots+\xi_r\otimes\bar \xi_r$,
and of course $\{\xi_1,\dots,\xi_r\}$ is linearly 
independent.  If $\eta_1,\dots,\eta_s$ is any other 
linearly independent set which represents $P$ in a
similar way 
$P=\eta_1\otimes\bar\eta_1+\dots+\eta_s\otimes\bar\eta_s$
then the $\eta_k$ do not have to be eigenvectors
of $P$, but Lemma 8.12 implies that $s=r$ and that 
$\{\eta_1,\dots,\eta_r\}$ is obtained from the 
eigenvector sequence $\{\xi_1,\dots,\xi_r\}$ by a 
simple rotation as in (8.11). 

We also point out that the proof of 
Lemma 8.12 is closely related 
to the construction of the metric operator space 
associated with a normal completely positive map 
$\phi: \Cal B(H)\to\Cal B(H)$ 
and in fact Lemma 8.12 extends in an appropriate way 
(with a proof similar to the one given)
to infinite sets of vectors \cite{2}.   
\endremark

In view of Lemma 8.12, the uniqueness of metric bases 
will follow from
\proclaim{Proposition 8.13}
Let $\Psi$ be a metric basis for a nonzero graded ideal 
$I$ in $A$, let 
$\Delta = (\bold 1-T_1T_1^*-\dots-T_dT_d^*)^{1/2}$
be the defect operator of $M=\bar I$, and let 
$E_n$ be the projection of $H^2$ 
onto the space of homogeneous 
polynomials of degree $n$.  Then 
$\Delta E_n=0$ iff $\Psi_n=\emptyset$, and
when $\Psi_n=\{\psi_1,\dots,\psi_s\}\neq \emptyset$
we have
$$
\sum_{k=1}^s\psi_k\otimes\bar \psi_k = \Delta ^2E_n.  
$$ 
\endproclaim
\demo{proof}
The idea of the proof is similar to that of 
Proposition 8.7, in that one constructs a graded minimal 
dilation of the Hilbert $A$-module 
$\bar I\subseteq H^2$ using the metric basis $\Psi$, 
and then simply verifies that 
the required property holds.  However, the details
are different enough that we include them.  

Let $G$ be the graded Hilbert space $\ell^2(\Psi)$
of all square summable functions $f: \Psi\to\Bbb C$
with gauge group $\Gamma_G$
$$
(\Gamma_G(\lambda)(f))(\psi)=\lambda^{\deg \psi}f(\psi),
\qquad \psi\in\Psi.  
$$
We identify the delta function 
concentrated at an element $\psi\in\Psi$ with 
$\psi$ itself, so that 
$\Psi$ becomes an orthonormal basis 
for $G$ consisting of eigenvectors of $\Gamma_G$. 

We will work with the graded free Hilbert module
$F=H^2\otimes G$, $\Gamma_F=\Gamma_0\otimes\Gamma_G$.  
Since $\Psi$ is an orthonormal basis for $G$, every 
element $\xi\in F$ has unique expansion
$$
\xi=\sum_{\psi\in\Psi} f_\psi\otimes\psi
$$
where $f_\psi$, $\psi\in\Psi$, is a sequence of 
elements of $H^2$ satisfying 
$$
\|\xi\|^2=\sum_{\psi\in\Psi} \|f_\psi\|^2.  
\tag{8.14}
$$

If $f_1,\dots,f_r$ is a set of polynomials
in $A$ and $\psi_1,\dots,\psi_r$ are distinct
elements of $\Psi$, then by (8.4) we have 
$$
\|f_1\cdot\psi_1+\dots+f_r\cdot\psi_r\|^2 \leq 
\|f_1\|^2+\dots+\|f_r\|^2 = 
\|f_1\otimes\psi_1+\dots+f_r\otimes\psi_r\|^2_F.  
$$
It follows that there is a unique bounded operator 
$U: F\to \bar I$ satisfying $U(f\otimes \psi)=f\cdot \psi$
for every $f\in A$, $\psi\in\Psi$, and of course 
$\|U\|\leq 1$ and $U\in\hom(F,\bar I)$.  

Note that 
$$
H^2\otimes G\underset U\to\longrightarrow \bar I\longrightarrow 0
$$ 
is a graded minimal dilation.  That is to say, 
$U\Gamma_F(\lambda)=\Gamma_0(\lambda)U$ and 
$U$ is a coisometry: $UU^*=\bold 1_{\bar I}$.  Indeed, 
one sees that $U$ is graded exactly as in the 
proof of Proposition 8.7, and the conditions of Definition
8.5 imply that $U$ is a partial isometry with dense range.  
In more detail, the range of $U$ obviously contains 
$I$, and for a given element $g\in I$ we can find a 
representation 
$$
g=f_1\cdot\psi_1+\dots+f_r\cdot\psi_r
$$
which satisfies the conditions (1) and (2) of Definition
8.5.  Setting 
$$
\xi=f_1\otimes\psi_1+\dots+f_r\otimes \psi_r\in F
=H^2\otimes G
$$
we find that $U\xi = g$ and 
$\|\xi\|^2=\|f_1\|^2+\dots+\|f_r\|^2=\|g\|^2$.  
Since $\|U\|\leq 1$, $U$ must be a partial isometry. 
That this dilation is minimal is clear from condition
(4) of Proposition 1.9.   

The left side of the formula of 
Proposition 8.13 is calculated as follows.  
Since $\Psi_n=\{\psi_1,\dots,\psi_s\}$
it follows that $\{\psi_1,\dots,\psi_s\}$ is an 
orthonormal basis for $G_n$.  Letting $D_n$ be the 
projection of $G$ onto $G_n$ and letting $1$ be the 
constant function of $H^2$, then 
$\{1\otimes \psi_k: k=1,\dots s\}$ is an orthonormal 
basis for the range of the projection 
$E_0\otimes D_n\in\Cal B(H^2\otimes G)$.  
Since $U(1\otimes \psi_k)=\psi_k\in \bar I$, 
we have
$$
\align
\sum_{k=1}^s\psi_k\otimes\bar\psi_k &= 
\sum_{k=1}^s U(1\otimes \psi_k)\otimes 
\overline{U(1\otimes\psi_k)}\\
&=U(\sum_{k=1}^s(1\otimes\psi_k)\otimes
\overline{(1\otimes\psi_k)})U^*=
U(E_0\otimes D_n)U^*.  
\endalign
$$
Since the $n$th spectral projection $\tilde E_n$
of $\Gamma_F$ is related to $E_0\otimes D_n$ by
$$
\tilde E_n(E_0\otimes\bold 1_G) = E_0\otimes D_n,
$$
the right side of the preceding formula is 
$$
U\tilde E_n(E_0\otimes \bold 1_G)U^* = 
E_nU(E_0\otimes \bold 1_G)U^*
$$
because $U$ is graded, $E_n$ being the projection 
of $H^2$ onto the space $H^2_n$.  We conclude 
that 
$$
\sum_{k=1}^s\psi_k\otimes\bar\psi_k=
E_nU(E_0\otimes\bold 1_G)U^*.
\tag{8.15}
$$

Recalling that $\Delta$ commutes with $E_n$, we can 
obtain the desired formula by showing that 
$$
U(E_0\otimes\bold 1_G)U^*=\Delta^2.  
\tag{8.16}
$$
Indeed, (8.16) is a general formula satisfied by 
all minimal dilations such as $U$.  Let 
$S_1,\dots,S_d$ be the $d$-shift acting on $H^2$.  
We have seen that 
$$
E_0=\bold 1_{H^2}-S_1S_1^*-\dots-S_dS_d^*.  
$$
Since the canonical operators of $F=H^2\otimes G$
are $S_1\otimes\bold 1_G,\dots,S_d\otimes\bold 1_G$, 
we have
$$
E_0\otimes \bold 1_G=
\bold 1_F-\sum_{k=1}^d 
(S_k\otimes\bold 1_G)(S_k\otimes\bold 1_G)^*.  
$$
Since $U\in\hom(F,\bar I)$ we find that 
for $T_1,\dots,T_d$ the canonical operators 
of $\bar I$ we have $U(S_k\otimes \bold 1_G) =T_kU$,
and thus  
$$
\align
U(E_0\otimes\bold 1_G)U^*&=
UU^*-T_1UU^*T_1^*-\dots -T_dUU^*T_d^* \\
&=
\bold 1_{\bar I}-T_1T_1^*-\dots T_dT_d^* =\Delta^2,
\endalign
$$
as required.   
\qedd
\enddemo

\proclaim{Corollary 1}
Any two metric bases for a graded ideal 
$I\subseteq \Bbb C[z_1,\dots,z_d]$ are equivalent.
\endproclaim
\demo{proof}
This is an immediate consequence of 
Proposition 8.13 and Lemma 8.12.
\qedd
\enddemo

\proclaim{Corollary 2}
Let $\Phi=\{\phi_1,\phi_2,\dots\}$ be a metric basis 
for a graded ideal $I\subseteq A$, and let 
$M=\bar I$ be its associated closed submodule of 
$H^2$.  Then $\Phi$ is an inner sequence 
for $M$ satisfying (8.1) and (8.2).  
\endproclaim
\demo{proof}
The minimal dilation 
$$
H^2\otimes G\underset U\to\longrightarrow M\longrightarrow 0
$$
exhibited in the proof of Proposition 8.13 shows that the 
sequence $\Phi$ satisfies (8.1) and (8.2), and Theorem E
implies that $\Phi$ is an inner sequence.  
\qedd
\enddemo

The following result implies that metric bases are 
almost always infinite.  

\proclaim{Theorem F}
Let $I$ be a proper graded ideal in $A$ having a 
finite metric basis.  
Then $I$ is of finite codimension in $A$ and 
the canonical generators $z_1,\dots,z_d$ of $A$ 
are all nilpotent modulo $I$.  
\endproclaim
\demo{proof}
Let $\Phi=\{\phi_1,\dots,\phi_n\}$ be a metric
basis for $I$ and let 
$s(z)=|\phi_1(z)|^2+\dots+|\phi_n(z)|^2$.  By 
Theorem E, $s(z)=1$ almost everywhere on the sphere 
$\partial B_d$ and since $s$ is continuous it must 
be identically $1$ on $\partial B_d$.  Thus the 
variety 
$V=\{z\in \Bbb C^d: \phi_1(z)=\dots=\phi_n(z)=0\}$
of common zeros does not intersect the unit sphere.  
$V$ cannot be empty since $I$ is proper.  
Since $V$ is a nonempty set 
invariant under multiplication by nonzero 
scalars which misses the unit sphere, 
it must consist of 
just the single point $(0,0,\dots,0)$.  

By Hilbert's Nullstellensatz there is 
an integer $p\geq 1$ such that $z_1^p,\dots,z_d^p\in I$
\cite{14}, Theorem 1.6.
Since the $A$-module $A/I$ has a cyclic vector 
$1+I$ and its canonical operators are nilpotent,
it follows that $A/I$ is finite dimensional.
\qedd
\enddemo

\proclaim{Corollary}
Let $M\neq \{0\}$ be a closed graded submodule of $H^2$ 
such that $H^2/M$ is infinite dimensional.  Then
$M$ is an infinite rank Hilbert $A$-module.  
\endproclaim
\demo{proof}
Contrapositively, assume 
that $\rank M<\infty$ and let 
$$
\Delta = (\bold 1_M-T_1T_1^*-\dots-T_dT_d^*)^{1/2}
$$
where $T_kf=z_k\cdot f$, $f\in M$. 
$\Delta$ is a positive finite rank operator
which commutes with the gauge group 
of $M$.  Thus the metric 
basis $\Phi$ for the ideal $I=M\cap A$ exhibited  
in (8.6) must be a finite set, and 
Theorem F implies that $A/I$ is 
finite dimensional.  Since $A$ is dense in 
$H^2$, the natural map 
of $A/I$ into $H^2/\bar I=H^2/M$ must have 
dense range and thus $H^2/M$ is finite dimensional,
contradicting the hypothesis on $M$.   
\qedd
\enddemo

\subheading{9. Examples}

In this section we give examples of pure rank-one
Hilbert modules illustrating (1) the failure of Theorem B
for ungraded modules, and (2) the computation of 
the degree in cases where the Euler 
characteristic vanishes.  We also give examples
of pure rank-two graded Hilbert modules 
illustrating (3) the computation of 
nonzero values of $K(H)=\chi(H)$.  

We begin with a discussion of the limits of Theorem 
B by presenting a class of examples for which 
$K(H)\neq \chi(H)$ (Proposition 9.2); a
concrete example of such a Hilbert $A$-module
is given in Example 9.3.  We then describe a natural
method for associating a graded Hilbert 
$A$-module with an algebraic  
variety in complex projective space, and we 
show that for {\it some} varieties one can 
calculate all numerical invariants of their 
associated Hilbert modules.  

\remark{Remark 9.1}
We will make use of the fact that if $K_1$ and 
$K_2$ are two closed submodules of the free 
Hilbert module $H^2$ for which $H^2/K_1$ is 
isomorphic to $H^2/K_2$, then $K_1=K_2$.  
In particular, no nontrivial quotient of 
$H^2$ of the form $H^2/K$ with $K\neq \{0\}$
can be a free Hilbert $A$-module.  

Indeed, this is part of Corollary 2 of 
\cite{1, Theorem 8.5}.  
One can also deduce it from the material summarized
in section 1 as follows.  
Let $U_j$ be the natural projection
of $H^2$ onto $H^2/K_j$, $j=1,2$.  
Then both diagrams 
$$
H^2\underset U_j\to\longrightarrow H^2/K_j\longrightarrow 0,
\qquad j=1,2
$$
define minimal dilations.  Applying Theorem 1.11, 
we see that any unitary isomorphism 
$W: H^2/K_1\to H^2/K_2$ of Hilbert $A$-modules 
gives rise to a unique unitary operator $\tilde W$ in 
the commutant of $C^*(H^2)$ such that the diagram
$$
\CD
H^2		@>>U_1>			 H^2/K_1\\
@V\tilde W VV						@VV W V \\
H^2 	@>>U_2>			H^2/K_2
\endCD
$$
commutes.  Since $C^*(H^2)$ is the 
(irreducible) Toeplitz \cstar\ 
\cite{1, Theorem 5.7}, 
$\tilde W$ must be a scalar multiple of the 
identity operator, and the assertion follows
from the fact that $\tilde W K_1=K_2$.  
\endremark

\proclaim{Proposition 9.2}
Let $K\neq \{0\}$ be a closed submodule of 
$H^2$ which contains no nonzero polynomials, 
and consider the pure rank-one module
$H=H^2/K$.  Then 
$$
0\leq K(H)<\chi(H)=1.  
$$
\endproclaim
\demo{proof}
We show first that $\chi(H)=1$ by proving that the 
algebraic submodule $M_H$ of $H$ is free.  
Let $U\in\hom(H^2,H)$ be the natural projection 
onto $H=H^2/K$.  The kernel of $U$ is $K$, and 
$U$ maps the dense linear subspace 
$A\subseteq H^2$ of polynomials onto $M_H$,
$U(A) = M_H$.  Since $A\cap K=\{0\}$, the
restriction of $U$ to $A$ gives
an isomorphism of $A$-modules 
$A\cong M_H$, and hence
$\chi(H) = \chi(A)=1$.  

On the other hand, if $K(H)$ were to equal 
$1=\rank(H)$ then by the extremal property 
(4.13) $H$ would be isomorphic to the free Hilbert
module $H^2$ of rank-one, which is impossible 
because of Remark 9.1.  
\qedd
\enddemo

\example{Problem}
We do not know if the curvature invariant $K(H)$ of 
a {\it pure} finite rank Hilbert $A$-module is 
always an integer.  Theorem B implies that this is 
the case for graded Hilbert modules, but Proposition 
9.2 shows that Theorem B does not always apply.  In 
particular, it is not known if $K(H)=0$ for the 
ungraded Hilbert modules $H$ of Prop. 9.2.  
In such cases, the 
equation $K(H)=0$ is equivalent to the existence 
of an ``inner sequence" for the invariant subspace
$K$ (see Theorem 4.16).  
\endexample

\example{Example 9.3}
It is easy to give concrete examples of submodules
$K$ of $H^2$ satisfying the hypothesis of Proposition 9.2. 
Consider, for example, the graph of the exponential 
function $G=\{(z,e^z): z\in\Bbb C\}\subseteq \Bbb C^2$.  
Take $d=2$, let $H^2=H^2(\Bbb C^2)$, and let $K$ be 
the submodule of all functions in $H^2$ which vanish 
on the intersection of $G$ with the unit ball
$$
K=\{f\in H^2: f\restriction_{G\cap B_d}=0 \}.  
$$ 
Since $f\in H^2\mapsto f(z)=\<f,u_z\>$ is 
a bounded linear functional for every $z\in B_d$ it 
follows that $K$ is closed, and it is clear that
$K\neq \{0\}$ (the function 
$f(z_1,z_2)=e^{z_1}-z_2$ belongs 
to $H^2$ and vanishes on $G\cap B_d$).  
After noting that the open unit disk about
$z=-1/2$ maps into $G\cap B_d$,
$$
\{(z,e^z): |z+1/2|<1\}\subseteq G\cap B_d
$$
an elementary argument (which we 
omit) establishes the obvious fact that no nonzero 
polynomial can vanish on $G\cap B_d$.  
\endexample

An algebraic set in complex projective space 
$\Bbb P^{d-1}$ can be described as the set of 
common zeros of a finite set of {\it homogeneous} 
polynomials $f_1,\dots, f_n\in\Bbb C[z_1,\dots,z_d]$,
$$
V=\{z\in \Bbb C^d: f_1(z)=\dots=f_n(z)=0\}
\tag{9.4}
$$
\cite{14}, pp 39--40.  One can associate 
with $V$ a graded rank-one 
Hilbert $A$-module 
in the following way.  Let $M_V$ be the 
submodule of $H^2=H^2(\Bbb C^d)$ defined by
$$
M_V = \{f\in H^2: f\restriction_{V\cap B_d}=0\}.  
$$
As in example 9.3, $M_V$ is a closed submodule 
of $H^2$, and because $\lambda V\subseteq V$ for 
complex scalars $\lambda$, $M_V$ is invariant 
under the action of the gauge group of 
$H^2$ and hence it is a graded 
submodule of $H^2$.  Thus, $H=H^2/M_V$ is a 
graded, pure, rank-one Hilbert $A$-module.  

We will show how to explicitly compute 
$H^2/M_V$ and its numerical invariants 
in certain cases, using operator-theoretic methods.  
The simplest member of this class of examples is the 
variety $V$ defined by the range of the quadratic 
polynomial 
$$
F:(x,y)\in\Bbb C^2 \mapsto 
(x^2,y^2,\sqrt2 xy)\in\Bbb C^3,   
$$
that is, 
$$
V=\{(x^2,y^2,\sqrt2 xy): x,y\in\Bbb C\}\subseteq \Bbb C^3.  
$$
However, one finds more interesting behavior
in the higher dimensional example
$$
V=\{(x^2,y^2,z^2,\sqrt2 xy,\sqrt2 xz,\sqrt2 yz): 
x,y,z\in \Bbb C\}\subseteq \Bbb C^6,
\tag{9.5}
$$
and we will discuss the example (9.5) in some detail.  
We describe a more general 
context for these examples in Remark 9.12.  

Notice first that $V$ can be described in the form 
(9.4) as the set 
$$
V=\{z\in\Bbb C^6: f_1(z)=f_2(z)=f_3(z)=f_4(z)=0\}
\tag{9.6}
$$
of common zeros of the four homogeneous polynomials 
$f_k: \Bbb C^6\to\Bbb C$,
$$
\align
f_1(z) &= z_4^2-2z_1z_2 = 0\\
f_2(z) &= z_5^2-2z_1z_3 = 0\\
f_3(z) &= z_6^2-2z_2z_3 = 0\\
f_4(z) &= z_4z_5z_6-2^{3/2}z_1z_2z_3 = 0.
\endalign
$$
The equivalence of (9.5) and (9.6) is an elementary 
computation which we omit.  Note, however, that the 
fourth equation $f_4(z)=0$ is necessary in order to 
exclude points such as 
$z=(1,1,1,-\sqrt2,\sqrt2,\sqrt2)$, which satisfy 
the first three equations 
$f_1(z)=f_2(z)=f_3(z)=0$ but which do not
belong to $V$.  Note too that $f_1,f_2, f_3$ are quadratic 
but that $f_4$ is cubic.  

We will describe the Hilbert module $H=H^2(\Bbb C^6)/M_V$ 
by identifying its associated $6$-contraction 
$(T_1,\dots,T_6)$.  These operators act on 
the {\it even} subspace $H$ of $H^2(\Bbb C^3)$, defined
as the closed linear span of all homogeneous polynomials
$f(z_1,z_2,z_3)$ of even degree $2n$, $n=0,1,2,\dots$.  
Let $S_1,S_2,S_3\in\Cal B(H^2(\Bbb C^3))$ be the 
$3$-shift.  The even subspace $H$ is not invariant 
under the $S_k$, but it is invariant under any product
of two of these operators $S_iS_j$, $1\leq i,j\leq 3$.  
Thus we can define a $6$-tuple of operators 
$T_1,\dots,T_6\in\Cal B(H)$ by 
$$
(T_1,\dots,T_6) = 
(S_1^2\restriction_H,
S_2^2\restriction_H,
S_3^2\restriction_H,
\sqrt2 S_1S_2\restriction_H,
\sqrt2 S_1S_3\restriction_H,
\sqrt2 S_2S_3\restriction_H).  
\tag{9.7}
$$
$(T_1,\dots,T_6)$ is a $6$-contraction because 
$$
\sum_{k=1}^6 T_kT_k^* = 
\sum_{i,j=1}^3S_iS_j(P_H)S_j^*S_i^*\leq \bold 1_H,  
$$
and in fact $H$ becomes a pure Hilbert 
$\Bbb C[z_1,\dots,z_6]$-module.  

If $f$ is a sum of homogeneous polynomials of even 
degrees then 
$$
\Gamma(e^{i\theta})f(z_1,z_2,z_3) =
f(e^{i\theta/2}z_1,e^{i\theta/2}z_2,e^{i\theta/2}z_3)
$$
gives a well-defined unitary action of the circle 
group on the subspace $H\subseteq H^2(\Bbb C^3)$, 
and $H$ becomes a graded Hilbert module.  

\proclaim{Proposition 9.8}
$H$ is a rank-one graded Hilbert 
$\Bbb C[z_1,\dots,z_6]$-module which is isomorphic to 
$H^2(\Bbb C^6)/M_V$.  The invariants of $H$ are 
given by $K(H)=\chi(H)=0$, $\deg(H)=4$, $\mu(H)=4$.  
\endproclaim
\demo{proof}
Let $\phi(A)=T_1AT_1^*+\dots+T_6AT_6^*$ be the 
canonical completely positive map of $\Cal B(H)$ 
and, considering $H$ as a subspace of $H^2(\Bbb C^3)$, 
let $\sigma: \Cal B(H^2)\to\Cal B(H^2)$ be the map
associated with the $3$-shift
$$
\sigma(B) = S_1BS_1^*+S_2BS_2^*+S_3BS_3^*.  
$$
$\phi$ and $\sigma$ are related in the following 
simple way: for every $A\in\Cal B(H)$ we have
$$
\phi(A) = \sum_{k=1}^6T_kAT_k^* = 
\sum_{i,j=1}^3S_iS_jAP_HS_j^*S_i^* = \sigma^2(AP_H).  
\tag{9.9}
$$
If $E_n\in\Cal B(H^2)$ denotes the projection onto 
the subspace of homogeneous polynomials of degree 
$n$, then 
$$
\phi(\bold 1_H)=\sigma^2(\sum_{n=0}^\infty E_{2n}) =
\sum_{n=0}^\infty E_{2n+2}.  
$$
It follows that 
$$
\Delta^2 = \bold 1_H-\phi(\bold 1_H) = E_0
$$
is the one-dimensional projection onto the 
space of constants.  Since 
$$
\phi^n(\bold 1_H)=\sigma^{2n}
(\sum_{p=0}^\infty E_{2p})= \sum_{p=n}^\infty E_{2p}
$$
obviously decreases to $0$ as
$n\to\infty$, we conclude that $H$ 
is a pure Hilbert module of rank one.  

Hence the minimal dilation of $H$ 
$$
H^2(\Bbb C^6)\underset U_0\to\longrightarrow 
H\longrightarrow 0
$$
is given by
$$
U_0(f) = f\cdot\Delta 1=f(T_1,\dots,T_6)\Delta 1.  
$$
If we evaluate this expression at a point
$z=(z_1,z_2,z_3)\in B_3$ we find that 
$$
U_0(f)(z_1,z_2,z_3)=
f(z_1^2,z_2^2,z_3^2,
\sqrt2 z_1z_2,\sqrt2 z_1z_3,\sqrt2 z_2z_3).  
$$
The argument on the right is a point in 
the ball $B_6$, and thus 
the preceding formula extends immediately to all 
$f\in H^2(\Bbb C^6)$.  
Notice too that $U_0$ is a {\it graded} 
morphism in that 
$U_0\Gamma_0(\lambda) = \Gamma(\lambda)U_0$, 
$\lambda\in\Bbb T$, where $\Gamma_0$ is the gauge 
group of $H^2(\Bbb C^6)$. The precding formula
shows that the kernel
of $U_0$ is $M_V$, and thus we conclude
that $H$ is isomorphic to $H^2(\Bbb C^6)/M_V$,
as asserted in Proposition 9.8.  

It remains to calculate the power series $\hat\phi(t)$
of Proposition 7.5 which determines the numerical 
invariants of $H$.  
Since $\bold 1_H-\phi^{n+1}(\bold 1_H)$ is the projection
$$
\bold 1_H-\phi^{n+1}(\bold 1_H)=
E_0+E_2+\dots+E_{2n},
$$ 
it follows that 
$$
\hat\phi(t) =\sum_{n=0}^\infty \dim(E_0+E_2+\dots+E_{2n})t^n,
$$
and therefore
$$
(1-t)\hat\phi(t) = \sum_{n=0}^\infty \dim E_{2n}t^n.  
\tag{9.10}
$$

Setting 
$$
\hat\sigma(t) = \sum_{p=0}^\infty \dim E_p t^p,
$$
we find that for $0<t<1$ 
$$
\align
\sum_{n=0}^\infty \dim E_{2n} t^n &=
1/2(\sum_{p=0}^\infty \dim E_p (\sqrt t)^p + 
\sum_{p=0}^\infty \dim E_p (-\sqrt t)^p) \\
&=1/2(\hat\sigma(\sqrt t)+\hat\sigma(-\sqrt t)),
\endalign
$$
and hence from (9.10) we have 
$$
\hat\phi(t) = \frac{\hat\sigma(\sqrt t)+\hat\sigma(-\sqrt t)}
{2(1-t)}, \qquad 0<t<1.  
\tag{9.11}
$$

The dimensions $\dim E_p$ were computed in \cite{1, Appendix A},
where it was shown that $\dim E_p = q_2(p)$, $q_2(x)$ being 
the polynomial defined by (3.3).  Thus
$$
\hat\sigma(t)=\sum q_2(n) t^n = \frac{1}{(1-t)^3};
$$
and finally (9.11) becomes
$$
\align
\hat\phi(t) = \frac{(1-\sqrt t)^{-3}+(1+\sqrt t)^{-3}}
{2(1-t)} = \frac{(1+\sqrt t)^3+(1-\sqrt t)^3}
{2(1-t)^4} = 
\frac{1+3t}{(1-t)^4}.  
\endalign
$$
The right side of the last equation can be rewritten 
$$
\hat\phi(t) = \frac{-3}{(1-t)^3}+\frac{4}{(1-t)^4},
$$
hence the coefficients $(c_0,c_1,\dots,c_6)$ of 
(7.1) are given by $(0,0,0,-3,4,0,0)$.  
One now reads off the numerical invariants listed
in Proposition 9.8.  
\qedd
\enddemo

\remark{Remark 9.12}
One can easily write down a class of related examples
(all having Euler characteristic zero) by considering 
powers $\sigma^N: \Cal B(H^2(\Bbb C^d))\to\Cal B(H^2(\Bbb C^d))$
of the completely positive map 
$\sigma(B)=S_1BS_1^*+\dots+S_dBS_d^*$
of the $d$-shift for 
arbitrary powers $N\geq 2$ and in arbitrary dimensions 
$d$.  The example (9.5) we have discussed is associated with 
the values $N=2$ and $d=3$.  To explain this briefly, 
notice that if one expands the expression
for $\sigma^2$ in dimension $d=3$ 
into a sum of the form 
$$
\sigma^2(B)=\sum_{k=1}^6T_kBT_k^*, 
\qquad B\in\Cal B(H^2(\Bbb C^3)), 
$$
then one finds that one set of candidates for 
the $6$ operators $T_1,\dots,T_6$ is given by 
$$
S_1^2,S_2^2,S_3^2, \sqrt2 S_1S_2,\sqrt2 S_1S_3,\sqrt2 S_2S_3.  
\tag{9.13}
$$
These operators are the most natural ones, 
but of course one can replace them with 
{\it certain} linear combinations to obtain other $6$-tuples 
which also serve to represent $\sigma^2$; any two 
such $6$-tuples must be related by a complex unitary 
$6\times 6$ matrix as in (8.9).  

Once one settles on a $6$-tuple such as (9.13), one 
finds that while it certainly defines a $6$-contraction, 
it is not an irreducible $6$-contraction 
because each of the operators $T_1,\dots,T_6$ leaves 
both the even subspace and its orthogonal complement
(the odd subspace) of 
$H^2(\Bbb C^3)$ invariant.  The 
example leading to (9.5) was obtained by restricting the 
$6$-tuple (9.13) 
to the irreducible even subspace.  

If one chooses other values for $N$ and the dimension
$d$, one finds that this method 
generates infinitely many higher dimensional examples
for which one can, in principle, explicitly 
calculate $\deg(H)$ and $\mu(H)$.  
\endremark

Finally, we compute nontrivial values of the curvature
invariant $K(H)$ for certain examples of pure 
rank-two graded Hilbert modules $H$.  Let 
$\phi$ be a homogeneous 
polynomial of degree $N=1,2,\dots$ 
in $A=\Bbb C[z_1,\dots,z_d]$ and let 
$M$ be the graph of its 
associated multiplication operator 
$$
M=\{(f,\phi\cdot f): f\in H^2\}\subseteq H^2\oplus H^2.  
$$
$M$ is a closed submodule of the free 
Hilbert module $F=H^2\oplus H^2$, 
and $H=F/M$ is a pure Hilbert 
module of rank $2$ whose minimal dilation 
$$
F\underset U\to\longrightarrow H\longrightarrow 0
\tag{9.14}
$$
is given by the natural projection $U$ of $F$ onto 
$H=F/M$.  

We make $H$ into a graded Hilbert module as follows.  
Let $\Gamma$ be the gauge group defined on $F=H^2\oplus H^2$
by
$$
\Gamma(\lambda)(f,g) = 
(\Gamma_0(\lambda)f,\lambda^{-N}\Gamma_0(\lambda)g), 
\qquad f,g\in H^2,
$$
where $\Gamma_0$ is the natural gauge group of $H^2$ defined
by 
$$\Gamma_0(\lambda)f(z_1,\dots,z_d)=
f(\lambda z_1,\dots,\lambda z_d).
$$
One verifies that 
$\Gamma (\lambda)M\subseteq M$, $\lambda\in\Bbb T$.  
Thus the action of $\Gamma$ can be promoted naturally
to the quotient $H=F/M$, and $H$ 
becomes a graded rank-two 
pure Hilbert module whose gauge group has spectrum 
$\{-N, -N+1,\dots\}$.   (9.14) becomes a graded 
dilation in that 
$U\Gamma(\lambda)=\Gamma(\lambda)U$ for 
all $\lambda\in\Bbb T$.  

\proclaim{Proposition 9.15}
For these examples we have $K(H)=\chi(H)=1$.  
\endproclaim
\demo{proof}
By Theorem B, $K(H)=\chi(H)$, and it suffices 
to show that $\chi(H)=1$.  

Let $H_n=\{\xi\in H: \Gamma(\lambda)\xi=\lambda^n\xi\}$, 
$n\in\Bbb Z$, be the spectral subspaces of $H$.  It is 
clear that $H_n=\{0\}$ if $n<-N$, and (9.14) implies 
that the algebraic submodule $M_H$ is given by
$M_H=U(A\oplus A)$.  Hence $M_H$ is the (algebraic) sum
$$
M_H=\sum_{n=-\infty}^\infty H_n.  
$$
Consider the proper filtration 
$M_1\subseteq M_2\subseteq\dots$ of $M_H$ defined by
$$
M_k=\sum_{n\leq k}H_n, \qquad k=1,2,\dots.  
$$
By the Corollary of Proposition 3.10 we have
$$
\chi(M_H) = d!\lim_{k\to\infty}\frac{\dim M_k}{k^d},
\tag{9.16}
$$
and thus we have to calculate the dimensions 
$$
\dim M_k = \dim (\sum_{n\leq k} H_n) =
\dim H_{-N}+\dots +\dim H_{k-1}+\dim H_k  
\tag{9.17}
$$
for $k=1,2,\dots$.  

In order to calculate the dimension of $H_n$ it is 
easier to realize $H$ as the orthogonal complement 
$M^\perp\subseteq F$, with canonical operators 
$T_1,\dots,T_d$ given by compressing the natural 
operators of $F=H^2\oplus H^2$ to $M^\perp$.  
Since $M$ is the graph of the multiplication operator
$M_\phi f=\phi\cdot f$, $f\in H^2$, $M^\perp$ is 
given by 
$$
M^\perp = \{(-M_\phi^* g,g): g\in H^2\}.  
$$
We compute 
$$
H_n=(M^\perp)_n = 
\{\xi\in M^\perp: \Gamma(\lambda)\xi=\lambda^n\xi,
\quad\lambda\in\Bbb T\}.  
$$
Since 
$\Gamma_0(\lambda)M_\phi^*\Gamma_0(\lambda)^{-1} =
(\Gamma_0(\lambda)M_\phi\Gamma_0(\lambda)^{-1})^*=
(\lambda^NM_\phi)^* = \lambda^{-N}M_\phi^*$, we have 
$$
\Gamma(\lambda)(-M_\phi^*g,g) =
(-\Gamma_0(\lambda)M_\phi^*g,\lambda^{-N}\Gamma_0(\lambda)g)
=(-\lambda^{-N}M_\phi^*\Gamma_0(\lambda)g,
\lambda^{-N}\Gamma_0(\lambda)g),  
$$
thus $\Gamma(\lambda)(-M_\phi^*g,g) = \lambda^n(-M_\phi^*g,g)$
iff $\Gamma_0(\lambda)g = \lambda^{n+N}g$, $\lambda\in \Bbb T$.  
For $n<-N$ there are no nonzero solutions of this equation, and 
for $n\geq -N$ the condition is satisfied iff $g$ is a homogeneous
polynomial of degree $n+N$.  

We conclude that $\dim H_n=0$ if $n<-N$ and 
$\dim H_n = \dim A_{n+N} = q_{d-1}(n+N)$ if $n\geq -N$.  Thus
for $k\geq -N$ we see from (9.17) that 
$$
\dim M_k = \sum_{n=-N}^k H_n = 
\sum_{n=-N}^k q_{d-1}(n+N).  
$$
The identity $q_{d-1}(x)=q_d(x)-q_d(x-1)$ of (3.2.2)
implies that the right 
side of the preceding formula telescopes to 
$q_d(k+N)-q_d(-1)=q_d(k+N)$.  Thus (9.16) 
implies that  
$$
\chi(H) = \chi(M_H)
= d!\lim_{k\to\infty}\frac{q_d(k+N)}{k^d} =
\lim_{k\to\infty}\frac{(k+N+1)\dots(k+N+d)}{k^d}=1,
$$
as asserted.
\qedd
\enddemo

\vfill
\pagebreak

\end
\Refs
%

\ref\no 1\by Arveson, W.\paper Subalgebras of $C^*$-algebras III:
Multivariable operator theory
\jour Acta Math\paperinfo to appear
\endref

\ref\no 2\bysame\paper The index of a quantum dynamical
semigroup\jour J. Funct. Anal. 
\vol 146, no. 2 \yr 1997 \pages 557--588
\endref

\ref
\no 3\bysame\paper Dynamical invariants for 
noncommutative flows\inbook Operator Algebras 
and Quantum Field Theory
\pages 476--514\bookinfo Doplicher, S.,
Longo, R., Roberts, J. and Zsido, L. editors
\publ International Press
\publaddr Accademia Nazionale dei Lincei, Roma
\yr 1997
\endref


\ref\no 4\by Athavale, A.\paper On the intertwining
of joint isometries
\jour Jour. Op. Th. \vol 23\yr 1990\pages 339--350
\endref



\ref\no 5\by Bunce, John W.\paper Models for $n$-tuples 
of noncommuting operators
\jour J. Funct. Anal. \yr 1984\vol 57\pages 21--30
\endref

\ref\no 6
\by Chern, S.-S.\paper A simple intrinsic proof 
of the Gauss-Bonnet formula for closed Riemannian
manifolds\jour Ann. Math. \vol 45\yr 1944
\pages 741--752
\endref


\ref\no 7
\by Cowen, M. and Douglas, R.\paper Complex geometry 
and operator theory\jour Acta Math. \vol 141\yr 1976
\pages 187--261
\endref


\ref\no 8\by Curto, R. and Vasilescu, F.-H.
\paper Automorphism invariance of the operator-valued
Poisson transform 
\jour Acta Sci. Math. (Szeged)\yr 1993 \vol57\pages 65--78
\endref

\ref\no 9
\by Do Carmo, M. P.\book Differential Geometry of 
Curves and Surfaces \publ Prentice Hall\yr 1976
\publaddr Englewood Cliffs, New Jersey
\endref

\ref\no 10
\by Douglas, R. \paper Hilbert modules for function
algebras \inbook Operator Theory: Advances and 
Applications\vol 17 \publ Birkhauser Verlag
\publaddr Basel\yr 1986\pages 125--139
\endref

\ref\no 11
\by Douglas, R. and Yan, K. \paper Hilbert-Samuel 
polynomials for Hilbert modules\jour Indiana Univ. 
Math. J.\vol 42, no. 3\yr 1993\pages 811--820
\endref

\ref\no 12\by Douglas, R. and Paulsen, V.
\book Hilbert Modules over Function Algebras
\bookinfo Pitman Research Notes in Mathematics
\vol 217 \publ Longman Scientific \& Technical
\publaddr Harlow, Essex, UK
\endref

\ref\no 13\by Drury, S.\paper A generalization of 
von Neumann's inequality to the complex ball
\jour Proc. AMS\yr 1978 \vol 68\pages 300--304
\endref

\ref\no 14\by Eisenbud, D.\book Commutative 
Algebra with a view toward algebraic geometry
\bookinfo Graduate Texts in Mathematics \publ
Springer-Verlag \vol 150\yr 1994
\endref


\ref\no 15\by Frazho, Arthur E.\paper Models for 
noncommuting operators
\jour J. Funct. Anal. \yr 1982 \vol 48\pages 1--11
\endref

\ref\no 16\by Hilbert, D. \yr 1890
\paper \"Uber die Theorie von algebraischen Formen
\jour Math. Ann. \vol 36\pages 473--534
\endref

\ref\no 17\bysame \yr 1893
\paper \"Uber die vollen Invariantensysteme
\jour Math. Ann. \vol 42\pages 313--373
\endref

\ref\no 18
\book Commutative Rings\by Kaplansky, I.
\publ Allyn and Bacon\publaddr Boston
\yr1970
\endref


\ref\no 19\by McConnell, J. C. and Robson, J. C.
\book Noncommutative Noetherian Rings
\publ Wiley\yr 1987\publaddr New York
\endref


\ref\no 20\by M\"uller, V. and Vasilescu, F.-H.
\paper Standard models for some commuting multioperators
\jour Proc. AMS\yr 1993 \vol 117\pages 979--989
\endref


\ref\no 21\by Pisier, G.\book Similarity problems
and completely bounded maps
\publ Springer Verlag Lecture Notes in Mathematics\vol 
1618\yr 1995
\endref


\ref\no 22\by Popescu, G.\paper von Neumann inequality
for $(\Cal B(H)^n)_1$
\jour Math. Scand.\yr 1991 \vol 68\pages 292--304
\endref

\ref\no 23\bysame\paper On intertwining dilations 
for sequences of noncommuting operators
\jour J. Math. Anal. Appl.\yr 1992 \vol 167\pages 382--402
\endref


\ref\no 24\bysame\paper Functional calculus for 
noncommuting operators
\jour Mich. Math. J.\yr 1995 \vol 42\pages 345--356
\endref

\ref\no 25\bysame\paper Multi-analytic operators on 
Fock space \jour Math. Ann.\yr 1995 \vol 303\pages 31--46
\endref


\ref\no 26\by Rudin, W.\book Function theory in the unit ball of 
$\Bbb C^n$\publ Springer Verlag\yr 1980
\endref

\ref\no 27\by Segal, I. E. \paper Tensor algebras over 
Hilbert spaces\jour Trans. Amer. Math. Soc. 
\yr 1956\pages 106--134
\endref

\ref\no 28\bysame\paper Tensor algebras over Hilbert 
spaces II\jour Annals Math.\vol 63, no. 1 \yr 1956\pages 160--175
\endref

\ref\no 29\by Serre, J.-P. \book Alg\`ebre, Locale
Multiplicit\'es \bookinfo
Lecture notes in Mathematics\publ Springer-Verlag
 \vol 11\yr 1965
\endref

\ref\no 30
\by Taylor, J. L.\paper A joint spectrum for several 
commuting operators\jour Jour. Funct. Anal. \vol 6
\yr 1970 \pages 172--191
\endref

\ref\no 31\by Vasilescu, F.-H.\paper An operator-valued
Poisson kernel 
\jour Jour. Funct. Anal.\vol 110\yr 1992\pages 47--72
\endref


\ref\no 32\bysame\paper Operator-valued 
Poisson kernels and standard models in several variables
\inbook Algebraic methods in operator theory\ed Ra\'ul 
Curto and Palle Jorgensen\yr 1994 \publ Birkh\"auser 
\pages 37--46
\endref

\endRefs
